\RequirePackage{fix-cm}
\documentclass[twocolumn]{svjour3}          
\smartqed  
\usepackage{subfigure}
\usepackage{booktabs}
\usepackage{caption}

\newlength{\tab}
\usepackage{amsfonts}

\usepackage{amsmath, amscd, graphicx,amssymb}

\DeclareMathAlphabet\gothic{U}{euf}{m}{n}

\newcommand*{\Scale}[2][4]{\scalebox{#1}{$#2$}}

\def\th{\theta}

\newcommand{\X}{X}
\newcommand{\Id}{\operatorname{Id}\nolimits}
\newcommand{\spann}{\operatorname{span}\nolimits}

\renewcommand{\Re}{\operatorname{Re}\nolimits}
 \newcommand{\sgn}{\operatorname{sgn}}
\newcommand{\ds}{\displaystyle}
\newcommand{\desda}{\Leftrightarrow}

\renewcommand{\i}{\mathrm{i}}
\newcommand{\eexp}[1]{\operatorname{e}^{#1}}

\newcommand{\tcusp}{t_\mathrm{cusp}}

\newcommand{\Arg}{\operatorname{arg}}
\newcommand{\const}{\operatorname{const}}

\newcommand{\sth}{\operatorname{s\th}}
\newcommand{\LL}{\operatorname{L}}
\newcommand{\smax}{s_\mathrm{max}}
\newcommand{\smaxE}{s^{E}_{max}}
\newcommand{\smaxL}{s^{L}_{max}}
\newcommand{\smaxH}{s^{H}_{max}}
\newcommand{\Pcurve}{\bf P_{curve}}
\newcommand{\Pmec}{\bf P_{mec}}
\newcommand{\PMEC}{\bf P_{mec}}
\newcommand{\PmecOne}{\bf P_{mec}^1}

\def\R{{\mathbb R}}
\def\C{{\mathbb C}}
\def\N{{\mathbb N}}
\def\Z{{\mathbb Z}}

\newcommand{\ul}{\mathbf}

\newcommand{\SO}{\operatorname{SO(3)}\nolimits}

\newcommand{\SE}{\operatorname{SE(2)}\nolimits}
\newcommand{\Htri}{\operatorname{H(3)}\nolimits}

\newcommand{\so}{\operatorname{so(3)}\nolimits}

\renewcommand{\sl}{\operatorname{sl(2)}\nolimits}

\renewcommand{\S}{S}
\renewcommand{\phi}{\varphi}

\newcommand{\M}{\operatorname{M}}

\newenvironment{sketchproof}{\textit{Proof of Theorem 4.1:}}{\hfill$\square$\par}

\DeclareMathOperator{\Exp}{Exp}
\DeclareMathOperator{\WF}{WF}

\begin{document}

\title{Tracking of Lines in Spherical Images via Sub-Riemannian Geodesics in $\SO$}

\author{A.~Mashtakov \and R.~Duits \and  Yu.~Sachkov \and  E.J.~Bekkers \and I.~Beschastnyi} 

\authorrunning{Mashtakov, Duits, Sachkov, Bekkers, Beschastnyi}

\institute{
A.~Mashtakov \at  CPRC,
              Program Systems Institute of RAS, \\
              Tel.: +7-965-216-81-04,
              \email{alexey.mashtakov@gmail.com}  \\         
                  \and
R.~Duits \at  CASA \& BMIA,
              Eindhoven University of Technology, \\
              Tel.: +31-40-2472859,
              \email{R.Duits@tue.nl}  \\         
                  \and
Yu.~Sachkov \at  RUDN University, \\
              Tel.: 007-4852-695-228,        %
              \email{sachkov@sys.botik.ru} \\
\and
E.J.~Bekkers \at BMIA \& CASA,
              Eindhoven University of Technology, \\
              Tel.: +31-40-2473037,        %
              \email{e.j.bekkers@tue.nl} \\
\and
I.~Beschastnyi \at MAMA,
              International School for Advanced Studies, \\
              Tel.:+39-40-3787442,        %
              \email{i.beschastnyi@gmail.com}
}

\maketitle

\begin{abstract}
In order to detect salient lines in spherical images, we consider the problem of minimizing the functional $\int \limits_0^l \gothic{C}(\gamma(s)) \sqrt{\xi^2 + k_g^2(s)} \, {\rm d}s$ for a curve $\gamma$ on a sphere with fixed boundary points and directions. The total length $l$ is free, $s$ denotes the spherical arclength, and $k_g$ denotes the geodesic curvature of~$\gamma$. Here the smooth external cost $\gothic{C}\geq \delta>0$ is obtained from spherical data. We lift this problem to the sub-Riemannian (SR) problem in Lie group $\SO$ and show that the spherical projection of certain SR geodesics provides a solution to our curve optimization problem. In fact, this holds only for the geodesics whose spherical projection does not exhibit a cusp. The problem is a spherical extension of a well-known contour perception model, where we extend the model by Boscain and Rossi to the general case $\xi > 0$, $\gothic{C} \neq 1$. For $\gothic{C}=1$, we derive SR geodesics and evaluate the first cusp time. We show that these curves have a simpler expression when they are parameterized by spherical arclength rather than by sub-Riemannian arclength. For case $\gothic{C} \neq 1$ (data-driven SR geodesics), we solve via a SR Fast Marching method. Finally, we show an experiment of vessel tracking in a spherical image of the retina and study the effect of including the spherical geometry in analysis of vessels curvature.
\keywords{Sub-Riemannian geodesics \and Geometric control \and Spherical image \and Lie group SO(3) \and Vessel tracking}
\end{abstract}
\section{Introduction}\label{sec:Intro}

In computer vision, it is common to extract salient curves in flat images via data-driven minimal paths or 
geodesics~\cite{PeyreCohen2010,Kimmel,mirebeau,sethian,osher}.
The minimizing geodesic is
defined as the curve that minimizes the length functional, which is typically weighted by a
cost function with high values at image locations with high curve saliency.

Another set of geodesic methods, partially inspired by the psychology of vision, was developed
in~\cite{Petitot_B,Citti_Sarti_B}. In these articles, sub-Riemannian (SR) geodesics in respectively the Heisenberg
$\Htri$ and the Euclidean motion group $\SE$ are proposed as a model for contour perception and contour integration.

The combination of such contour perception models with data-driven  geodesic methods has been presented in~\cite{SIAM}. There, a
computational framework for tracking of lines via globally optimal data-driven sub-Riemannian geodesics on the Euclidean motion group $\SE$ has been presented with comparisons to exact solutions~\cite{cut_sre2}.

\begin{figure}[ht]
\centering
\includegraphics[width=0.9\hsize]{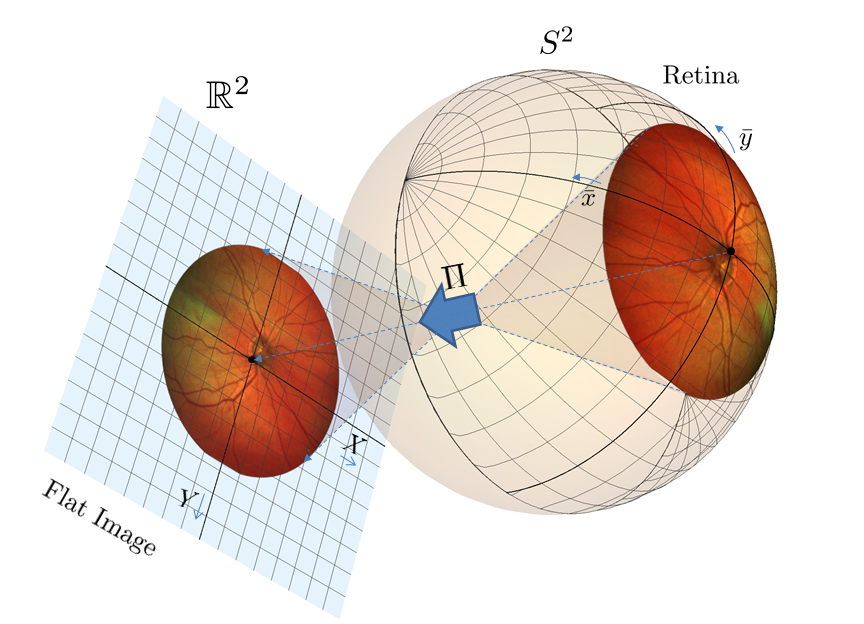}
\caption{ Photography of the retina. A part of the retina is projected onto the image plane. The camera coordinates are denoted by $(X,Y)$, and object coordinates are denoted by $(\bar{x},\bar{y})$.}
\label{fig:RetinaSphere}
\end{figure}

In this work, we extend this framework for tracking of lines in \emph{spherical} images (e.g., images of the retina, cf. Fig.~\ref{fig:RetinaSphere}). This requires a sub-Riemannian manifold structure in a different Lie group, namely the group $\SO$ (consisting of 3D rotations) acting transitively on the 2-sphere $S^2$.

\begin{figure}[ht]
\centering
\includegraphics[width=0.50\hsize]{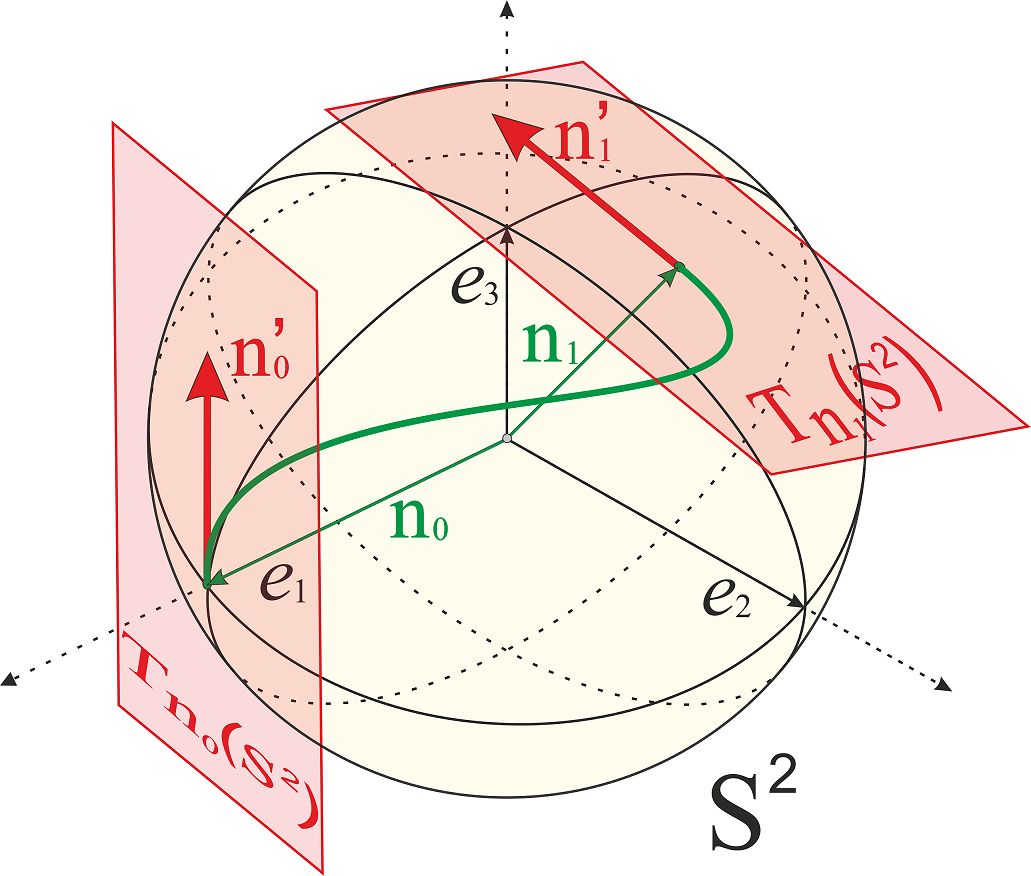}
~
\includegraphics[width=0.45\hsize]{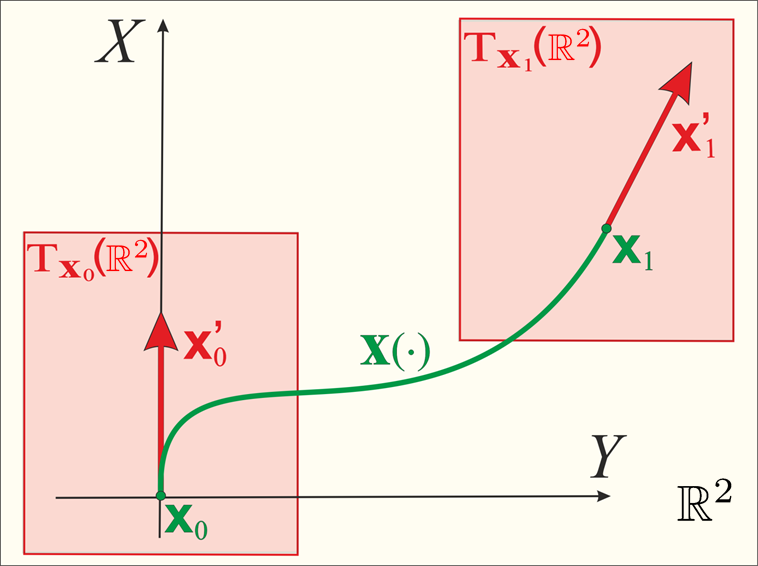}
\caption{
\textbf{Left:} Problem $\Pcurve(S^2)$: for given boundary conditions on a 2D sphere (both positions and velocities), we aim to find a curve  minimizing the functional compromising length and geodesic curvature. In the optimization functional we also include an external cost induced by spherical image data.
\textbf{Right: } Problem $\Pcurve(\R^2)$~\cite{Boscain2014,DuitsJMIV2014}: for given boundary conditions on a 2D plane, to find a curve  minimizing the compromise between length and curvature. The external cost factor is added for adaptation to flat image data (see ~\cite{SIAM}).
\label{fig:PcurveStat}}
\end{figure}

Here we study the problem $\Pcurve(S^2)$ of finding a smooth curve $\ul{n}(\cdot)$ on a 
 unit sphere $S^2$ that satisfies gi\-ven boundary conditions (both positions and velocities)
\[\ul{n}(0) = \ul{n}_0, \, \, \ul{n}(l) = \ul{n}_{1}, \quad \ul{n}'(0) = \ul{n}_0', \, \, \ul{n}'(l) = \ul{n}_1' \]
 and minimizes the functional
\[\int_0^l \gothic{C}(\ul{n}(s))\sqrt{\xi^2 + k_g^2(s)} \, {\rm d} s,\]
where $k_g(\cdot)$ denotes the geodesic curvature of $\ul{n}(\cdot)$, $s$~denotes the spherical arclength, and total length $l$ is free, see Fig.~\ref{fig:PcurveStat}. Furthermore, we include a curve-stiffness parameter $\xi>0$. In the optimization functional, we also include an external cost factor $\gothic{C}: S^2 \to \R^+$ adding for adaptation to a given spherical image data. We state this problem $\Pcurve(S^2)$ more explicitly in Sect.~\ref{subsec:Pcurvestat}.

The problem  $\Pcurve(S^2)$ is a spherical analog of a well-known problem $\Pcurve(\R^2)$~\cite{Boscain2014,DuitsJMIV2014} (cf. Fig.~\ref{fig:PcurveStat}) of finding a smooth curve $\ul{x}(\cdot)$ on a plane $\R^2$ that satisfies given boundary conditions
\begin{equation*}
\begin{array}{c c}
\ul{x}(0) = \ul{x}_0 = (X_0,Y_0), \, \, & \, \, \ul{x}'(0) = \ul{x}_0'= (\cos\Theta_0, \sin\Theta_0), \\
\ul{x}(l) = \ul{x}_{1}= (X_1,Y_1), \, \, & \, \, \ul{x}'(l) = \ul{x}_1'= (\cos\Theta_1, \sin\Theta_1),
\end{array}
\end{equation*}
 and minimizes the functional
\begin{equation*}
 \int_0^l \gothic{c}(\ul{x}(s))\sqrt{\xi^2 + k^2(s)} \, {\rm d} s,
\end{equation*}
where $k(s)$ denotes the curvature, $l$ denotes the total length, and curve stiffness is regulated by $\xi>0$. The external cost factor $\gothic{c}: \R^2 \to \R^+$ is added for adaptation to a given flat image data (see ~\cite{SIAM}).

There are three motivations for our study.
The first motivation comes from a medical image analysis application, where automatic extraction of the vascular retinal tree on images is helpful for early detection of many diseases, such as diabetic retinopathy, glaucoma, atherosclerosis, and others (see, e.g.,~\cite{Ikram,Remco_Eric_B}). Optical retinal images are mostly acquired by flat cameras, and as a result, distortion appears. Such distortion could lead to questionable (distorted) geometrical features (vessel curvature, thickness, etc.) that are used as biomarkers~\cite{Ikram,Kalitzeos,Erik_Curv} for different diseases. We will show that the distortion that appears near the boundary of a flat image can play a significant role in the quantitative analysis of the vascular structure and its curvature. 
As the retina itself is not flat, we should include the spherical geometry both in image processing and in image representation in order to avoid distortion. Such distortion comes from the central projection of the physical retinal surface to the image plane. It is illustrated in Fig.~\ref{fig:RetinaSphere}, where $(X,Y)$ denote Cartesian camera coordinates and $(\bar{x},\bar{y})$ denote spherical object coordinates. This motivates us to study data-driven SR geodesics on the  rigid body rotation group $\SO$ and their spherical projections;
likewise, it was done (see~\cite{SSVM,SIAM}) for flat images and lifting to the roto-translation group $\SE$.

The second motivation comes from models of the visual system of mammals. As mentioned 
by Boscain and Rossi~\cite{Boscain_Rossi_B}, the problem $\Pcurve(S^2)$ can be considered as a spherical extension of a (flat) cortical-based model $\Pcurve(\R^2)$ for perceptual completion, proposed by Citti and Sarti~\cite{Citti_Sarti_B}, and Petitot~\cite{Petitot_B}. Such a spherical extension is again motivated by the fact that the retina is not flat. By the same argument, cuspless SR geodesics in $\SO$ could provide a model of association fields in the psychology of vision (see~\cite{DuitsJMIV2014}). Note that nonuniform distribution of photoreceptors on the retina can be included in our model by putting an external cost $\gothic{C}$ on a hemisphere. 

The third motivation for this study is that in geometric control theory optimal synthesis for the SR problem in $\SO$ has not been achieved in the general case (not even for the case of uniform cost $\gothic{C}=1$),  despite many strong efforts in this direction~\cite{Berestovskiy1,Berestovskiy2,Schcerbakova,Boscain_Rossi_B,Bonnard,Ivan,Markina1,Markina2}. In this paper, we will not provide optimal synthesis analytically, but instead we do provide a HJB theory for computing globally optimal (data-driven) geodesics. In previous works~\cite{SSVM,SIAM}, we achieved this for SR geodesics in $\SE$, used for tracking of blood vessels  in flat 2D images.

In view of these three motivations, we lift the problem $\Pcurve$ on the set $\S^2$ to a SR problem on the group $\SO$, which we will call $\Pmec$,  and we provide explicit formulas for SR geodesics. 
This allows us to describe the set of endpoints in $\SO$ reachable by geodesics whose spherical projections do not contain cusps. Furthermore, we present a Hamilton--Jacobi--Bellman PDE theory that allows us to numerically compute the SR distance map, from which a steepest descent backtracking (via the Pontryagin maximum principle) provides only the \emph{globally optimal} geodesics for general external cost and general $\xi>0$. We verify our numerical solution, by comparison with exact geodesics in the case $\gothic{C}=1$.  Finally, we use these results in a vessel tracking algorithm in spherical images of the retina, without central projection distortion.

The main contributions of this article are:
\begin{itemize}
\item[$\bullet$] New formulas for the geodesics of $\Pcurve(S^2)$ for the uniform cost case $\gothic{C}=1$.
\item[$\bullet$] Analysis and parametrization of cusp points arising in spherical projections of geodesics of $\Pmec(\SO)$.
\item[$\bullet$] A new HJB-PDE theory for the numeric computation of globally optimal geodesics of $\Pcurve(S^2)$ and $\Pmec(\SO)$ for the general external cost case, with verification for $\gothic{C}=1$.
\item[$\bullet$] Tracking of lines in spherical image data (e.g., retinal images) without central projection distortion, with comparison to tracking of lines in flat image data.
\end{itemize}

\subsection{Vessel Tracking in Spherical Images of the Retina}\label{subsec:Motivation1}
In this subsection, we first describe the mapping from object coordinates on the retina to camera coordinates, and then we discuss the relevance of considering spherical images of the retina rather than flat images, which are commonly used in the medical application (see Fig.~\ref{fig:RetinaSphere}). We show that distortion appears inevitably on flat images, with a significant relative error (up to 7\%) in length measures. Even larger relative errors (over 20\%) appear in the application of differential operators (used for vessel detection).

\begin{figure}[ht]
\centering
\includegraphics[width=0.9\hsize]{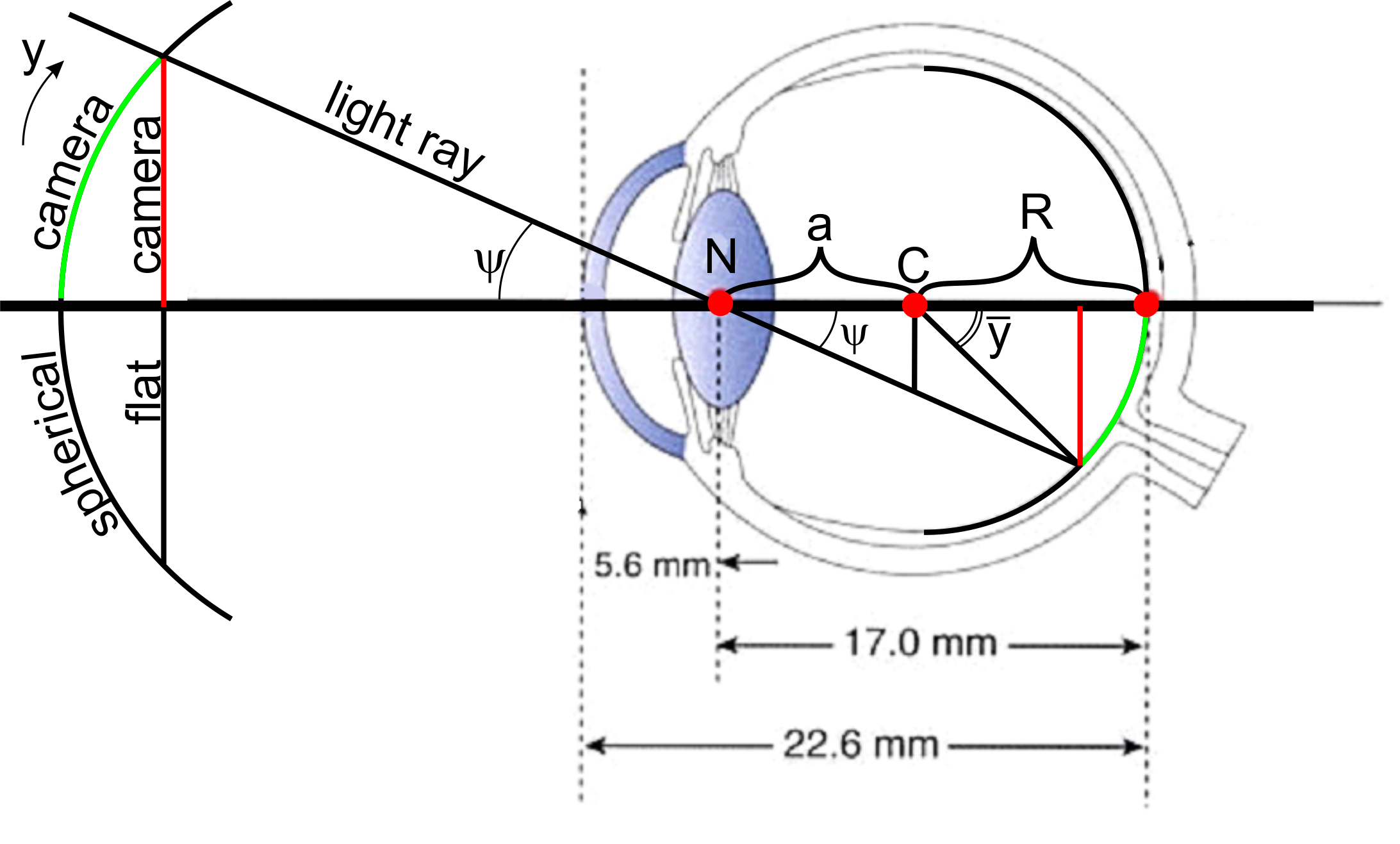}\\
\includegraphics[width=0.5\hsize]{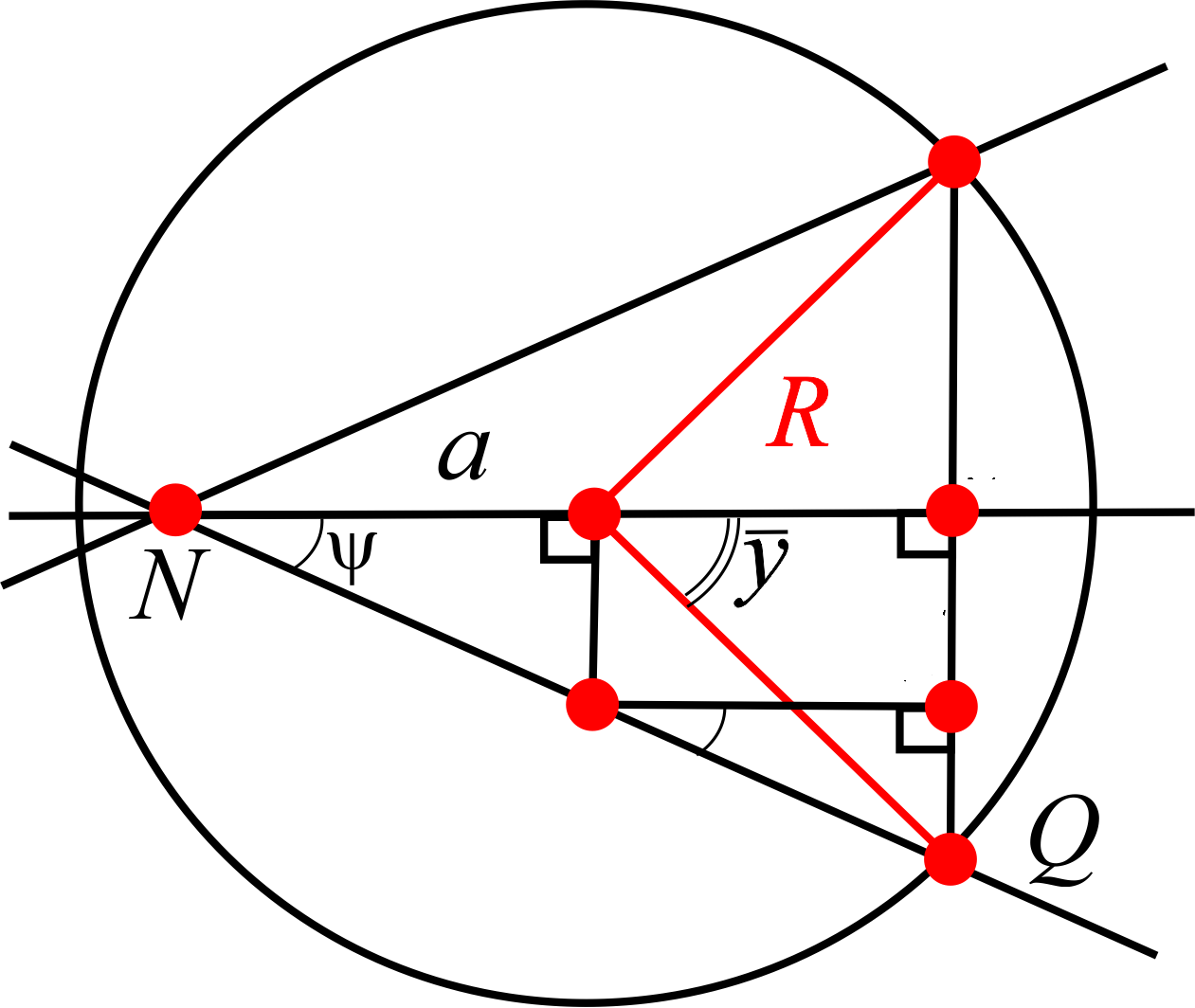}
\caption{
\textbf{Top:} Schematic eye~\cite{ClinOpt} and central projection of images onto the retina. Here $R \approx 10.5$mm.\\
\textbf{Bottom:}
Schematic eye, enlarged to support Eq.~(\ref{eq:ymax})}
\label{fig:Eye2}
\end{figure}

We base our computations on the reduced schematic eye model (see Fig.~\ref{fig:Eye2}), which is commonly used in clinical ophthalmology (see, e.g.,~\cite{ClinOpt}). Let $R$ be the radius of an eyeball, $a$ be the distance from the nodal point $N$ to the center $C$, and $\psi$ be the angle between visual axis and a light ray passing through $N$. Here we consider a simplified model, where the optical axis (the best approximation of a line passing through the optical center of the cornea, lens, and fovea) coincides with the visual axis (the line connecting fixation point and the fovea).\footnote{There is small difference between these two axes (c.f. Fig.~33~\cite{ClinOpt}) which we neglect in our basic model} The average radius of a human eye is $R \approx 10.5$ mm, and the maximum distance between nodal point $N$ and the central point $C$ is $a_\mathrm{max} = 17\, \text{mm} - 10.5\, \text{mm} = 6.5\, \text{mm}$.

Now we switch to mathematical object coordinates of the retina where we use the eyeball radius to express lengths, i.e., we have $R=1$ and $a_\mathrm{max}= \frac{6.5\, \text {mm}}{10.5 \,\text {mm}}\cdot R = \frac{13}{21}$ in dimensionless coordinates.

In order to compute the maximum absolute angle ${\bar{y}}_\mathrm{max}$ let us express the angle $|\bar{y}|$ with respect to center of the eyeball (see right Fig.~\ref{fig:Eye2}). Expressing the squared length of segment $NQ$  yields
 \begin{equation*}\label{eq:ymaxPrelim}
 (a+R \cos |\bar{y}|)^2 \cos^{-2}\psi  = (a+R\cos |\bar{y}|)^2 + R^2 \sin^2 |\bar{y}|.
 \end{equation*}
 Solving this equation with respect to $\cos|\bar{y}|$, we obtain unique nonnegative solution:
\begin{equation}
\label{eq:ymax}
\Scale[0.94]{|\bar{y}(a,R,\psi)|=\arccos \left(\cos\psi\ \sqrt{1-\frac{a^2 \sin^2\psi}{R^2}}-\frac{a \sin ^2\psi}{R}\right).}
\end{equation}

A standard fundus camera used for producing of the retinal images has the angular range $\psi \in [-\frac{\pi}{8}, \frac{\pi}{8}]$. Thus, substitution $R = 1$, $a = \frac{13}{21}$, and $\psi = \frac{\pi}{8}$ in Eq.~(\ref{eq:ymax})
gives the maximum angle
\begin{equation}
\label{eq:ymaxfondus}
{\bar{y}}_\mathrm{max} \approx 0.63 \text{rad} \approx 36^{\circ}.
\end{equation}

\begin{figure}[ht]
\centering
\includegraphics[width=0.8\hsize]{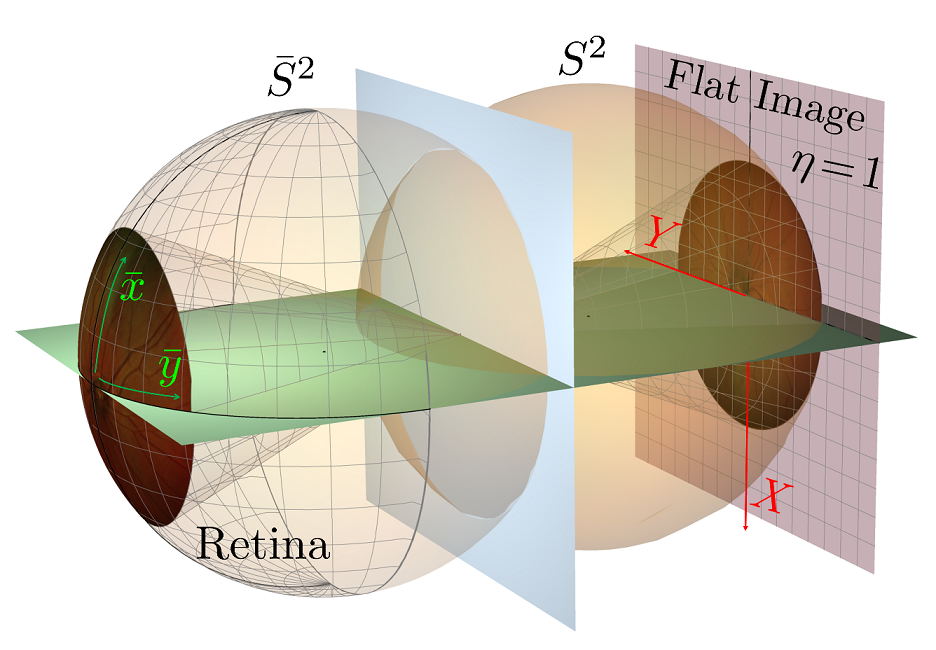}\\
\includegraphics[width=0.9\hsize]{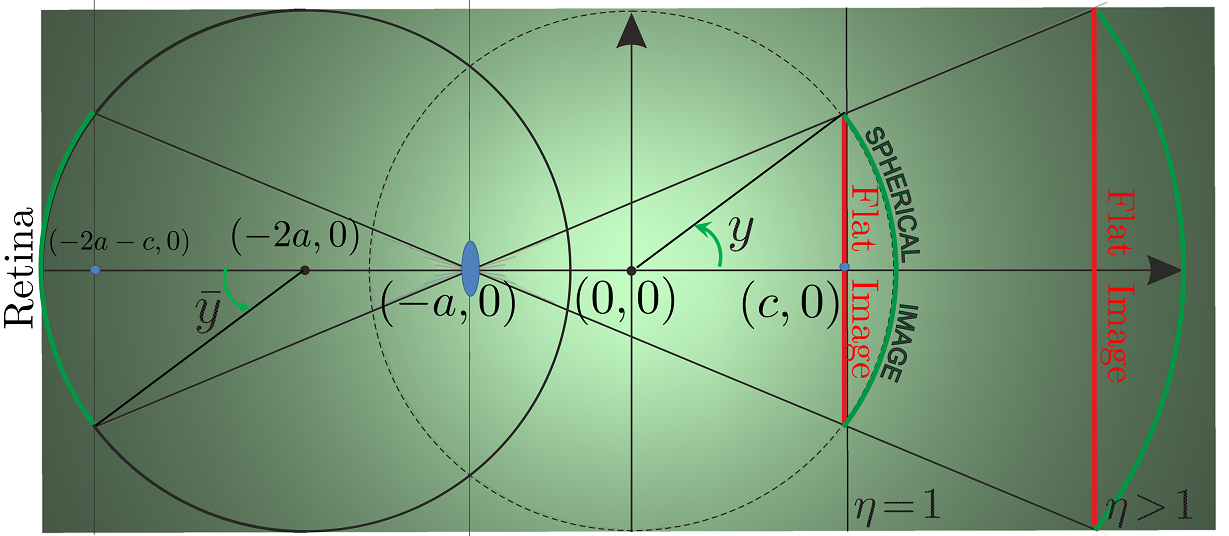}
\caption{Spherical object coordinates $(\bar{x},\bar{y})$ on a retina, cartesian camera coordinates $(X,Y)$ on a flat image of the retina, and spherical camera coordinates $(x,y)$ on a spherical image of the retina. 
}
\label{fig:Twospheres}
\end{figure}

We rely on the following parameterization of the image sphere $S^2$ and the retinal sphere $\bar{S}^2$ (see Fig.~~\!\ref{fig:Twospheres}):
\begin{equation*}\label{eq:S2param}
\begin{array}{l}
S^2 \ni \ul{n}(x,y) = \left(\cos x \cos y, \, \cos x \sin y, \, \sin x \right)^T,\\
\bar{S}^2 \ni \left(-2 a, 0,0\right)^T-\left(\cos \bar{x} \cos \bar{y}, \, \cos \bar{x} \sin \bar{y}, \, \sin \bar{x} \right)^T,\\
%
\end{array}
\end{equation*}
with $x, \bar{x} \in [-\frac{\pi}{2}, \frac{\pi}{2}]$ and $y, \bar{y} \in  [-\pi, \pi]$.

Next we present the explicit relation between the object coordinates $(\bar{x},\bar{y})$ \--- spherical coordinates on a unit sphere $\bar{S}^2$ representing the surface of the eyeball; flat photo coordinates $(X,Y)$ --- Cartesian coordinates on the image plane and the spherical coordinates $(x,y)$ on the image sphere, see Figs.~\!\ref{fig:RetinaSphere} and~\!\ref{fig:Twospheres}.

To take into account the distance from the eyeball to the camera in our model, we introduce a parameter $\eta>0$.
In Fig.~\!\ref{fig:Twospheres}, by setting $\eta=1$ we fix the distance equal to $(a + c)$ radiuses of the eyeball.
This corresponds to the case when the image sphere $S^2$ is obtained by reflection of the physical retinal sphere $\bar{S}^2$ through the point $(-a,0,0)^T \in \R^3$. In this case, we have
\begin{equation}\label{eq:xytobarxbary}
(x,y) = (\bar{x},\bar{y}),
\end{equation}
and we will always rely on this identification in the sequel. The general case $\eta>0$ can be taken into account by congruency and scaling $X \mapsto \eta R X$, $Y \mapsto \eta R Y$.

The central projection $\Pi$ (cf. Fig.~\ref{fig:RetinaSphere}) from $(x,y)$ to $(X,Y)$ including the scaling factor $\eta>0$ (with physical dimension length in units of $R$) is given by
\begin{equation} \label{eq:xytoXY}
X = \frac{(a+c)\sin x}{a+ \cos{x} \cos{y}} \eta, \quad Y = \frac{(a+c)\cos x \sin y}{a+\cos x \cos y} \eta.
\end{equation}

The inverse mapping $\Pi^{-1}$ from $(X,Y)$ to $(x,y)$ for $\eta =1$ is given by
\begin{equation} \label{eq:XYtoxy}
\begin{array}{l}
x = \arcsin(X \bar{p}(X,Y)), \\
y = \arg(p_1(X,Y) + \i \; Y \bar{p}(X,Y)),
\end{array}
\end{equation}
where
\begin{eqnarray*}
&&\bar{p}(X,Y) = \frac{a(a+c) +\Xi_{a,c}(X,Y)}{(X^2+Y^2) + (a+c)^2}, \\
&&p_1(X,Y) = \frac{(a+c)\Xi_{a,c}(X,Y)-a(X^2+Y^2)}{(X^2+Y^2)+(a+c)^2},
\end{eqnarray*}
with $\Xi_{a,c}(X,Y) = \sqrt{(X^2+Y^2)(1-a^2) + (a+c)^2}$.

In these formulas, we need to substitute $a=a_\mathrm{max} = \frac{13}{21}$ and $c=\frac{4}{5}<R=1$, depicted in Fig.~\ref{fig:Twospheres}, where we work in dimensionless coordinates.
\subsubsection*{Quantification of Local and Global Deformation}
The local deformation from spherical photo coordinates $(x,y)$ to planar photo coordinates $(X,Y)$ is now given by the Jacobian
\begin{equation*}
\label{eq:localdef}
\begin{array}{l l}
J(x,y)  & =
\det \left(\begin{array}{c c}
\frac{\partial{X}}{\partial x}(x,y) & \frac{\partial{X}}{\partial y}(x,y) \\
\frac{\partial{Y}}{\partial x}(x,y) & \frac{\partial{Y}}{\partial y}(x,y)
\end{array} \right)  \\
 & = \displaystyle \frac{(a + c)^2 \cos x (1 + a \cos x \cos y)}{(a + \cos x \cos y)^3} \eta^2.
\end{array}
\end{equation*}
In mathematical analysis, we can set $\eta =1$; however, in experiments $\eta$ is to be taken into consideration. Note that for $\eta=1$ we have
\[0.77 \approx J(0,0) \leq J(x,y) \leq  J(y_\mathrm{max},y_\mathrm{max}) \approx 1.1,\]
showing that local deformation plays a considerable role and varies in $(x,y)$.

Next we consider the global distortion along the line $x=0$. It is defined as $\mathcal{GD}(y) = \frac{|y-Y(0,y)|}{|y|}$, and it has a maximum when $y = y_\mathrm{max}$. We have
$$
0\leq \mathcal{GD}(y) \leq \mathcal{GD}(y_\mathrm{max}) \approx 0.07,$$
and we see the distortion up to 7\% along the line $x=0$. The same holds along the line $y=0$.

We conclude that it makes a considerable difference to study $\Pcurve(S^2)$ or $\Pcurve(\R^2)$ in the retinal imaging application. In the sequel, we will write $\Pcurve$ instead of $\Pcurve(S^2)$ as we will always be concerned with the case where the base manifold equals $S^2$.
\section{Problem $\Pcurve$ on $\S^2$ and $\PMEC$ in $\SO$ }\label{sec:PmecPcurveAndConnect} 
In this section, we first provide preliminaries on group theoretical tools and notation, and then we show that the spherical projection of certain SR geodesics $\gamma(\cdot)$ on the Lie group $\SO$ provides solution curves to the problem $\Pcurve$ on $S^2$, which we formulate next.
\subsection{Mathematical Foundation and Notations}\label{subsec:mathnot}
The Lie group $\SO$ is the group of all rotations about the origin in $\R^3$. We shall denote a counterclockwise rotation around axis $\ul{a} \in \S^{2}$ with angle $\phi$ via $R_{\ul{a},\phi}$, in particular for rotations around standard axes
$$\ul{e}_{1}=(1,0,0)^{T}\!\!,\quad \ul{e}_{2}=(0,1,0)^T\!\!,\quad \ul{e}_{3}=(0,0,1)^T.$$
We use representation of $\SO$ by $3 \times 3$ matrices
\begin{equation*}
R(x,y,\th) = R_{\ul{e}_{3},y} R_{\ul{e}_{2},-x} R_{\ul{e}_1,\theta} =
\end{equation*}
\begin{eqnarray}
\label{matexp}
\Scale[1]{\left(\begin{array}{ccc}
cx\, cy \ & \  -sx \, cy \, s\th - sy \, c\th \ & \ sy \, s\th - sx \, cy \, c\th  \\
cx\, sy \ & \ cy \, c\th - sx \, sy \, s\th \ & \  -cy \, s\th - sx \, sy \, c\th \\
sx\, \ & \ cx\, s\th \ & \  cx \, c\th
\end{array}\right),}
\end{eqnarray}
where we denote $cx = \cos x$, $cy = \cos y$, $c\th= \cos \theta$, $sx = \sin x$, $sy = \sin y$, $\sth= \sin \theta$, and where
\begin{equation}\label{eq:xythdom}
(x,y,\theta) \in \left[-\frac{\pi}{2}, \frac{\pi}{2}\right] \times \R/\{2 \pi \Z\} \times \R/\{2 \pi \Z\}.
\end{equation}
The Lie group $\SO$ defines an associated Lie algebra $$\so = T_{\Id}(\SO) = \spann(A_1, A_2, A_3),$$
$$\Scale[0.9]{A_1 = \left(\begin{array}{ccc}
0 & 0 & 0 \\
0 & 0 & -1 \\
0 & 1 & 0
\end{array}\right), \quad
A_2 = \left(\begin{array}{ccc}
0 & 0 & 1 \\
0 & 0 & 0 \\
-1 & 0 & 0
\end{array}\right),\quad
A_3 = \left(\begin{array}{ccc}
0 & -1 & 0 \\
1 & 0 & 0 \\
0 & 0 & 0
\end{array}\right),}$$
where $T_{\Id}(\SO)$ denotes the tangent space at the unity element $e \in \SO$, represented by identity matrix $\Id$, which corresponds to the origin in the parameter space
$$e \, \, \sim \, \,  \Id \, \, \sim \, \, \{(x,y,\th) = (0,0,0)\}.$$
\begin{remark}\label{remark:AlgebraGen}
The formula for $A_i$ follows by
$$\Scale[0.94]{A_1 = \frac{\partial R(x,y,\th)}{\partial \th}\big|_{\Id}, \,\,\,   A_2 = -\frac{\partial R(x,y,\th)}{\partial x}|_{\Id}, \,\,\, A_3 = \frac{\partial R(x,y,\th)}{\partial y}|_{\Id}.}$$
\end{remark}

The nonzero Lie brackets are given by
\begin{equation}\label{eq:comtab}
[A_1, A_2] = A_3, \quad [A_1, A_3] = -A_2, \quad [A_2, A_3] = A_1.
\end{equation}
There is a natural isomorphism between  $\so$ and the Lie algebra $\LL$ of left-invariant vector fields in $\SO$, where commutators of the vector fields in $\LL$ correspond to the matrix commutators in $\so$
\begin{equation}\label{eq:commutrelat}
[R A, R B] = R [A,B], \quad A,B \in \so, \, R \in \SO.
\end{equation}
We express $\LL$ in matrix form as
\begin{eqnarray}\label{eq:vfXi}
\Scale[0.9]{
\LL = \spann(\X_1, \X_2, \X_3), \quad
\begin{cases} \X_1(x, y, \th) = -R(x, y, \th) A_2,\\  \X_2(x, y, \th) = R(x, y, \th) A_1 , \\  \X_3(x, y, \th) =  R(x, y, \th) A_3.\end{cases}
}
\end{eqnarray}
\begin{remark}
Note that at the unity one has
$$\Scale[0.94]{X_1\big|_{\Id} = -A_2, \quad X_2\big|_{\Id} = A_1, \quad X_3\big|_{\Id} = A_3.}$$
The formulas for $X_i$ in (\ref{eq:vfXi}) follow by (cf. Remark~\ref{remark:AlgebraGen})
$$\Scale[0.9]{X_1\big|_{R} = R \, (\partial_{x} R)|_{\Id},\,\, X_2\big|_{R} = R \, (\partial_{\th} R)|_{\Id}, \,\, X_3\big|_{R} = R \, (\partial_{y} R)|_{\Id}.}$$
\end{remark}
We also use the isomorphism between $\so$ and $\R^3$
\begin{equation}
\label{so3r3iso}
A_i \sim \ul{e}_i,\quad R A_i R^{-1} \sim R \ul{e}_i,
\end{equation}
where $A_i \in \so$, $R \in \SO$, $\ul{e}_i \in \R^3, \, i = 1,2,3.$

Note that~(\ref{matexp}) is a product of matrix exponentials:
\begin{equation}
\label{Rexp}
R(x,y,\th) = \eexp{y A_3} \eexp{-x A_2} \eexp{\th A_1}.
\end{equation}
We choose to rely on this parameterization to keep the analogy with previous $\SE$ models~\cite{DuitsJMIV2014,cut_sre2,SIAM} mentioned in introduction.
\subsection{Statement of the Problem $\Pcurve$}\label{subsec:Pcurvestat}
Let $\S^{2} = \{\ul{n} \in \R^3 \, \big| \, \|\ul{n}\| = 1\}$ be a sphere of unit radius. We consider the problem $\Pcurve$ (see Fig.~\!\ref{fig:PcurveStat}), which is for given boundary points $\ul{n}_0, \ul{n}_1 \in \S^2$ and directions $\ul{n}_0' \in T_{\ul{n}_0}(\S^2)$, $\ul{n}_1' \in T_{\ul{n}_1}(\S^2)$, $\|\ul{n}_0'\| = \|\ul{n}_1'\|=1$ to find a smooth
curve $\ul{n}(\cdot):[0,l]\to \S^2$ that satisfies the boundary conditions
\begin{eqnarray}\label{pcurveboundsgen}
\ul{n}(0) = \ul{n}_0, \ \ul{n}(l) = \ul{n}_{1}, \quad \ul{n}'(0) = \ul{n}_0', \ \ul{n}'(l) = \ul{n}_1',
\end{eqnarray}
and for given $\xi>0$ minimizes the functional
\begin{eqnarray} \label{EPc}
\mathcal{L}(\ul{n}(\cdot)) := \int_0^l \gothic{C}(\ul{n}(s))\sqrt{\xi^2 + k_g^2(s)} \, {\rm d} s,
\end{eqnarray}
where $k_g(s)$ denotes the geodesic curvature  on $\S^{2}$ of $\ul{n}(\cdot)$ evaluated in time $s$, and $\gothic{C}: S^2 \to [\delta,+\infty)$, $\delta>0$, is a smooth function that we call ``external cost.''

Here the total length $l$ is free and $s=\int_0^s 1 \, {\rm d} \sigma = \int_0^s \|\ul{n}'(\sigma)\| \, {\rm d} \sigma$ denotes the spherical arclength. Thus, we have $\|\ul{n}'(s)\|=1$ and the Gauss--Bonnet formula
\begin{equation}
\label{eq:GaussBonnet}
k_g(s) = \ul{n}''(s)\cdot (\ul{n}(s) \times \ul{n}'(s)).
\end{equation}
\begin{remark}\label{remark:20}
In introduction, we provided two motivations for this problem, where the external cost $\gothic{C}$ plays a different role. Firstly, there is a motivation coming from retinal images, where $\gothic{C}$ is adapted to spherical image data. Secondly, there is a motivation from the modeling of human vision, where the nonuniform distribution of photoreceptors can be modeled by the external cost.
\end{remark}
\begin{remark}
\label{rm1}
Without loss of generality, we can fix the boundary conditions as
\begin{equation}\label{pcurveboundsfix}
\ul{n}(0) = \ul{e}_1, \, \, \ul{n}(l) = \ul{n}_{1} , \quad
\ul{n}'(0) = \ul{e}_3, \, \, \ul{n}'(l) = \ul{n}_1',
\end{equation}
since the problem is left-invariant
w.r.t. the natural left action of $\SO$ onto $\S^{2}$ for $\gothic{C}=1$. But also for $\gothic{C}\neq 1$ we can fix these boundary conditions by shifting the external cost. 
This boundary value convention is used throughout this article for $\Pcurve$. 
\end{remark}
\begin{remark}
\label{rmNonExist}
Following the previous works \cite{Boscain2014,DuitsJMIV2014} in $\SE$ group, we do not expect existence of minimizers for the problem $\Pcurve(S^2)$ in general.
\end{remark}
\subsection{Statement of $\PMEC$ Problem in $\SO$}\label{subsec:Pmecstat}
We call $\PMEC$ the following SR problem in $\SO$:
\begin{eqnarray}
&\dot R = - u_1 R A_2 + u_2 R A_1, \label{eq:so3sys} \\
&R(0) = \Id, \qquad R(T) = R_{fin}, \label{eq:so3bounds}\\
&\mathcal{L}(R(\cdot)) := \int \limits_{0}^{T} C(R(t))\sqrt{\xi^2 u_1^2(t) +  u_2^2(t)}\, {\rm d} t \to \min, \label{intsrlength}\\
&\qquad R \in \SO, \quad (u_1, u_2) \in \R^2, \quad \xi>0, \label{eq:so3domain}
\end{eqnarray}
with $T>0$ free.

The external cost $C:\SO \to [\delta, +\infty)$, $\delta>0$, is a smooth function that is typically obtained by lifting the external cost $\gothic{C}$ from the sphere $S^2$ to the group $\SO$, i.e., $C(R) = \gothic{C}(R \, \ul{e}_1)$.

\begin{figure}[ht]
\centering
\includegraphics[width=0.55\hsize]{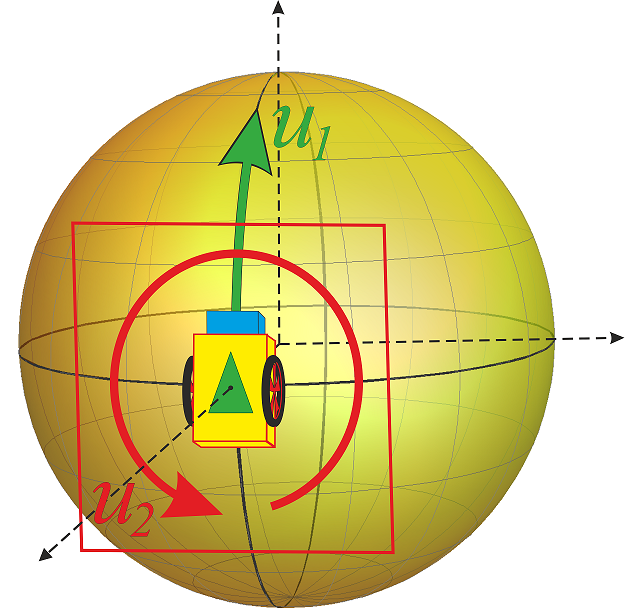}
\caption{The controls $u_1$ and $u_2$ along the ``spatial generator'' $X_1$ and the ``angular generator'' $X_2$ (cf. Remark~\!\ref{remark:ReedsShepp}).
\label{fig:ReedsShepp}}
\end{figure}

We study the problem $\PMEC$ for $C=1$ (case of uniform external cost) in Sect.~\ref{sec:PMEC}, but let us first consider some preliminaries.
\begin{remark}
In $\Pmec$, we only have two velocity controls $u_1$ and $u_2$ in a 3D tangent space $T_{R(t)}(\SO)$ (cf. Fig.~\ref{fig:ReedsShepp}).
\end{remark}
\begin{remark}\label{rm:existence}
The existence of minimizers for the problem $\PMEC$ is guaranteed by Chow-Rashevskii and Filippov theorems on sub-Riemannian manifolds~\cite{agrachev_sachkov}.
Moreover, the geodesics of $\PMEC$ are smooth, since there are no abnormal extremals (see Example 1.3.13 in \cite{Rifford}).
\end{remark}
\begin{remark}\label{remark:SRManifold}
Sub-Riemannian manifolds are commonly defined by a triplet $(\mathcal{M}, \Delta, \mathcal{G})$, with manifold $\mathcal{M}$, distribution $\Delta \subset T(\mathcal{M})$ and metric tensor $\mathcal{G}$. In our case
\begin{equation}\label{eq:subriemmanifold}
\begin{array}{c}
\mathcal{M}= \SO,\quad \Delta = \spann\{R A_1, R A_2\}, \\[3pt]
\mathcal{G}(\dot{R}, \dot{R}) = C^2(\cdot)\left(\xi^2 u_1^2 + u_2^2\right),
\end{array}
\end{equation}
where the controls $u_1$ and $u_2$ are components of velocity vector w.r.t. moving frame of reference, see~(\ref{eq:so3sys}). The choice of the sub-Riemannian structure~(\ref{eq:subriemmanifold}) is determined by the initial conditions $n(0) = e_1$ and $n'(0) = e_3$ in~(\ref{pcurveboundsfix}).
\end{remark}
\begin{remark}\label{CauchySchwarz}
By virtue of Cauchy--Schwarz inequality, the problem of minimizing the sub-Riemannian length $\mathcal{L}$ is equivalent to the problem of minimizing the action functional
\begin{equation}\label{eq:so3int}
\mathcal{E}(R(\cdot))= \int_{0}^{T} C^2(R(t)) \frac{\xi^2 u_1^2(t) +  u_2^2(t)}{2}\, dt \to \min,
\end{equation}
with $T$ fixed (see, e.g.,~\cite{Montgomery}). 
\end{remark}
\begin{remark}\label{remark:ReedsShepp}
In analogy with the SR problem $\Pmec$ in $\SE$, cf.\!~\cite{DuitsJMIV2014,Boscain2014},
we sometimes call $X_{1}=-R A_{2}$ the ``spatial generator'' and $X_{2}=R A_{1}$ the ``angular generator'', despite the fact that $X_{1}$ and $X_{2}$ are both angular generators on $S^{2}$. The problem $\Pmec(\SO)$ can be seen as a model of the Reeds-Shepp car on a sphere. The Reeds-Shepp car can move forward and backward and rotate on a place. The input $u_1$ controls the motion along $X_1$, and the input $u_2$ controls the motion along $X_2$, see Fig.~\ref{fig:ReedsShepp}. 
\end{remark}
\subsection{Relation Between the Problems $\Pcurve$ and $\PMEC$}\label{subsec:connect}
\begin{figure}[ht]
\centering
\includegraphics[width=0.49\hsize]{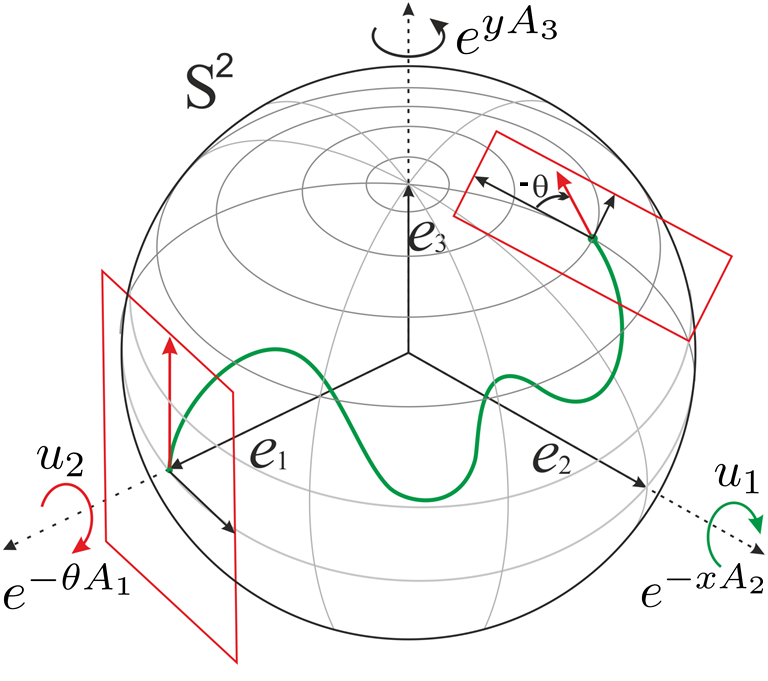}~
\includegraphics[width=0.45\hsize]{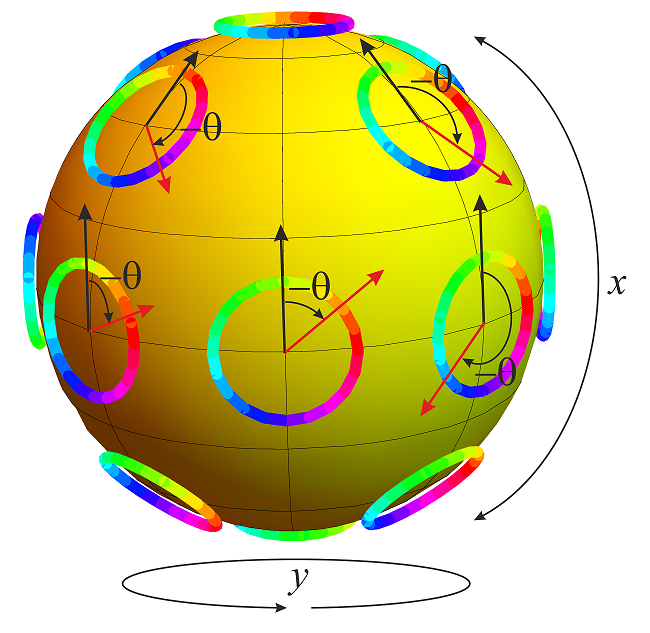}
\caption{\textbf{Left:} Illustration of the parameterizations used in $\Pmec$ and $\Pcurve$. The rotations are parameterized by (\ref{Rexp}), i.e. by angles $x$, $y$ and $\th$ of rotation about basis axes, see~(\ref{matexp}).
\textbf{Right:} Illustration of the Hopf fibration~\cite{Markina2} --- circle bundle with the base $S^2$.
\label{fig:x}}
\end{figure}

We call a spherical projection the following projection map from $\SO$ onto $\S^{2}$ (see Fig.~\ref{fig:x}):
\begin{equation} \label{proj}
\SO \ni R \mapsto R \, \ul{e}_{1} \in \S^2.
\end{equation}
In coordinates $(x,y,\th)$ defined by (\ref{matexp}), we have
\begin{equation}
\label{eq:proj}
R(x,y,\theta) \, \ul{e}_{1}=
\left(
\begin{array}{c}
\cos x \cos y \\
\cos x \sin y \\
\sin x
\end{array}
\right) = \ul{n}(x, y) \in \S^2.
\end{equation}
So we see that $(x,y)$ are spherical coordinates on $\S^2$.

In problem $\Pcurve$, one is interested in a curve
$\ul{n}(s) = R(t(s)) \, \ul{e}_1$,
which satisfies~(\ref{pcurveboundsgen}) and minimizes (\ref{EPc}). Here $R(t) = R(x(t),y(t),\th(t))$, and
\begin{equation}\label{eq:tfroms}
t(s) = \int_0^s \gothic{C}(\ul{n}(\sigma))\sqrt{\xi^2 + k_g^2(\sigma)} {\rm d} \, \sigma.
\end{equation}

Next we show that the spherical projection (\ref{proj}) of certain minimizers of $\Pmec$  provides solution of problem $\Pcurve$.
More precisely, this only holds for the minimizers whose spherical projection does not have a cusp.
\begin{definition}
The spherical projection of a minimizer of $\Pmec$ is said to have a cusp at $t=t^n_\mathrm{cusp}$ if there exists $\epsilon>0$, s.t. $u_{1}(a) u_{1}(b)<0$ for all
$a \in (t^n_\mathrm{cusp}-\epsilon, t^n_\mathrm{cusp})$ and $b \in (t^n_\mathrm{cusp},t^n_\mathrm{cusp}+\epsilon)$. That is, if the control in front
of the ``spatial'' generator $X_{1}$ switches sign locally. We are interested in the first cusp time $t_\mathrm{cusp} = \min \limits_{n \in \N} \{t^n_\mathrm{cusp}>0\}$, and we call $\smax$ the corresponding value of spherical arclength, s.t. $t_\mathrm{cusp} = t(\smax)$ via (\ref{eq:tfroms}).
\end{definition}
Notice that if $u_1 \equiv 0$ then the trajectory of $\Pmec$ is projected at a single point on $S^2$ which does not provide a solution to $\Pcurve$.
This allows us to define $\smax$ as
\begin{equation}\label{eq:smaxdef}
\smax = \inf\{s>0\, \big| \, u_1(t(s)) = 0\},
\end{equation}
so that $t_\mathrm{cusp}=t(\smax)$.

The following theorem states that minimizers $\ul{n}(s)$ of $\Pcurve$ for $s\in[0,\smax)$ are given by spherical projection of the minimizers $R(t)$ of $\Pmec$ for $t \in [0, \tcusp)$. 
Here $s$ denotes the spherical arclength parameter, defined by $\|\frac{\mathrm{d}}{\mathrm{d} s} \ul{n}(s)\| \equiv 1$, and $t$ denotes the SR arclength, defined by $C(R(t))\sqrt{\xi^2 u_1^2(t) + u_2^2(t)} \equiv 1$.
\begin{theorem}\label{prop:31}
Let $R(t)$, $t \in [0,T]$ be a minimizer of $\Pmec$ parameterized by SR-arclength, and let the corresponding optimal control $(u_1(t),u_2(t))$ satisfy the inequality $u_1(t)>0$ for all $t\in[0,T]$. Set
\begin{equation}
\label{pcurveRbounds}
\begin{cases}
\ul{n}_0 = \ul{e}_1, \qquad \ul{n}_1 = R(T)\, \ul{e}_1,\\
\ul{n}_0' = \ul{e}_3, \qquad \ul{n}_1' = R(T)\, \ul{e}_3.
\end{cases}
\end{equation}
Then for such boundary conditions, $\Pcurve$ has a minimizer $\ul{n}(s)$, along which we have
\begin{eqnarray}
&& \ul{n}(s) = R(t(s)) \, \ul{e}_1, \quad
\begin{cases}u_1(t)= \frac{{\mathrm{d}} s}{\mathrm{d} t}(t),\\ u_2(t)= 
k_g(s(t)) \frac{{\mathrm{d}} s}{\mathrm{d} t}(t) ,\end{cases}  \label{Pcurvecontrols}\\
&& t(s) = \int_0^s \gothic{C}(\ul{n}(\sigma)) \sqrt{\xi^2 + k_g^2(\sigma)} \, {\rm d} \, \sigma, \label{eq:ts}
\end{eqnarray}
for $0\leq s\leq l< \smax$, and $T = t(l)$.
\end{theorem}
\textbf{Proof } See Appendix~\ref{app:C}.
$\hfill \Box$
\section{The Sub-Riemannian Problem $\PmecOne$ in $\SO$}\label{sec:PMEC}
In this section, we study the problem $\Pmec$ in the special case of uniform external cost $C=1$, which we call $\PmecOne$. This problem can be seen as a left-invariant SR problem in Lie group $\SO$, and it was tackled by many authors (see \cite{Boscain_Rossi_B,Berestovskiy1,Berestovskiy2,Schcerbakova,Bonnard,Markina1,Markina2}) in different contexts. It is remarkable that the statement of the problem is very simple, but due to the high complexity of the formulas describing geodesics, the minimizers are still unknown for general $\xi>0$. The case $\xi=1$ corresponds to a symmetric SR structure in $\SO$, and it was completely solved in~\cite{Boscain_Rossi_B,Boscain2014,Berestovskiy1,Markina1,Markina2}. In this section, we consider the general case $\xi>0$ and derive analytic formulas for the geodesics
using the parametrization~(\ref{matexp}) for $\SO$. Further, in Sect.~\ref{sec:NumSol} these formulas allow us to verify our numerical approach to compute the global minimizers of $\Pmec$, presented for the first time by knowledge of the authors.

In our analytical study, first we introduce a coordinate chart in $\SO$ and write the optimal control problem $\PmecOne$ in this chart. Then we apply Pontryagin maximum principle (PMP), which is the first-order necessary condition for optimality~\cite{agrachev_sachkov}, and we derive the explicit formulas for the sub-Riemannian geodesics. 
Afterward we explain the notion of sub-Riemannian wavefront and sphere.
Finally, we discuss what is known about optimality of the geodesics. 
\subsection{$\PmecOne$ as Optimal Control Problem on $\S^2 \times \S^1$}\label{subsec:PMECchart}
In this section, we introduce the chart in $\SO$
and formulate problem $\PmecOne$ given by (\ref{eq:so3sys})--(\ref{eq:so3domain}), see Section \ref{subsec:Pmecstat}, as an optimal control problem on $\mathcal{M}=\S^2_{x,y}\times \S^1_{\th}$, where $(x,y) \in [-\frac\pi2,\frac\pi2] \times \R/\{2 \pi \Z\}$ are the spherical coordinates, cf.~(\ref{eq:proj}), and $\th \in \R/\{2 \pi \Z\}$. We use this specific chart to stress the analogy with the closely related problem in $\SE$ group (see~\cite{DuitsJMIV2014,Boscain2014}).

Consider a smooth curve $\gamma(\cdot) = (x (\cdot),y (\cdot),\th (\cdot)) \in C^{\infty}(\R, \S^2_{x,y} \times \S^1_{\th})$ departing from the origin, i.e., $\gamma(0) = e= (0,0,0)$. In Sect.~\ref{subsec:mathnot} we parameterized the group $\SO$ by the angles $x$, $y$, $\th$, recall~(\ref{matexp}). A smooth curve $\gamma(\cdot)$ corresponds to a smooth curve $R(\cdot)$ in $\SO$ defined by the one-parameter family of matrices
\begin{equation}
\label{Rcurvet}
R(t) = R(x(t),y(t),\th(t))\in \SO,
\end{equation}
depending smoothly on the parameter $t\in \R$ and satisfying the initial condition $R(0) = \Id$. A tangent vector of $R(t) = R(x(t),y(t),\th(t))$ is expressed as
$$\dot{R} = \frac{\partial R}{\partial x} \dot x + \frac{\partial R}{\partial y} \dot y + \frac{\partial R}{\partial \th} \dot \th.$$
Therefore, the control system (\ref{eq:so3sys})  can be rewritten as
\begin{equation}
\label{eq:controlSystem}
\dot{\gamma}=
\left(\begin{array}{c}
\dot{x}\\
\dot{y}\\
\dot{\th}
\end{array}\right)=\left(\begin{array}{c}
\cos \th \\
-\sec x \sin\th\\
 \tan x \sin\th
\end{array}\right)u_{1}+\left(\begin{array}{c}
0\\
0\\
1
\end{array}\right)u_{2}.
\end{equation}
Vector fields near the controls $u_i$ are two of three of the basis left--invariant vector fields in $\SO$ expressed in our chart. We denote them by $X_1$, $X_2$, and the third basis left invariant vector field is given by their commutator $X_3 = [X_1,X_2]$. Thus we have
\begin{equation}
\label{X1X2}
\Scale[0.89]{
\begin{array}{c}
X_1\!=\!
\left(\!\!\begin{array}{c}
\cos \th\\
-\sec x \sin \th\\
 \tan x \sin\th
\end{array}\!\!\right)\!\!, \, \,
X_2\!=\!\left(\!\begin{array}{c}
0\\
0\\
1
\end{array}\!\right)\!\!, \,\,
X_3\!=\!\left(\!\!\begin{array}{c}
\sin\th\\
\sec x \cos\th \\
-\tan x \cos\th
\end{array}\!\!\right).
\end{array}
}
\end{equation}
Then problem $\PmecOne$ given by~(\ref{eq:so3sys})--(\ref{eq:so3domain}) for $C=1$, taking into account Remark~\ref{CauchySchwarz}, is equivalent to the following optimal control problem:
\begin{eqnarray}
&\dot \gamma = u_1 X_1 + u_2 X_2, \quad (u_1, u_2) \in \R^2\label{eq:so3vsys} \\
&\gamma(0) = e,\quad \gamma(T) = \nu_{fin} = (x_1, y_1,\th_1) \in \S^2 \times \S^1, \label{eq:so3vbounds}\\
&\mathcal{E} = \frac12 \int_{0}^{T}\left(\xi^2 u_1^2 +  u_2^2\right)\, {\rm d} \,t \to \min, \quad \xi>0. \label{intaction}
\end{eqnarray}

\begin{remark}
\label{rmProjective}
A projective version of problem $\Pmec^1$ with $\xi=1$ is studied in \cite{Boscain_Rossi_S2}. It is obtained by identification of antipodal directions $\th$ and $\th + \pi$; thus, its solution is given by minimum between two geodesics with terminal configurations $h = (x_1,y_1,\th_1)$  and $\bar{h} = (x_1,y_1,\th_1 +\pi)$. The associated SR distance is given by  
$$
d^P(g,h) =
\min \limits_{ \scriptsize  \begin{array}{c}\gamma(\cdot)\in \operatorname{Lip}([0,T],\SO),\\ \gamma(0)=g, \gamma(T) \in \{h,\bar{h}\},\\ \dot{\gamma}(\cdot) \in \Delta\end{array}} \mathcal{L}(\gamma(\cdot)).
$$ The model in~\cite{Boscain_Rossi_S2} is connected to problem $\Pmec^1$ via right shift $g \mapsto g\cdot e^{\pi A_3}e^{-\frac{\pi}{2} A_2}$. The local chart is obtained by linear transformation $a = \frac{\pi}{2} - x$, $b = y$, $\theta = \xi + \pi$. The existence of optimal trajectories whose spherical projection contains a cusp is mentioned in~\cite{Boscain_Rossi_S2}.
\end{remark}
\subsection{Application of the Pontryagin Maximum Principle}
In this subsection, we apply the Pontryagin maximum principle (PMP) to problem $\PmecOne$ given by~(\ref{eq:so3vsys})--(\ref{intaction}) and derive the Hamiltonian system of PMP. Here we follow standard geometric control theory (see, e.g.,~\cite{agrachev_sachkov}). 

Due to the absence of abnormal extremals (see Remark~\ref{rm:existence}), we consider only the normal case, where the control-dependent Hamiltonian of PMP reads as
\[\begin{array}{l}
H_u(\lambda,\nu) = \langle\lambda, u_1 X_1 + u_2 X_2 \rangle -\frac{\xi^2 u_1^2 +  u_2^2}{2}, \\[3pt]
\left. \qquad \qquad \qquad \qquad \qquad \text{ with } \lambda = \sum \limits_{k=1}^{3} \lambda_k {\rm d} \, \nu^k \in T^{\ast}_\nu \mathcal{M}. \right.
\end{array}\]
Here $\nu = (x,y,\theta)$, $X_i$ are given by~(\ref{X1X2}), $\langle\cdot,\cdot\rangle$ denotes the action of a covector on a vector, and $({\rm d} \, \nu^1, {\rm d} \, \nu^2, {\rm d} \, \nu^3) = ({\rm d} \,  x, {\rm d} \,  y, {\rm d} \,  \theta)$ are basis one forms.

The maximization condition of PMP reads as
\begin{eqnarray*}
&& H_{u(t)}(\lambda(t),\gamma(t)) = \\
&& \Scale[0.9]{\displaystyle{\max_{(u_1,u_2)\in\R^2}}\left(\langle\lambda(t), u_1 X_1|_{\gamma(t)} + u_2 X_2|_{\gamma(t)} \rangle -\frac{1}{2}\left(\xi^2 u_1^2 +  u_2^2\right)\right),}
\end{eqnarray*}
where $u(t)$ denotes an extremal control and $(\lambda(t),\gamma(t))$ is an extremal.

Introduce left-invariant Hamiltonians $h_i = \langle\lambda,X_i \rangle$, $i = 1,2,3$. The maximization condition gives the expression for the extremal controls
\begin{equation} \label{eq:extremalcontrols}
u_1 = \frac{h_1}{\xi^2},\qquad u_2 = h_2.
\end{equation}
Then the maximized Hamiltonian, which we call ``the Hamiltonian'', reads as
\begin{equation}\label{eq:theHamiltonian}
H = \frac12 \left(\frac{h_1^2}{\xi^2} +h_2^2 \right).
\end{equation}
The vertical part (for momentum components $h_i$) of the Hamiltonian system of PMP is given by
\begin{equation}\label{eq:poiswithHam}
\dot{h}_i = \left\{H,h_i\right\} = \sum_{j=1}^3 \left( \frac{\partial H}{\partial \nu_j} \frac{\partial h_i}{\partial \lambda_j} - \frac{\partial h_i}{\partial \nu_j} \frac{\partial H}{\partial \lambda_j }\right),
\end{equation}
where $\{\cdot,\cdot\}$ denotes the Poisson bracket (see~\cite{agrachev_sachkov}).

Using the standard relation between Poisson and Lie brackets $\{h_i,h_j\} = \langle \lambda,[X_i,X_j]\rangle$, we get
\begin{eqnarray*}
&&\left\{H,h_1\right\}  =  h_2 \left\{h_2,h_1\right\} =  - h_2 h_3, \\
&&\left\{H, h_2\right\}  = \frac{1}{\xi^2}\,  h_1 \left\{h_1,h_2\right\} = \frac{1}{\xi^2}\, h_1 h_3,\\
&&\left\{H,h_3\right\}  \!=\! \frac{1}{\xi^2}h_1 \left\{h_1,h_3\right\} + h_2 \left\{h_2,h_3\right\} \!=\! \left(\!1-\frac{1}{\xi^2}\!\right) h_1 h_2.
\end{eqnarray*}

Thus the Hamiltonian system with the Hamiltonian $H$ follows by the last three identities and~(\ref{eq:poiswithHam}):
\begin{equation}
\Scale[0.95]{
\label{eq:hamsys}
\begin{array}{ll}
\begin{cases}
\dot{h}_1 = - h_2 h_3,\\
\dot{h}_2 = \frac{1}{\xi^2}\, h_1 h_3, \\
\dot{h}_3 = \left(1-\frac{1}{\xi^2}\right) h_1 h_2
\end{cases} 
& \begin{cases}
\dot{x} = \frac{h_1}{\xi^2} \cos \th , \\
\dot{y} = -\frac{h_1}{\xi^2} \sec x \sin \th, \\
\dot{\th} = \frac{h_1}{\xi^2} \sin\th \tan x  + h_2
\end{cases} \\ 
 \text{--- vertical part,}
&  \text{--- horizontal part,}
\end{array}
}
\end{equation}
with the boundary conditions
\begin{equation}\label{eq:hamsysBound}
\begin{array}{l l l}
h_1(0) = h_1^0, \,& \,h_2(0) = h_2^0, \, & \, h_3(0) = h_3^0, \\
x(0) = 0, \,&\, y(0) = 0,  \,&\, \th(0) = 0.
\end{array}
\end{equation}
In~(\ref{eq:hamsys}), the horizontal part follows from~(\ref{eq:extremalcontrols}) and~(\ref{eq:controlSystem}).

It is known that the Hamiltonian system~(\ref{eq:hamsys}) is Liouville integrable~\cite{Mashtakov_Sachkov_Integr,Kruglikov}, and it was explicitly integrated in~\cite{Ivan,Schcerbakova}. In the next subsections, we classify the possible solutions by values of the parameter $\xi$, and we adapt the explicit solution to our coordinate chart $(x,y,\th) \in \mathcal{M}$, where we follow the analogy to the closely related problem in $\SE$. This allows us to obtain simpler formulas than in~\cite{Ivan,Schcerbakova}.
\subsection{Classification of Different Types of Solutions}\label{subsec:classification}
\begin{figure}[ht]
\begin{minipage}[b]{0.31\linewidth}
\centering
\includegraphics[width=\textwidth]{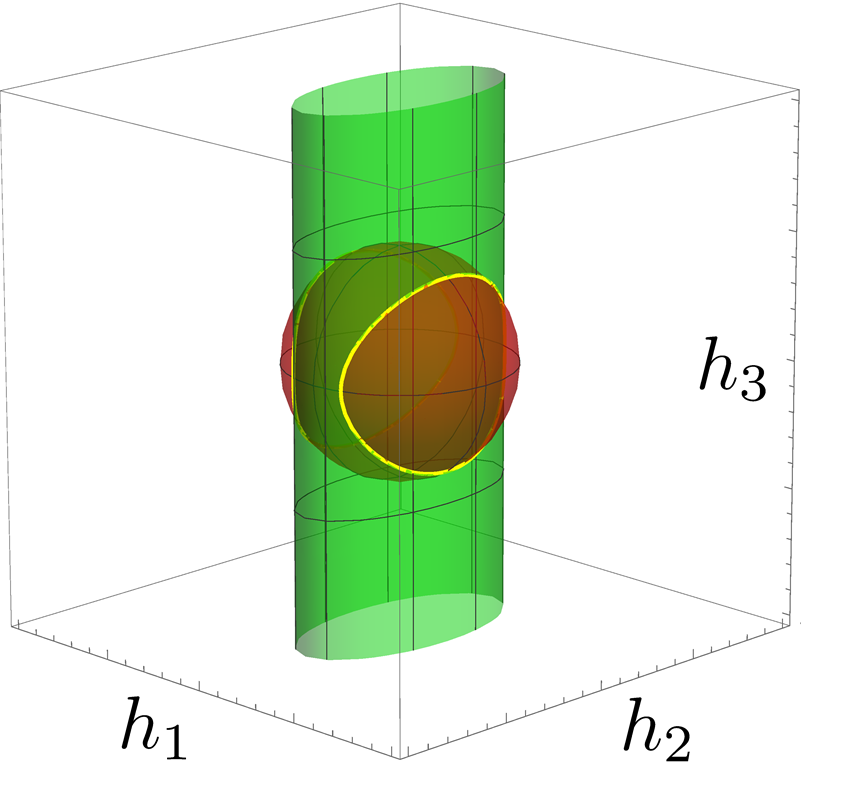}
Elliptic $(\xi<1)$
\end{minipage}
~
\begin{minipage}[b]{0.31\linewidth}
\centering
\includegraphics[width=\textwidth]{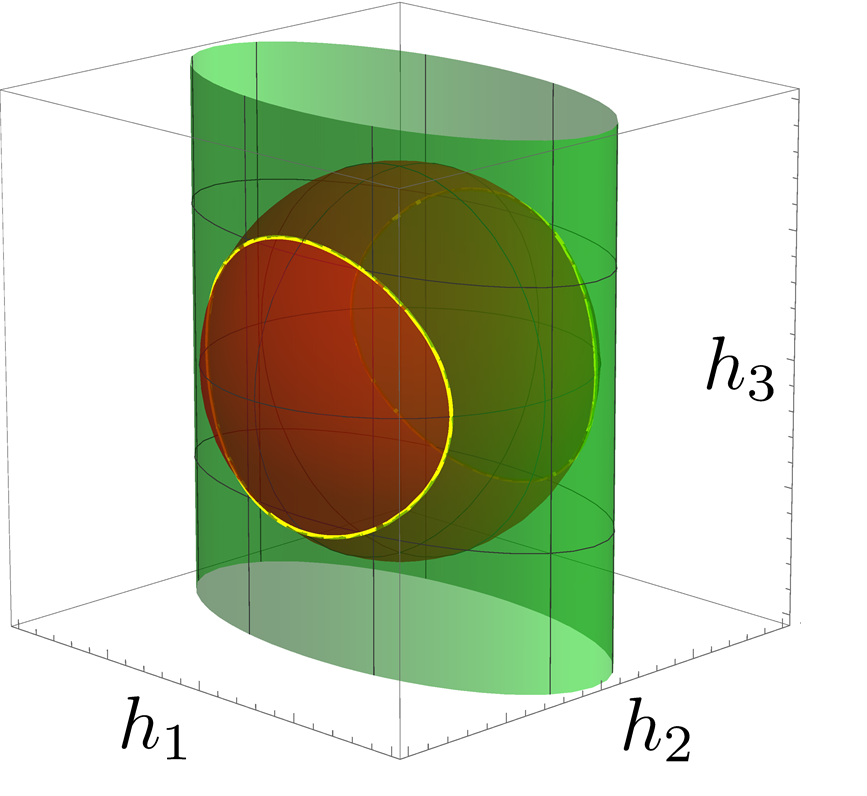}
Linear $(\xi=1)$
\end{minipage}
~
\begin{minipage}[b]{0.31\linewidth}
\centering
\includegraphics[width=\textwidth]{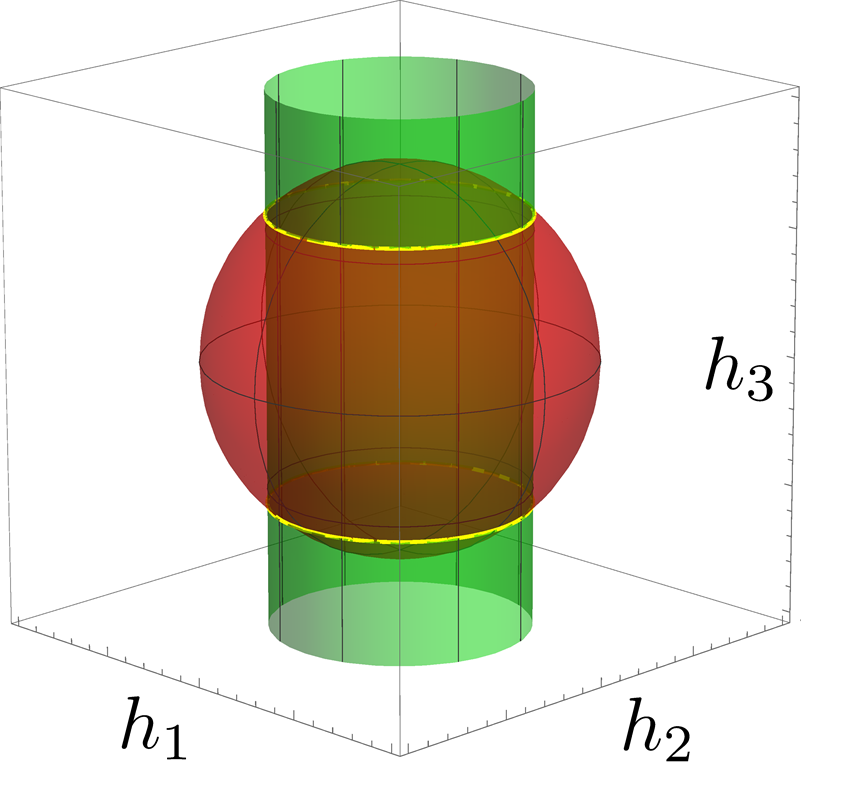}
Hyperbolic $(\xi>1)$
\end{minipage}
\centering
\newline
\begin{minipage}[b]{0.31\linewidth}
\centering
\includegraphics[width=\textwidth]{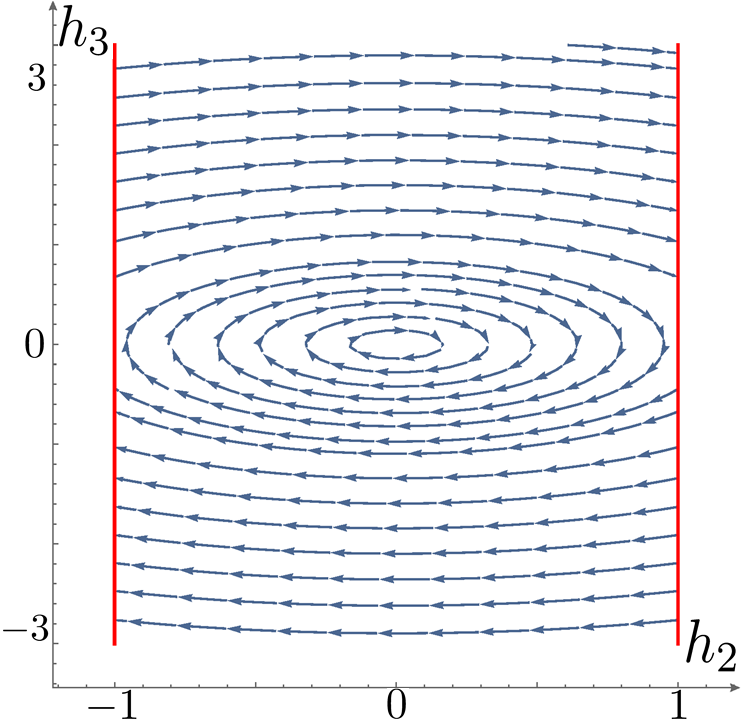}
\end{minipage}
~
\begin{minipage}[b]{0.31\linewidth}
\centering
\includegraphics[width=\textwidth]{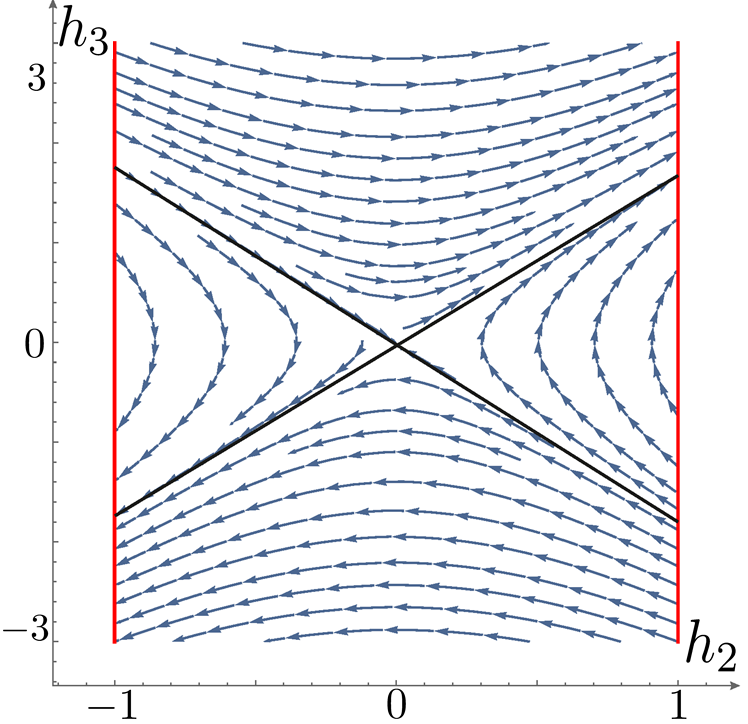}
\end{minipage}
~
\begin{minipage}[b]{0.31\linewidth}
\centering
\includegraphics[width=\textwidth]{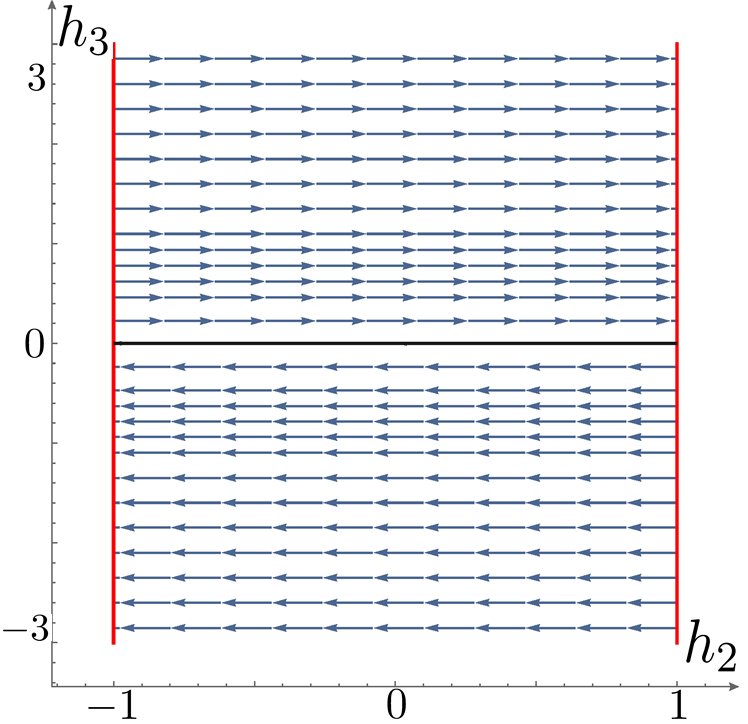}
\end{minipage}
\caption{\textbf{Top:} Integral manifolds of the vertical part of~(\ref{eq:hamsys}). The red sphere is a coadjoint orbit $M^2=h^2_1+h^2_2+h^2_3$ and the green cylinder is the Hamiltonian $2 H= \frac{h^2_1}{\xi^2} + h^2_2 = 1$.
\newline
\textbf{Bottom:} Vector field plot on the $(h_2,h_3)$ plane.}
\label{fig:IntegralsVP}
\end{figure}

Here we describe the domains of parameter $\xi>0$, which correspond to different dynamics of the Hamiltonian system (\ref{eq:hamsys}). It is well known (see, e.g.,~\cite{Schcerbakova}) that this system has two integrals of motion that depend only on momentum components $h_i$ (see illustration in Fig.~\!\ref{fig:IntegralsVP})
\begin{eqnarray}
&& 2 H = \frac{h_1^2}{\xi^2} + h_2^2 = 1 \text{ - the Hamiltonian, }\label{eq:theHamiltonianTwice}\\
&& M^2 = h_1^2 + h_2^2 + h_3^2 \text{ - Casimir function.}\label{eq:coadjointorbit}
\end{eqnarray}
Different values of the Hamiltonian $H$ correspond to different speed along extremal trajectories $\gamma(\cdot)$. By fixing the value $H = \frac12$ in (\ref{eq:theHamiltonianTwice}), we use SR arclength parameter
$$t = \int_0^t \sqrt{\mathcal{G}|_{\gamma(\tau)}(\dot{\gamma}(\tau),\dot{\gamma}(\tau))}\, {\rm d} \,\tau = \int_0^t 1 \, {\rm d} \,\tau$$ along extremal trajectories. The momentum trajectory $h(t)$ can be seen as the curve formed by intersection of the cylinder $H=\frac12$ and the sphere $M = \text{const}$ (see the yellow line in Fig.~\!\ref{fig:IntegralsVP}). 

Elimination of $h_1$ from (\ref{eq:theHamiltonianTwice}) and (\ref{eq:coadjointorbit}) yields
$$
h_2^2(1 - \xi^2) + h_3^2 = M^2-\xi^2.
$$
Since $h^2_2\leq1$, there exist the following types of solution of the vertical part of~(\ref{eq:hamsys}):
\begin{itemize}
\item[1)] Elliptic for $\xi<1$, $M^2>\xi^2$;
\item[2)] Linear for $\xi=1$, $M^2>\xi^2$;
\item[3)] Hyperbolic for $\xi>1$, $M^2\neq\xi^2$, $M>1$;
\item[4)] A point (0,0) for $\xi<1$, $M^2=\xi^2$;
\item[5)] A segment of a straight line for $\xi = 1$, $M=1$;
\item[6)] Two crossing segments of straight lines for $\xi>1$, $M^2=\xi^2$;
\item[7)] No solution otherwise.
\end{itemize}

According to this classification, we will refer to the case $\xi<1$ as the elliptic case, $\xi=1$ as the linear case, and $\xi>1$ as the hyperbolic case.
\subsection{Geodesics in SR Arclength Parametrization}\label{subsec:geodesicstpar}
In this subsection, we provide explicit formulas for the SR geodesics, where we reexamine results from~\cite{Ivan,Schcerbakova}.

In the coordinates $(\beta, c) \in \R/\{4\pi\Z\} \times \R$, defined by
\begin{equation}\label{eq:polarcoords}
h_1 = \xi \cos\frac{\beta}{2}, \; \; h_2 = \sin\frac{\beta}{2}, \; \; h_3 = \frac{\xi c}{2},
\end{equation}
 the vertical part of (\ref{eq:hamsys}) becomes a system of mathematical pendulum
\begin{equation}
\label{pendeqsysc}
\begin{cases}
\dot{\beta} =  c, & \beta(0) = 2 \arg\left(\frac{h_1^0}{\xi} + \i h_2^0\right),\\
\dot{c} = -r \sin \beta, & c(0) = \frac{2 h_3^0}{\xi},
\end{cases} 
\end{equation}
where $r = \frac{1}{\xi^2} -1 \in (-1, +\infty)$.
Thus $r \in (-1, 0)$ corresponds to the elliptic, $r = 0$ to the linear, and $r \in (0,+\infty)$ to the hyperbolic case.

The system~(\ref{pendeqsysc}) has a symmetry $(\tilde{\beta}(t) = \beta(t) + \pi, \; \tilde{r} = -r)$, which allows us to study the elliptic case and obtain the hyperbolic case by symmetry observations. Note also that the linear case should be treated separately, but the equations in this case are much simpler. Therefore, in the remainder of this manuscript we restrict ourselves to the case $r\geq 0 \desda 0<\xi\leq 1$.


In the following theorem, we provide explicit formulas for SR geodesics in SR arclength parametrization. 
\begin{theorem}\label{thmGeogTpar}
A solution of the system~(\ref{eq:hamsys}) reads as
\begin{equation}
\label{xythsol}
\Scale[0.97]{
\begin{array}{l}
h_1(t) =  \xi \cos\frac{\beta(t)}{2}, \ \ h_2(t) =  \sin\frac{\beta(t)}{2}, \ \ h_3(t) = \frac{\xi c(t)}{2}, \\[4pt]
x(t) = \arg(\sqrt{R_{11}^2(t) + R_{21}^2(t)} + \i R_{31}(t)),\\
y(t) = \arg(R_{11}(t) + \i R_{21}(t)),\\
\th(t) = \arg(R_{33}(t) + \i R_{32}(t)),
\end{array}
}
\end{equation}
where $(\beta(t), c(t))$ is a trajectory of~(\ref{pendeqsysc}), and
\begin{eqnarray}
\Scale[0.87]{R(t)\!=\!
\Scale[0.82]
{
\left(\begin{array}{c c c}
R_{11}(t) & R_{12}(t) & R_{13}(t)\\
R_{21}(t) & R_{22}(t) & R_{23}(t)\\
R_{31}(t) & R_{32}(t) & R_{33}(t)\\
\end{array}\right)
}
\!=\!
D_0^{T} \eexp{\tilde{y}(t) A_3} \eexp{-\tilde{x}(t) A_2} \eexp{\tilde{\th}(t) A_1},} \label{eq:RtotildeR}
\end{eqnarray}
\begin{eqnarray*} \label{D0}
&&\begin{array}{l}
\text{ where } D_{0}= \frac{1}{M} \left(
\begin{array}{ccc}
\mu & \frac{h_2^0 h_1^0}{\mu} & -\frac{h_2^0 h_3^0}{\mu} \\
0 & \frac{M h_3^0}{\mu} & \frac{M h_1^0}{\mu} \\
h_2^0 & -h_1^0 & h_3^0
\end{array}
\right),
\end{array}\\
&& \left. \right.  \quad \text{ with~~~   } \displaystyle{\mu = \sqrt{M^2 - (h_2^0)^2}}, \\
&& \left. \right. \qquad \quad \quad \displaystyle{M = \sqrt{(h_1^0)^2+(h_2^0)^2+(h_3^0)^2}}, \\
&& \text{ and where }
\end{eqnarray*}
\begin{equation*}
\label{solstilde}
\begin{array}{l}
\left. \right.  \quad \displaystyle{\tilde{x}(t)= \arg \left(\sqrt{h^2_{1}(t)+h^2_3(t)}+ \i \, h_{2}(t) \right),} \\
\left. \right.  \quad \displaystyle{\tilde{y}(t)= \frac{M^2}{\xi^2}\left(t - \int_0^t \frac{h_3^2(\tau)}{M^2 - h_2^2(\tau)} d \tau \right),}\\
\left. \right.  \quad \displaystyle{\tilde{\th}(t)= \arg\left(h_3(t) - \i \, h_1(t)\right).}
\end{array}
\end{equation*}
\end{theorem}
\textbf{Proof } See Appendix~\ref{app:B}.
$\hfill \Box$

\begin{figure}[ht]
~~~~
\begin{minipage}[b]{0.35\linewidth}
\centering
\includegraphics[width=\textwidth]{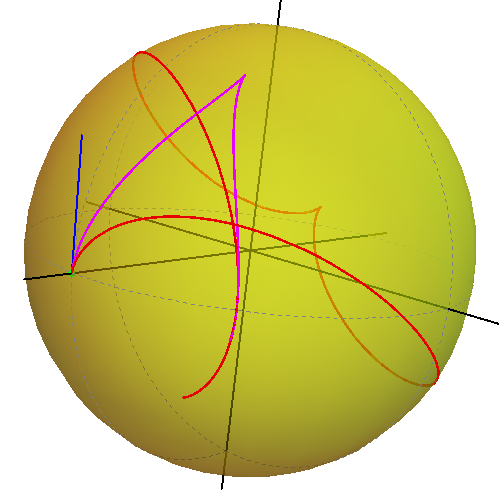}
\end{minipage}
\begin{minipage}[b]{0.35\linewidth}
\centering
\includegraphics[width=\textwidth]{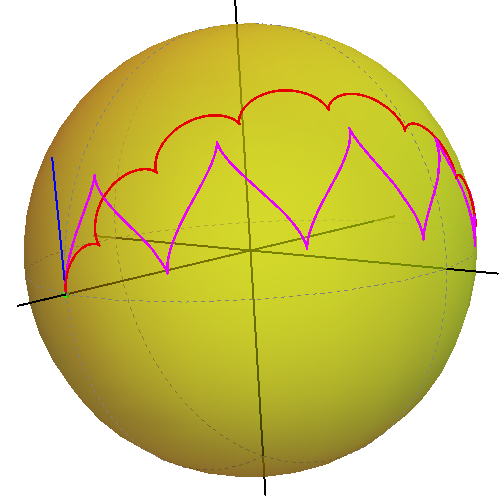}
\end{minipage}\\
\centering
\begin{minipage}[b]{0.41\linewidth}
\centering
\includegraphics[width=\textwidth]{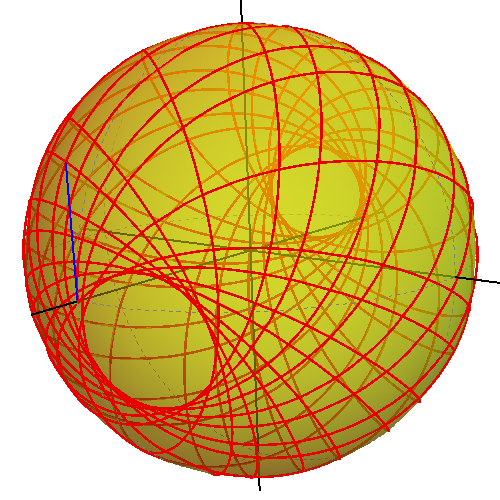}
\end{minipage}
\caption{Spherical projection of various sub-Riemannian geodesics in $\SO$.}\label{fig:extremalcurvesTpar}
\end{figure}

Note that $(\beta(t), c(t))$, and $\tilde{y}(t)$ admit explicit expression in Jacobi elliptic functions and elliptic integrals of the first, second and third kind~\cite{Ivan}.
We present plots of spherical projections~(\ref{proj}) of extremal trajectories in Fig.~\!\ref{fig:extremalcurvesTpar}. Here it is remarkable that
the spherical projections of SR geodesics in $\SO$ can represent cusps, 
which can be undesirable for image analysis applications. This motivates us to study SR geodesics before the first cusp in their spherical projections. For short, we call such curves ``cuspless geodesics.'' In Sect.~\ref{sec:pcurve} we describe the possible end conditions, which can be connected by a cuspless geodesic, and in Sect.~\ref{sec:NumSol} we provide PDE-based approach for solving the boundary value problem (BVP) for data-driven SR geodesics, where we reformulate the problem as a solution to the SR-eikonal equation.
\subsection{Sub-Riemannian Wavefronts and Spheres} \label{subsec:optimality} 
A useful tool for studying geodesics in left-invariant optimal control problems is {\it the exponential map}\footnote{Not to be confused with the exponential map from Lie algebra to Lie group.} that maps an initial momentum $h(0)$ and a time $t$ to the end point of corresponding geodesic $\gamma(\cdot)$ [i.e., the exponential map integrates the Hamiltonian system~(\ref{eq:hamsys})]:
\begin{equation}
\label{eq:expmap}
\Exp: \begin{array}{r l} T^*_e (\SO) \times \R & \to \SO, \\ (h(0),t) & \mapsto \gamma(t). \end{array}
\end{equation}
Now we explain the notion of {\it wavefront}. By definition, the wavefront consists of endpoints of all the geodesics of the same length $T$:
\begin{equation}
\label{eq:wavefront}
\WF(T) = \bigg\{ \Exp(h(0),T)\, \bigg| \, \begin{array}{l}h(0) \in T^{*}_{e} (\SO), \\ H(h(0))=\frac12 \end{array} \bigg\}.
\end{equation}
The outer surface of a wavefront forms a {\it sub-Riemannian sphere}, which is a set of endpoints in $\mathcal{M}$ equidistant from $e$:
\begin{eqnarray*}
\mathcal{S}(T) &&= \big\{ \Exp(h(0),T) \, | \, h(0) \in T^{*}_{e} (\SO),  \\
               && \qquad \qquad \qquad \qquad H(h(0))=\frac12, \, t_{cut}(h(0))\geq T\big\} \\
              &&= \left\{g \in \SO \, | \, d(e,g) = T\right\},
\end{eqnarray*}
where $d(q_1,q_2)$ denotes the sub-Riemannian distance 
between $q_1$ and $q_2$, and $t_{cut}$ denotes the cut time, i.e., an instance of time, when a geodesic loses its optimality.
\subsubsection*{What Is Known About Optimality?}
In the previous subsection, we computed SR geodesics or in other words extremal trajectories. It is known (see, e.g.,~\cite{Montgomery,agrachev_sachkov}) that sufficiently short arcs of SR geodesics are SR length minimizers (optimal trajectories). It is also known that in general a geodesic loses its optimality after a cut point. The corresponding instance of time is called cut time, and the set of all cut points in configuration space forms the so-called cut locus.
In general, it is complicated to derive cut loci~\cite{agrachev_sachkov}.
\begin{remark}\label{remark:MaxwellConj}
There are two reasons for a geodesic to lose optimality:\\
$\left. \quad \right.$ 1) Local optimality is lost at a conjugate point (critical point of exponential map, that integrates the Hamiltonian system);\\
$\left. \quad \right.$ 2) Global optimality is lost at a first Maxwell point (when two geodesics meet with the same length for the first time).
\end{remark}

In problem $\PmecOne$, both the Maxwell set and the conjugate locus are not known. The estimation of the first Maxwell time was given in~\cite{Ivan}, but still neither the cut locus, nor the SR length minimizers are known yet for general $\xi>0$. The cut locus in a special case of the bi-invariant metric (for $\xi=1$) was obtained in~\cite{Boscain_Rossi_B}. In~\cite{Schcerbakova}, conjugate locus in the Riemannian case was constructed, and by considering a SR metric as a limiting case of the Riemannian metric, the corresponding formulas for the conjugate locus in SR case could be obtained.

Here we do not provide such a limiting procedure. Instead, motivated by applications,  in Sect.~\ref{sec:NumSol} we propose a numerical solution to compute the SR length minimizers and spheres. In Sect.~\ref{sec:pcurve}, motivated by the study~\cite{DuitsJMIV2014}, we also give some statements on optimality of cuspless geodesics. We show that in contrast to $\SE$ case~\cite{DuitsJMIV2014,Boscain2014} there exist nonoptimal cuspless SR geodesics.
\section{Sub-Riemannian Geodesics in $\SO$ with Cuspless Spherical Projection}\label{sec:pcurve}
In this section, we study cuspless SR geodesics. Such geodesics allow parametrization by spherical arclength, which leads to simpler formulas than using SR arclength parametrization (see Sect.~\ref{sec:PMEC}).
Spherical arclength parameterization breaks down, when the spherical projection of a geodesic exhibits a cusp for the first time. So a natural question arises to characterize the set of end conditions $\mathcal{R} \subset \SO$ reachable by cuspless geodesics, similar to the $\SE$ case studied in~\cite{DuitsJMIV2014}.

The cuspless SR geodesics are projections in $\SO$ of the trajectories of the Hamiltonian system (\ref{eq:hamsys}) that are going in positive direction of $X_1$ (i.e., $u_1(t)>0$) before the first cusp time $t<\tcusp$. Thus, by virtue of~(\ref{eq:extremalcontrols}) the cuspless constraint is given by
$$h_1(t)>0, \text{ for all } t \in [0,T].$$
We shall often rely on short notation $h_{i}(s):=h_{i}(t(s))$, where we stress our notations
$$h'_{i}= \frac{\mathrm{d}}{\mathrm{d}s}h_{i} \text{ and } \dot{h}_{i}=\frac{\mathrm{d}}{\mathrm{d}t}h_i$$
in order to avoid confusion with the chain law for differentiation. Recall~(\ref{eq:smaxdef}) and note that $\smax$
for a given geodesics is determined by its initial momentum. We write $\smax(h(0))$ when we want to stress this dependence.
We derive an explicit formula for $\smax(h(0))$ in Sect.~\ref{subsec:smax}.
\begin{theorem} \label{thmHamSysspar}
For any $s \in [0,\smax(h(0)))$, $h_1(0)>0$ the system (\ref{eq:hamsys}) is equivalent to the following system (see corresponding vector flow in fig.~\ref{fig:IntegralsVP}):
\begin{equation}\label{eq:hamsysspar}
\Scale[0.88]{\!\!
\begin{array}{ll}
\begin{cases}
h_1(s) = \xi^2\frac{\mathrm{d} s}{\mathrm{d} t} \geq 0,  \\
h_2'(s) =  h_3(s), \\
h_3'(s) = (\xi^2-1) h_2(s),
\end{cases}
\!\!&\!\!
\begin{cases}
x'(s) =  \cos \th(s),\\
y'(s) = - \sec x(s) \sin \th(s),\\
\th'(s) = \sin\th(s) \tan x(s)  + \\
\qquad \, \, \, \,\, \xi h_2(s)/\sqrt{1-h_2^2(s)},\\ 
\end{cases} \\ 
 \text{--- vertical part,}
&  \text{--- horizontal part,}
\end{array}\!\!
}
\end{equation}
with the boundary conditions
\begin{equation}\label{eq:hamsyssparBound}
\begin{array}{l}
h_2(0) = h_2^0, \qquad h_3(0) = h_3^0, \\
x(0) = 0, \,\, y(0) = 0,  \,\, \th(0) = 0.
\end{array}
\end{equation}
\end{theorem}
\textbf{Proof }
Express the dynamics (\ref{eq:hamsys}) in the spherical arclength parameter $s$. We have
\begin{equation*}
\Scale[0.95]{
\begin{cases}
\begin{array}{l l}
h_2'(s) = \frac{\mathrm{d} h_2(s)}{\mathrm{d} s} = \frac{\dot{h}_2(t(s))}{h_1(t(s))}\xi^2 = h_3(s), \,&\, h_2(0) = h_2^0,\\
h_3'(s) = \frac{\mathrm{d} h_3(s)}{\mathrm{d} s} = \frac{\dot{h}_3(t(s))}{h_1(t(s))}\xi^2 = \left(\! \xi^2- 1\!\right) h_2(s), \,&\, h_3(0) = h_3^0.
\end{array}
\end{cases}
}
\end{equation*}
Here we recall that $u_{1} = \frac{\mathrm{d} s}{\mathrm{d} t}$, see (\ref{Pcurvecontrols}), and $h_1 = \xi^2 u_1$, see (\ref{eq:extremalcontrols}).
Thus we obtained the vertical part of (\ref{eq:hamsysspar}). Doing the same for the horizontal part in (\ref{eq:hamsys}), we obtain the full system (\ref{eq:hamsysspar}), where the Hamiltonian $H(s) = \frac{\xi^{-2}h_1^2(s) + h_2^2(s)}{2} = \frac12$ together with the cuspless constraint $h_1(s)>0$ allows us to write $h_1(s) = \xi \sqrt{1- h_2^2(s)}$. \\
Finally, the boundary conditions (\ref{eq:hamsyssparBound}) follow from (\ref{eq:hamsysBound}), since $t = 0 \Leftrightarrow s =0$.
$\hfill \Box$

Note that in contrast to the mathematical pendulum ODEs~(\ref{pendeqsysc}) expressed in SR arclength $t$, where the Jacobi elliptic functions appear, the vertical part in $s$ parameterization now is a simple linear system of ODEs integrated in elementary functions. To obtain simpler formulas, we define the parameter $\chi \in \C$ (where one should take principal square root):
\begin{equation}
\label{xitochi}
\chi=\sqrt{\xi^2-1} = \begin{cases}\sqrt{1-\xi^2} \, \i, \text{ for } 0<\xi<1,\\ \sqrt{\xi^2-1} , \text{ for } \xi \geq 1.  \end{cases}
\end{equation}
\begin{theorem} \label{thmVertPartSParSol}
The general solution of the vertical part in (\ref{eq:hamsysspar}) for all $\xi\neq 1$ is given by
\begin{equation}
\label{vertpartchi}
\begin{cases}
h_2(s) = h_2^0 \cosh s \chi + \frac{h_3^0}{\chi} \sinh s \chi,\\
h_3(s) = h_3^0 \cosh s \chi + \chi h_2^0 \sinh s \chi,\\
h_1(s) = \xi \sqrt{1- h_2^2(s)},
\end{cases}
\end{equation}
and for the case $\xi=1$, we find straight lines parallel to $h_2$ axis in the $(h_2,h_3)$-phase portrait:
\begin{equation}\label{pendlinsolxi1}
\begin{cases}
h_2(s) = h_2^0+h_3^0\, s,   \\
h_3(s) = h_3^0, \\
h_1(s) =  \sqrt{1- h_2^2(s)}.
\end{cases}
\end{equation}
\end{theorem}
\textbf{Proof }
The solution to the system of ODEs is unique and it is readily checked (\ref{vertpartchi}) is a solution of the vertical part of (\ref{eq:hamsysspar}) for all $\xi\neq1$.  For $\xi<1$, we take main values for the complex square root and complex $\sinh$ and $\cosh$. For both $\xi=1$ and $\xi \downarrow 1$ and $\xi \uparrow 1$, we
have $h_2(s) = h_2^0+h_3^0 s$ and $h_3(s) = h_3^0$.
$\hfill \Box$
\subsection{Computation of the First Cusp Time}\label{subsec:smax}
To analyze the cusp points in problem $\Pcurve$, we need to determine $\smax(h(0))$. It is given, recall~(\ref{eq:smaxdef}) and~(\ref{eq:extremalcontrols}), by the minimal positive root of equation $h_1(s) = 0$:
\[ \smax(h(0)) = \min\{s>0\,|\, h_1(s) = 0\}.\]
\begin{theorem}\label{thmSmax}
When moving along a SR geodesic $t \mapsto \gamma(t)$ the first cusp time is computed as $\tcusp(h(0)) = t(\smax(h(0))$, recall (\ref{eq:tfroms}), where 
\begin{eqnarray}
\label{eq:smax}
\Scale[0.95]{
\smax(h(0))\!=\!
\begin{cases}
\begin{array}{l l}
\frac{\sgn(h_3^0)-h_2^0}{h_3^0} & \text{ for } \chi = 0, \, h_3^0 \neq 0,\\[3pt]
\frac{1}{\chi}\log \left(\frac{s_1\left(\sqrt{\kappa}+\chi\right)}{h_2^0 \chi +h_3^0}\right) & \text{ for } \begin{array}{l}
                                                                                                            \chi\neq0, \, \kappa\geq 0, \\
                                                                                                                h_2^0 \chi +h_3^0\neq 0, \end{array}\\[3pt]
+\infty &\text{ otherwise},
\end{array}
\end{cases}
}
\end{eqnarray}
with
$s_1 = \sgn\left(\Re\left(h_2^0 \chi +h_3^0\right)\right)$ \\ and $\kappa = (h_3^0)^2+\left(1-(h_2^0)^2\right) \chi ^2 \in \R.$

As a result, in the SR manifold $(\SO, \Delta, \mathcal{G})$, recall~(\ref{eq:subriemmanifold}), there do exist nonoptimal cuspless geodesics.
\end{theorem}
\textbf{Proof } See Appendix~\ref{app:smax}.
$\hfill \Box$

The presence of nonoptimal cuspless geodesics is remarkable, as in the $\SE$ group every cuspless geodesic is globally optimal~\cite{DuitsJMIV2014}.
\subsection{Geodesics in Spherical Arclength Parametrization}
\begin{figure}[ht]
\begin{minipage}[b]{0.31\linewidth}
\centering
$\xi<1$\\
\includegraphics[width=\textwidth]{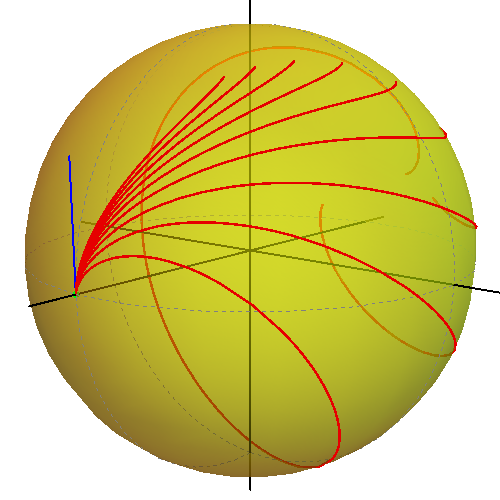}
\end{minipage}
~
\begin{minipage}[b]{0.31\linewidth}
\centering
$\xi=1$\\
\includegraphics[width=\textwidth]{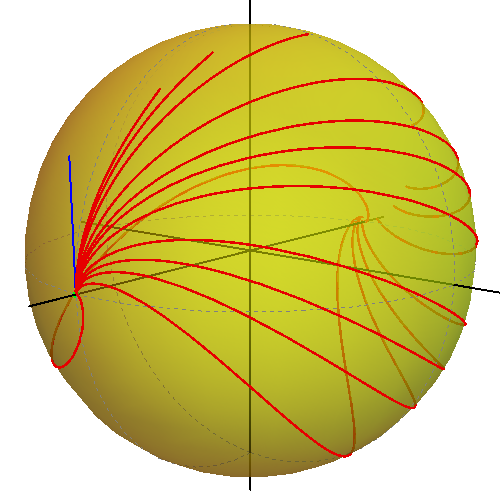}
\end{minipage}
~
\begin{minipage}[b]{0.31\linewidth}
\centering
$\xi>1$\\
\includegraphics[width=\textwidth]{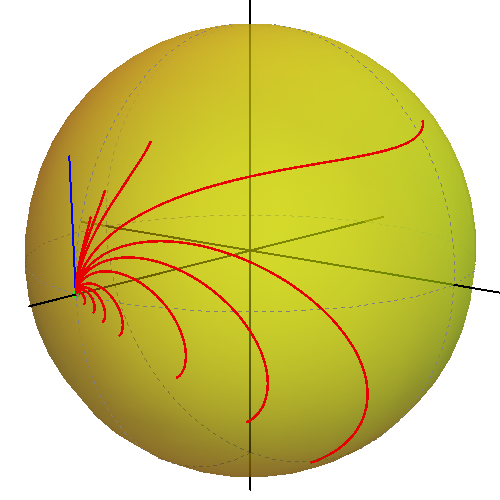}
\end{minipage}
\caption{Spherical projection of cuspless SR geodesics in $\SO$ in elliptic, linear and hyperbolic cases.
\label{fig:extremalcurveE}}
\end{figure}

In this subsection, we derive the formulas for cuspless SR geodesics in $s$ paramete\-rization.
\begin{theorem}\label{thmGeogSpar}
The unique solution of ~(\ref{eq:hamsysspar}) is defined for $s \in [0, \smax(h(0)))$, where $\smax(h(0))$ is given by~(\ref{eq:smax}). The solution to the vertical part is given by Theorem~\ref{thmVertPartSParSol} and the solution to the horizontal part is given by
\begin{equation}
\label{xythsolspar}
\begin{array}{l}
x(s) = \arg(\sqrt{R_{11}^2(s) + R_{21}^2(s)} + \i R_{31}(s)), \\
y(s) = \arg(R_{11}(s) + \i R_{21}(s)), \\
\th(s) = \arg(R_{33}(s) + \i R_{32}(s)),
\end{array}
\end{equation}
where\\
$\Scale[0.97]{
R(s)\!=\!
\Scale[0.85]{
\left(\begin{array}{c c c}
R_{11}(s) & R_{12}(s) & R_{13}(s)\\
R_{21}(s) & R_{22}(s) & R_{23}(s)\\
R_{31}(s) & R_{32}(s) & R_{33}(s)\\
\end{array}\right)}
\!=\!D_0^{T} \eexp{\tilde{y}(s) A_3} \eexp{-\tilde{x}(s) A_2} \eexp{\tilde{\th}(s) A_1},
}$
\begin{equation*} \label{D0}
\begin{array}{l}
D_{0}= \frac{1}{M} \left(
\begin{array}{ccc}
\mu & \frac{\xi h_2^0 \sqrt{1-(h_2^0)^2}}{\mu} & -\frac{h_2^0 h_3^0}{\mu} \\
0 & \frac{M h_3^0}{\mu} & \frac{\xi M \sqrt{1-(h_2^0)^2}}{\mu} \\
h_2^0 & -\xi \sqrt{1-(h_2^0)^2} & h_3^0
\end{array}
\right),
\end{array}
 \text{ and }
\end{equation*}
\begin{equation*}
\label{solstildeSpar}
\begin{array}{l}
\displaystyle{\tilde{x}(s)= \arg \left(\sqrt{M^2-h^2_{2}(s)}+ \i h_{2}(s) \right),} \\
\displaystyle{\tilde{y}(s)= \xi M^2\left(\int \limits_{0}^{s} \frac{\sqrt{1- h^2_{2}(\sigma)}}{M^{2}-h^2_{2}(\sigma)}{\rm d}\sigma\right),}\\
\displaystyle{\tilde{\th}(s)= \arg\left(h_{3}(s)- \i \xi\sqrt{1-h^2_{2}(s)}\right),}
\end{array}
\end{equation*}
with $\Scale[0.89]{\mu = \sqrt{M^2 - (h_2^0)^2}, \ \M = \sqrt{\xi^2(1-(h_2^0)^2)+(h_2^0)^2+(h_3^0)^2}}$.
\end{theorem}
\textbf{Proof }
It follows from Theoremes~\ref{thmGeogTpar},~\ref{thmHamSysspar}, and~\ref{thmVertPartSParSol}. Momentum component $h_{1}$ is expressed in $h_{2}$ via the Hamiltonian~(\ref{eq:theHamiltonian}) which equals to $\frac12$ along SR geodesics.
The sign of $h_{1}(0)$ is equal to the sign of $u_{1}(0)=\frac{\mathrm{d}s}{\mathrm{d}t}(0)$ and therefore nonnegative.\\
Regarding the integration constraints, we note that $M^2=h^2_{1}+h^2_{2}+h^2_{3}$
and $0\leq h_{1}=\xi\sqrt{1-h^2_{2}}$ and
\[
\frac{\mathrm{d}\tilde{y}}{\mathrm{d}s}= \frac{\xi^2}{h_{1}} \frac{\mathrm{d} \tilde{y}}{\mathrm{d} t}=
 M^2 \left(\frac{h_{1}}{M^2- h^2_{2}}\right)= \xi M^2
\frac{\sqrt{1-h^2_2}}{M^{2}-h^2_{2}}
\]
from which the result follows.
$\hfill \Box$
\begin{remark}
In contrast to general SR geodesics in $\SO$ given in $t$ parameterization, where Jacobi elliptic functions appear, the cuspless SR geodesics in $\SO$ given in $s$ parameterization involve only a single elliptic integral $\int \limits_{0}^{s} \frac{\sqrt{1- h^2_{2}(\sigma)}}{M^{2}-h^2_{2}(\sigma)}{\rm d}\sigma$. This integral can be expressed in terms of standard elliptic integrals. For example, in the elliptic case $\xi<1$, $M^2>\xi^2$ 
we have
$$h_2(s) 
       = \rho \sin \Psi(s),$$
where $\rho=\sqrt{\frac{M^{2}-\xi^2}{1-\xi^2}}$, $\Psi(s)= s\sqrt{1-\xi^2}  + \Psi_0$,\\
with $\Psi_0=\arg ( \frac{h_{3}(0)}{\sqrt{1-\xi^2}} + \i \,h_{2}(0))$, which yields
\[
\begin{array}{ll}
\tilde{y}(s) =& M  \,\int \limits_{0}^{s} \frac{\sqrt{1-\rho^{2}\sin^{2}(\Psi(\sigma))}}{M^2-\rho^{2} \sin^{2}(\Psi(\sigma))} {\rm d}\sigma  =\\
&\frac{M}{\sqrt{1-\xi^2}} \int \limits_{\Psi_0}^{\Psi(s)} \frac{\sqrt{1-\rho^{2}
\sin^{2}(\Psi)}}{M^2-\rho^{2} \sin^{2}(\Psi)} {\rm d}\Psi = \\
 &\frac{M}{\sqrt{1-\xi^2}}\Bigg(F\left(\Psi(s),\rho^2\right)- F(\Psi_0,\rho^2) -\left(1-\frac{1}{M^{2}}\right) \times \\
 &\left(\Pi \left(\frac{\rho^2}{M^2},\Psi(s),\rho^2\right)-\Pi \left(\frac{\rho^2}{M^2},\Psi_0,\rho^2 \right)\right)\Bigg),
\end{array}
\]
where $F$ denotes an elliptic integral of the first kind \\ and $\Pi$ denotes an elliptic integral of the third kind.
\end{remark}
See plots of projected cuspless geodesics $\ul{n}(s) = R(s)\, \ul{e}_1$ for $s \in (0,\smax)$  in Fig.~\!\ref{fig:extremalcurveE}.
\section{PDE Approach for Data-Driven SR Geodesics in $\SO$ }\label{sec:NumSol}
In this section, we adapt the PDE approach for data-driven SR geodesics in $\SE$~\cite{SSVM,SIAM} to the $\SO$ group.
 Here we consider the basis left-invariant vector fields $X_i$ as differential operators of the first order, and we write $X_i(\mathcal{W})$ for the derivative of a function $\mathcal{W}:\SO \to \R$ along $X_i$.

We aim to solve the following optimal control problem:
\begin{eqnarray*}
&\dot \gamma(t) = \sum \limits_{i=1}^{2} u_i X_i|_{\gamma(t)}, \text{ for } t \in [0,T], \\
&\gamma(t) \in \SO, \quad \gamma(0) = e, \ \gamma(T) = g_1, \quad (u_1, u_2) \in \R^2,\\
&l(\gamma(\cdot)) = \int_0^{T} C(\gamma(t))\sqrt{\xi^2 u_1^2(t) + u_2^2(t)} \, \mathrm{d} t \to \min.
\end{eqnarray*}
Here the terminal time $T$ is free; and $C : \SO \to [\delta, +\infty)$, $\delta > 0$ is the external cost.

By rescaling of time $\tau = \frac{t}{T} \in [0,1]$ simultaneously with controls  $\bar{u}_i(\tau) = T u_i(T \tau)$,
we write down the explicit solutions as
\begin{equation} \label{SRdist}
\begin{array}{l}
d_{\SO}^0(g,e):= \\\inf \limits_{
{\scriptsize
\begin{array}{c}
\gamma(0)=e, \gamma(1)=g, \\
\dot{\gamma}(\cdot) \in \left.\Delta\right|_{\gamma(\cdot)}, \\
\gamma \in \textrm{Lip}([0,1],\SO)
\end{array}
}
}\! \int \limits_{0}^{1}\! \sqrt{\left.\mathcal{G}\right|_{\gamma(\tau)}(\dot{\gamma}(\tau),\dot{\gamma}(\tau))} \, {\rm d}\tau .
\end{array}
\end{equation}
Here the sub-Riemannian metric tensor
\begin{equation} \label{SRMT}
\mathcal{G} = C^2(\cdot)\left(\xi^2 \omega^1\otimes \omega^1 + \omega^2\otimes \omega^2\right)
\end{equation}
is defined only on the distribution $\Delta$, recall Remark~\ref{remark:SRManifold}.

\subsubsection*{Sub-Riemannian Fast Marching in $\SO$}\label{subsec:FastMarching}
Here we propose a SR-FM method (sub-Riemannian fast marching) for the computation of data-driven SR length minimizers (not necessarily cuspless) in $\SO$ group, as a solution to the SR-eikonal system~(\ref{eq:eikonalsystfull}).
This method was successfully used in~\cite{CIARP} for the computation of data-driven SR length minimizers in $\SE$ group. The method is based on a Riemannian approximation of sub-Riemannian manifold, and computing Riemannian geodesics in highly anisotropic space, which becomes the SR manifold in the limiting case as aniso\-tropy tends to infinity.

Here we follow the explanation in~\cite{CIARP}, where we work now in new settings of the $\SO$ group and use the coordinate chart $(x,y,\th)$ defined in Sect.~\ref{subsec:PMECchart}. Recall that the basis left-invariant vector fields\footnote{In previous works~\cite{SSVM,SIAM}, $X_i$ was denoted by $\mathcal{A}_i$} $X_i$ in $\SO$ are given by the following differential operators:
\begin{eqnarray*}\label{eq:vecfield_diffop_again}
&& X_1 = \cos\th \, \partial_x - \sec x \sin\th \, \partial_y + \tan x \sin \th \, \partial_\th,\\
&& X_2 = \partial_\th,\\
&& X_3 = \sin\th \, \partial_x + \sec x \cos\th \, \partial_y - \tan x \cos \th \, \partial_\th,
\end{eqnarray*}
and corresponding basis left-invariant one forms $\omega^i$, satisfying $\langle \omega^i, X_j \rangle = \delta_i^j$, are expressed as
\begin{equation}\label{eq:basisforms}
\begin{array}{l}
\omega^1 = \cos\th \, {\rm d} x - \cos x \sin\th \, {\rm d} y,\\
\omega^2 = \sin x \, {\rm d} y + {\rm d}\th,\\
\omega^3 = \sin\th \, {\rm d} x + \cos x \cos\th \, {\rm d} y.
\end{array}
\end{equation}
We approximate the sub-Riemannian manifold by a Riemannian manifold by fixing a small $\varepsilon > 0$. Moreover, the SR-eikonal equation~(\ref{eq:eikonalsystfull}) is well defined and it can be derived as a limiting case of the eikonal equation on a Riemannian manifold via the \emph{inverse} metric tensor,
see Appendix~\ref{app:D}.  Let us denote the Riemannian distance as follows:
\begin{equation} \label{Rdist}
\begin{array}{l}
d_{\SO}^\varepsilon(g,e):= \\
\inf \limits_{
{\scriptsize
\begin{array}{c}
\gamma(0)=e, \gamma(1)=g, \\
\gamma \in \textrm{Lip}([0,1],\SO)
\end{array}
}
}\!\int \limits_{0}^{1}\! \sqrt{\left.\mathcal{G}_{\varepsilon}\right|_{\gamma(\tau)}(\dot{\gamma}(\tau),\dot{\gamma}(\tau))}\, {\rm d}\tau .
\end{array}
\end{equation}
with Riemannian metric tensor $\mathcal{G}_{\varepsilon}$ given by
\begin{equation} \label{RMT}
\Scale[0.99]{
\mathcal{G}_{\varepsilon} := C^2(\cdot)\left(\xi^2 \omega^1\otimes \omega^1 + \omega^2\otimes \omega^2 + \frac{\xi^2}{\varepsilon^2}\, \omega^3 \otimes \omega^3\right).
}
\end{equation}
\begin{remark}\label{remark:viscos}
In our approach, we shall rely on standard notion 
of viscosity solution \cite{Lions,Crandall,Evans} of the Riemannian eikonal equation, which admits a generalization to the sub-Rieman\-nian case. See details in Appendix~\ref{app:D}.
\end{remark}
The following theorem summarizes our approach for the computation of data-driven sub-Riemannian length minimizers in $\SO$.
\begin{theorem} \label{th:PDE}
Let $g \neq e \in \SO$ be chosen such that there exists a \underline{unique} minimizer $\gamma_{\varepsilon}^*:[0,1] \to \SO$
of $d_{\SO}^{\varepsilon}(g,e)$
such that $\gamma_{\varepsilon}^*(\tau)$ is not a conjugate point for all
$\tau \in [0,1]$
and all $\varepsilon \geq 0$ sufficiently small.

Then
$\tau \mapsto d_{\SO}^{0}(e,\gamma^*_0(\tau))$ is smooth and the minimizer $\gamma_{0}^*(\tau)$ is given by
$\gamma_{0}^*(\tau)=\gamma_{b}^*(1-\tau)$, with
\begin{equation}\label{eq:backtrackmain}
\begin{cases}
\begin{array}{l}
\dot{\gamma}_b(\tau)=-\bar{u}_1(\tau) \left.X_{1}\right|_{\gamma_b(\tau)}\!-\! \bar{u}_2(\tau) \left.X_{2}\right|_{\gamma_b(\tau)}, \\
\gamma_{b}(0)=g,
\end{array}
\end{cases}
\end{equation}
with
$$\Scale[0.95]{
(\bar{u}_1(\tau), \bar{u}_2(\tau)) = \frac{\mathcal{W}(g)}{(C(\gamma_b(\tau)))^2} \left( \frac{X_{1}|_{\gamma_b(\tau)}(\mathcal{W})}{\xi^2}, X_{2}|_{\gamma_b(\tau)}(\mathcal{W})\right),}$$
and with $\mathcal{W}(g)$ denoting the viscosity solution of the following boundary value problem: 
\begin{equation} \label{eq:eikonalsystfull}
\left\{
\begin{array}{l}
\sqrt{
\frac{\left(X_{1}|_g(\mathcal{W})\right)^2}{\xi^{2}} + \left(X_{2}|_g(\mathcal{W})\right)^2}=C(g),  \textrm{ for }g\neq e, \\
\mathcal{W}(e)=0.
\end{array}
\right.
\end{equation}
\end{theorem}
For an outline of the proof, see Appendix~\ref{app:D},
where we rely on a similar approximation approach as in \\\mbox{\cite[ch.5,app.A]{CHENDAPhDthesis}} for an elastica functional
and in \\
\mbox{\cite[Thm.2,Cor.2,app.A]{DuitsMeestersMirebeauPortegies}} for a SR problem on $SE(2)$.

\begin{remark}\label{remark:cusplessModification}
The approach in Theorem~\ref{th:PDE} can be adapted for producing only the cuspless minimizers.
See Appendix~\ref{app:D1}.
\end{remark}
\begin{remark}
The SR spheres of radius $t$ centered at $e$ are given by $\mathcal{S}_{t}=\! \{g \! \in \! \SO \,|\, \mathcal{W}(g)=t \}$.
\end{remark}
\begin{remark}
Recall that conjugate points are points where local optimality is lost. For a formal definition and characterization, see \cite[def.8.4.3,cor.8.4.5]{AgrachevBarilariBoscain}.
We would conjecture that the assumption on the conjugate points in the above theorem is not really needed,
as even the conjugate points that are limits of first Maxwell points do not seem to cause problems in the backtracking procedure,
akin to the $SE(2)$ case \cite[App.D]{SIAM}, but nevertheless our proof in Appendix~\ref{app:D} does rely on this assumption.
\end{remark}
From the derivations above Theorem~\ref{th:PDE}, we see that
the fast marching approach for computing SR geodesics in $\SE$, cf.~\cite{CIARP}, is easily generalized to the $\SO$ case.
To this end, we replace the matrix representation for $\mathcal{G}_{\varepsilon}$ expressed in the fixed $(x,y,\theta)$ Cartesian coordinate frame. In the $\SO$ case, it equals
\begin{equation*}
M_{\varepsilon}=
\gothic{R} \left(
\begin{array}{ccc}
 C^2(\cdot)\, \xi^2 & 0 & 0 \\
 0 &  C^2(\cdot) &0 \\
 0& 0& C^2(\cdot)\, \xi^2\, \varepsilon^{-2} 
\end{array}
\right) \gothic{R}^T,
\end{equation*}
with
\begin{equation*}
\gothic{R} = \left(
\begin{array}{c c c}
 \cos\th & \ 0 \ & \sin \th \\
 -\cos x\, \sin \th & \ \sin x \ & \cos x \, \cos\th\\
 0& \ 1 \ & 0
\end{array}
\right)\ . 
\end{equation*}
Here the diagonal matrix in the middle encodes the anisotropy between the $X_i$
directions, while the matrix $\gothic{R}$ is the basis transformation from the moving coframe $\{\omega^1,\, \omega^2,\, \omega^3\}$ to the fixed coframe $\{{\rm d} x,\, {\rm d} y,\, {\rm d} \th\}$, recall~(\ref{eq:basisforms}), in which the fast marching implementation via special anisotropic stencils~\cite{mirebeau} is used.

In Sect.~\ref{subsec:Experiment_UniformCost}, we will show that the thereby obtained fast marching approach already presents reasonable precision for
$\varepsilon=0.1$. Experiments in Sect.~\ref{subsec:Retinal_Example} show the application of the method (with data-adaptive nonuniform cost) to tracking of blood vessels in retinal images.

\section{Experiments}\label{sec:Experiments}
In Sect.~\ref{subsec:Experiment_UniformCost}, we verify the SR-FM method by comparison of SR length minimizers/spheres obtained via SR-FM with the exact SR geodesics/wavefronts (cf. Sect.~\ref{sec:PMEC}) for the case of uniform external cost (i.e. $C=1$).
In Sect.~\ref{subsec:ComparisonSE2SO3_UniformCost}, we compare SR geodesics and wavefronts in the groups $\SE$ and $\SO$ (Sect.~\ref{subsec:ComparisonSE2SO3geodesics}) for the uniform external cost case.
In Sect.~\ref{subsec:Retinal_Example}, we provide experiments of vessel tracking by SR geodesics in $\SO$ when the external cost $C$ is induced by spherical data, and compare them to the result of vessel tracking on the corresponding flat image by SR geodesics in $\SE$.
\subsection{Verification of the Fast Marching Method in the Case of Uniform External Cost}\label{subsec:Experiment_UniformCost}
In this subsection, we perform the experiments to validate the SR Fast Marching (SR-FM) method proposed in Sect.~\ref{subsec:FastMarching}.
The goal of the experiments is to check that the method produces an accurate approximate solution to the SR problem in $\SO$ group in the case of uniform external cost $C = 1$.
In all the experiments, we fixed the anisotropy parameter of the Riemannian approximation as $\varepsilon = 0.1$.

In the first experiment, we compare the geodesics $\gamma^{FM}(\cdot)$ obtained via SR-FM with the exact cuspless geodesics $\gamma(\cdot)$ computed via analytic formulas in Theorem~\ref{thmGeogSpar}. We perform the comparison as follows:\\
1) Fix $\xi>0$ and the initial momenta $h(0)$.
\\
2) Compute the first cusp time $\smax(h(0))$ corresponding to $h(0)$, and set $s_\mathrm{end} = \min\{\smax(h(0)),\frac{\pi}{2}\}$.
\\
3) Compute the geodesic $\gamma(s)$, $s\in[0,s_\mathrm{end}]$ via Thm.~\!\ref{thmGeogSpar}.
\\
4) Compute the distance function $\mathcal{W}(g)$ in the domain $g = (x,y,\theta) \in [-\frac{\pi}{2},\frac{\pi}{2}]\times[-\pi,\pi]\times[-\pi,\pi]$ via SR-FM. Here we compute the distance function in the grid of $201 \times 401 \times 401$ points and then interpolate it using third-order Hermite interpolation.
\\
5) Compute the geodesic $\gamma^{FM}(t)$, $t \in [0,\mathcal{W}(g_\mathrm{end})]$ via backtracking~(\ref{eq:backtrackmain}) from the endpoint $g_{1} = \gamma(s_\mathrm{end})$.
\\
6) Plot the spherical projections of $\gamma(s)$ and $\gamma^{FM}(t)$ and compare them.

\begin{figure}[ht]
\centering
\includegraphics[width=0.65\linewidth]{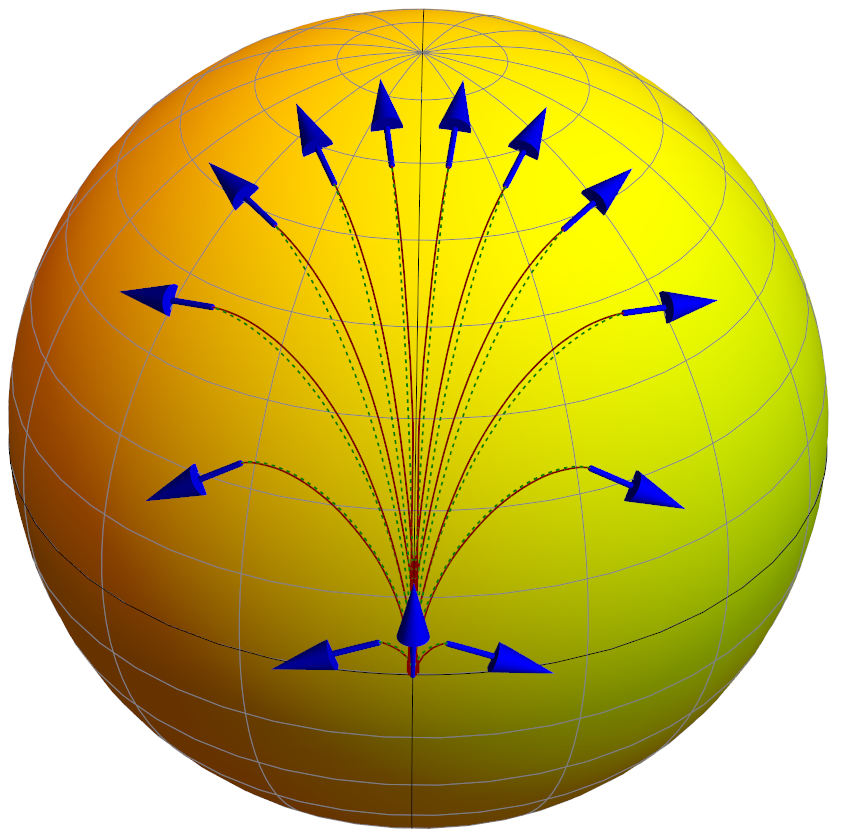}\\
\includegraphics[width=0.70\linewidth]{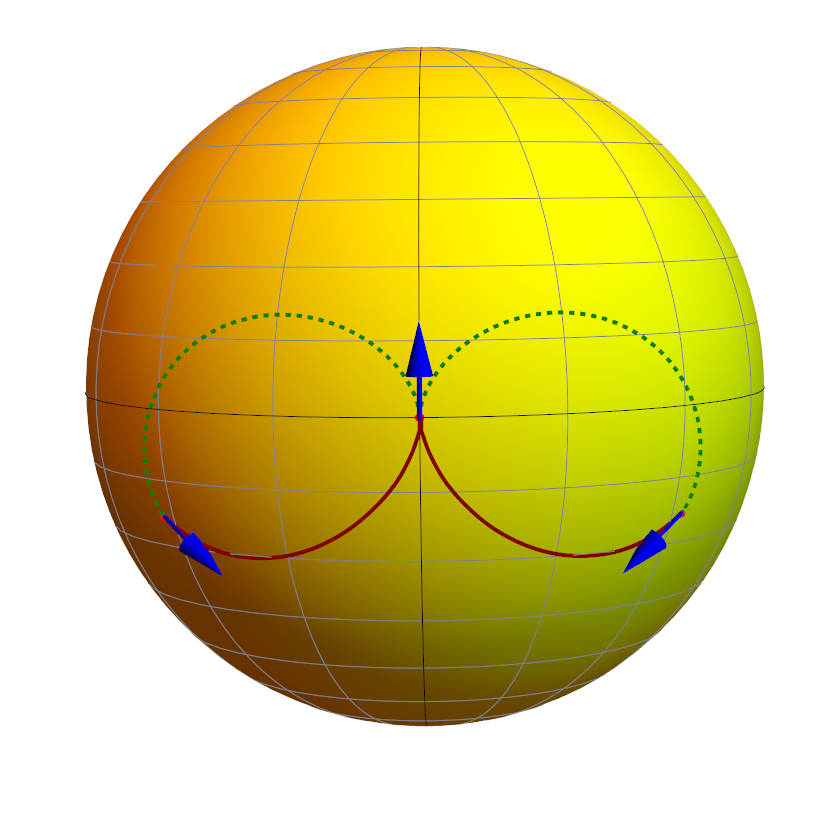}
\caption{\textbf{Top:}
The spherical projection of cuspless SR geodesics in $\SO$ computed by exact formulas in Theorem~\ref{thmGeogSpar} (green dashed lines) and our numerical Fast Marching approximations (red lines).\\
\textbf{Bottom:} Example of nonoptimal cuspless geodesics.} 
\label{fig:ExactGeodVsNumeric}
\end{figure}

A typical result of the comparison is shown in Fig.~\!\ref{fig:ExactGeodVsNumeric}, where we set $\xi = 1.5$, $h_3(0) = 0$ for all the curves and varied initial momentum $h_2(0) \in  \{-0.99$, $-0.81$, $-0.63$, $-0.45$, $-0.27$, $-0.09$, $0.09$, $0.27$, $0.45$, $0.63$, $0.81$, $0.99\}$. As a result, we see that the geodesics computed numerically via SR-FM accurately follow the exact geodesics.

We have performed a series of such experiments and always obtained similar results when the geodesics $\gamma^i(s)$ were optimal for $s \in [0, s^i_\mathrm{end}]$. It was also remarkable that SR-FM resulted into different curves (length minimizers) when the geodesic $\gamma^i(s)$ was not optimal for $s \in [0, s^i_\mathrm{end}]$. Such an example is illustrated in Fig.~\!\ref{fig:ExactGeodVsNumeric} (bottom). The question of optimality of SR geodesics in $\SO$ in the general case $\xi>0$ is an open important problem, recall Sect.~\ref{subsec:optimality}. Here we provide a numerical SR-FM method for computing only the optimal geodesics. In analogy with how it was done in~\cite{SIAM}, it is possible to compute Maxwell sets numerically. We have a strong conjecture that optimality is lost at the first Maxwell point induced by reflectional symmetries along the geodesic. We leave the computation of Maxwell points and analysis of optimality of the geodesics as a direction for future work.

\begin{figure}[ht]
\centering
\includegraphics[width=0.65\linewidth]{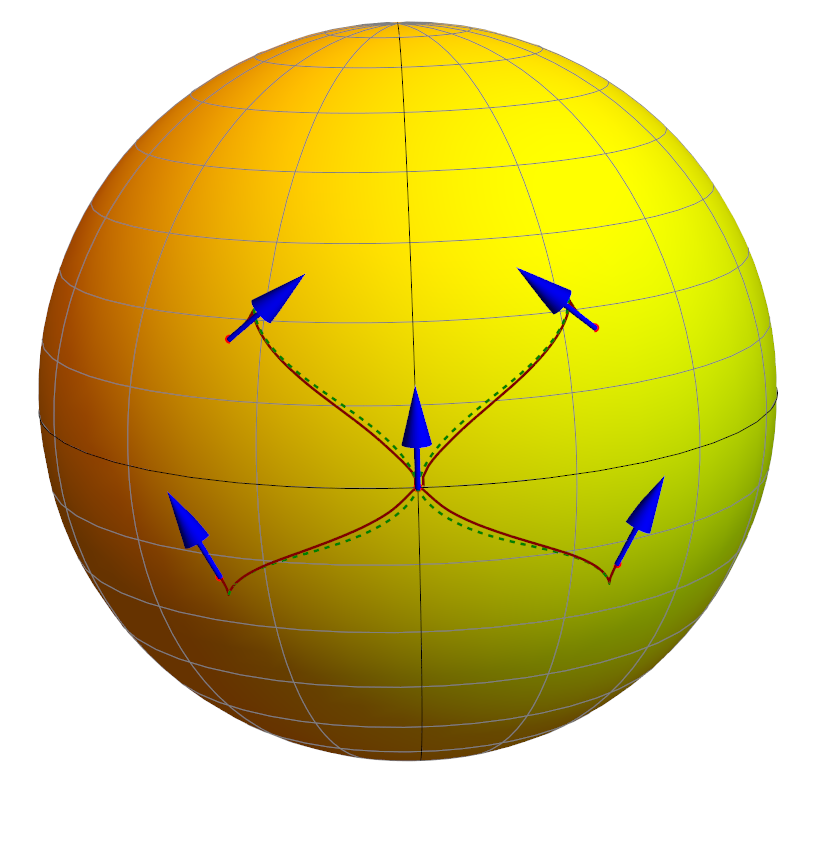}\\
\includegraphics[width=0.65\linewidth]{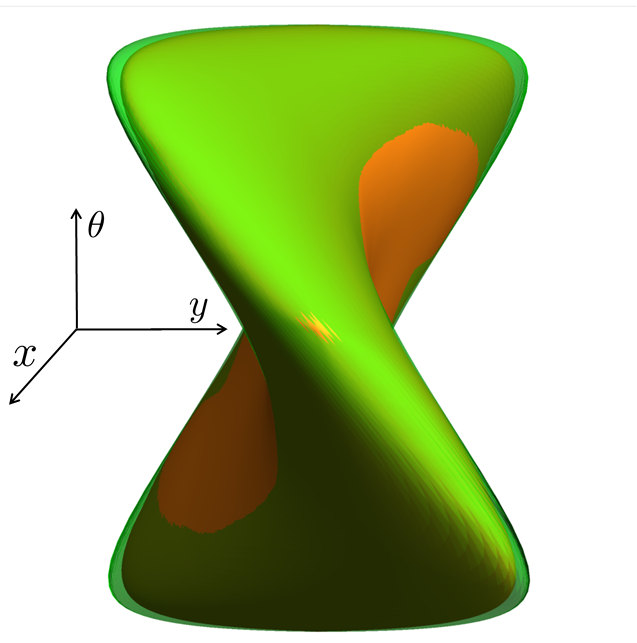}
\caption{\textbf{Top: } Comparison of exact SR geodesics in $\SO$ (green dashed lines) obtained via analytic formulas (Theorem~\ref{thmGeogTpar}) and SR geodesics computed numerically via SR-FM (red lines). Here $\xi =4.5$, $t_{end} = \frac{3\pi}{2}$, and initial momenta $(h^i_2(0),h^i_3(0)) = \{(-5,-0.99),(-5,0.99),(5,-0.99),(5,0.99)\}$.
\textbf{Bottom:} Comparison for $\xi=1$ of wavefront in $\SO$ (green transparent surface) obtained via the analytic formulas of Theorem~\ref{thmGeogTpar} for $T = \frac{15}{32}\pi$ and the SR-sphere (orange solid surface) for the same $T$, obtained via SR-FM. }
\label{fig:CompWithCusp}
\end{figure}

To verify the SR-FM method we also perform experiments with comparison of the geodesics obtained via SR-FM and the general geodesics (not necessarily cuspless) given by Theorem~\ref{thmGeogTpar}. A typical result is presented in Fig.~\!\ref{fig:CompWithCusp} (top). Here we again observe an accurate result of SR-FM, but now for the geodesics whose spherical projections have cusps.

To conclude this section, we provide one more experiment, where we compare the exact sub-Riemannian wavefronts and SR spheres computed via SR-FM.
We show that the SR-FM method provides a distance function $\mathcal{W}(g)$ that closely approximates the SR distance $d(e,g)\approx \mathcal{W}(g)$. In Fig.~\!\ref{fig:CompWithCusp}, we show the comparison of the exact wavefront $\WF(\frac{15}{32}\pi)$ and isosurface  $\mathcal{S}^{FM}(\frac{15}{32}\pi) = \{g \in \SO|\mathcal{W}(g) = \frac{15}{32}\pi\}$. Here we see that the isosurface computed via SR-FM accurately follows the outer surface of the exact wavefront (i.e., exact SR sphere). A similar picture was obtained for different radii $T$.
\subsection{Comparison of SR Geodesics and Wavefronts in $\SE$ and $\SO$ for $C=1$}\label{subsec:ComparisonSE2SO3_UniformCost}
In this subsection, we compare SR geodesics and wavefronts in $\SE$ and $\SO$ and analyze their applicability in image processing. We also include a short discussion about optimality of geodesics, which is closely related to the analysis of self-intersections of wavefronts and the study of cut loci.
To this end, we provide accurate plot of the wavefronts near their singularities.
\subsubsection{Comparison of SR Wavefronts in $\SO$ and $\SE$ for $C=1$}\label{subsec:ComparisonSE2SO3geodesics}
In this subsection, we show a comparison of the SR wavefronts in $\SO$ and $\SE$ in the case of uniform external cost $C=1$.
 Here we employ the fact that the coordinate chart $(x,y,\th)$ in $S^2_{x,y} \times S^1_{\th}$ was chosen in the way to obtain the analogy with $(X,Y,\Theta) \in \R^2 \times S^1$. This allows us to plot SR wavefronts of $\SO$ and $\SE$ in the same 3D plot.

 In this comparison, we use the SR arclength parameterization $t$, where the geodesic $\gamma^{\SO}(\cdot)$ is given by exponential map $\Exp(h^0,\cdot)$, recall (\ref{eq:expmap}), i.e., by the projection in $\SO$ of the solution of the Hamiltonian system~(\ref{eq:hamsys}), with initial condition $h(0) = h^0$, $\gamma^{\SO}(0) = e = (0,0,0)$.
 To establish the comparison of wavefronts we switch to polar coordinates $(\beta, c)$, recall~(\ref{eq:polarcoords}), where the Hamiltonian system in $\SO$ reads as:
 \begin{equation*} \label{eq:hamsysSO3inpolar}
\begin{array}{ll}
\begin{cases}
\dot{\beta} = c,\\
\dot{c} = \frac12 \frac{\xi^2-1}{\xi^2} \sin \left(2 \beta\right), \\
\end{cases} 
& \begin{cases}
\dot{x} = \frac{1}{\xi}\cos \beta \cos\theta, \\
\dot{y} = -\frac{1}{\xi} \cos\beta \sec x \sin \th, \\
\dot{\th} = \sin\beta + \frac{1}{\xi} \cos\beta \tan x \sin\th,
\end{cases} \\ 
 \text{--- vertical part,}
&  \text{--- horizontal part.}
\end{array}
\end{equation*}
A solution to this system is given by Theorem~\ref{thmGeogTpar}.

The Hamiltonian system for SR geodesics in $\SE$ (see, e.g.,~\cite{max_sre}
) reads as:
\begin{equation*} \label{eq:hamsysSE2}
\begin{array}{ll}
\begin{cases}
\dot{\beta} = c,\\
\dot{c} = \frac12 \sin \left(2 \beta\right), \\
\end{cases} 
& \begin{cases}
\dot{X} = \frac{1}{\xi}\cos \beta \cos\Theta, \\
\dot{Y} = -\frac{1}{\xi} \cos\beta \sin\Theta, \\
\dot{\Theta} = \sin\beta,
\end{cases} \\ 
 \text{--- vertical part,}
&  \text{--- horizontal part.}
\end{array}
\end{equation*}
A solution to this system is given in~\cite{max_sre}. Indeed, we observe a clear similarity between the $\SO$ and $\SE$ case, when using parametrization~(\ref{matexp}).

\begin{figure}[ht]
\centering
\includegraphics[width=0.45\linewidth]{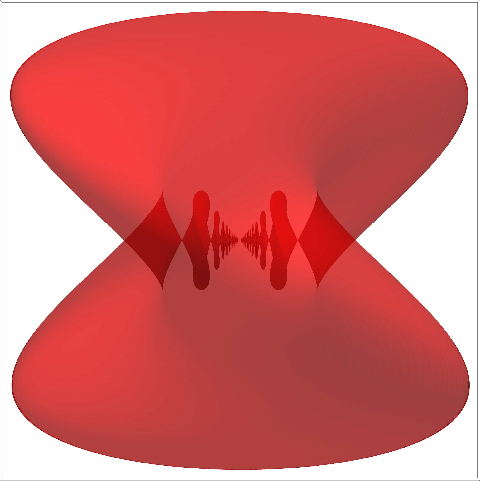}
~
\includegraphics[width=0.45\linewidth]{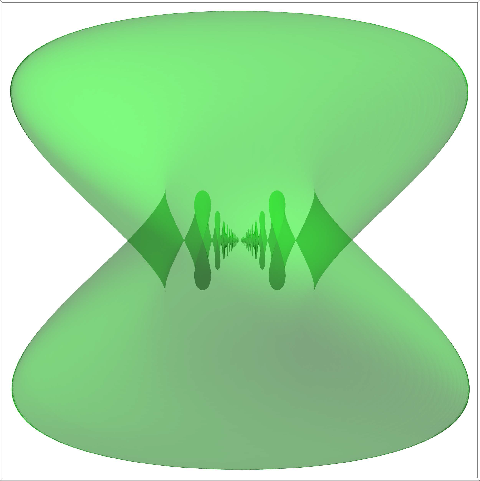}
\caption{Comparison of SR wavefronts in $\SE$ (red) and $\SO$ (green) for $T=0.5$, $\xi = 1$, $\eta = 1$. A zoomed in picture for the blue square is provided in Fig.~\ref{fig:WFSingcomparison}}
\label{fig:WFcomparison}
\end{figure}

In Fig.~\!\ref{fig:WFcomparison}, we plot the SR wavefronts in $\SO$ and $\SE$ for several values of end time $T$. In these plots, we identify $(x,y,\th)$ and $(X,Y,\Theta)$, so that the red surfaces $\WF^{\SE}$ are the SR wavefronts of $\SE$ and the green surfaces $\WF^{\SO}$ are the SR wavefronts of $\SO$. We see a very similar shape of $\WF^{\SO}$ and $\WF^{\SE}$ for small radii $T$, but the difference increases when $T$ increases. Thus, we conclude that the SR geodesics in $\SO$ can be locally approximated by the SR geodesics in $\SE$, but globally they are considerably different.

\begin{figure*}[ht]
\centering
\includegraphics[width=\linewidth]{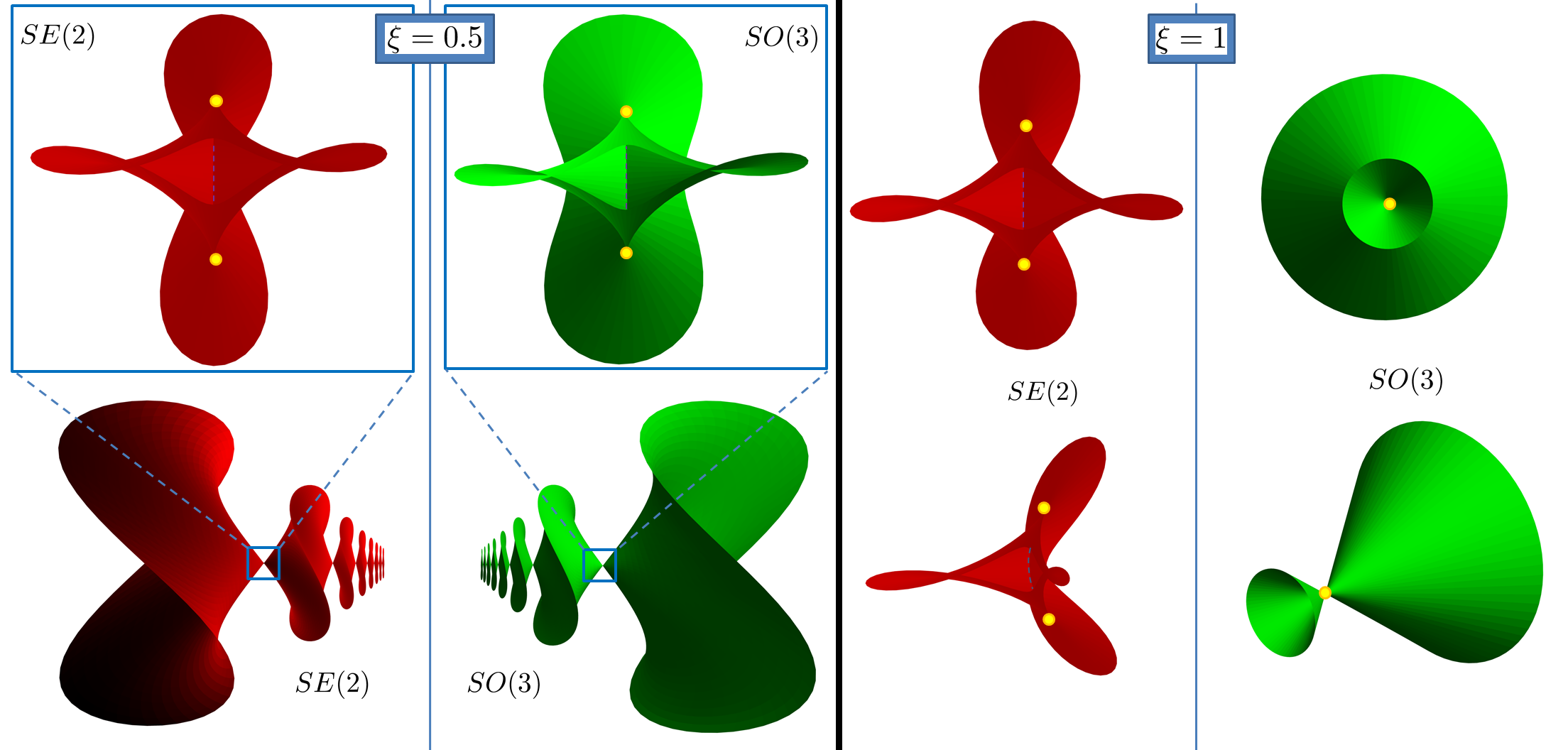}
\caption{Comparison of selfintersections of SR wavefronts in $\SE$ (red) and $\SO$ (green) for $\xi=1$ (linear case) and $\xi = 0.5$ (elliptic case). Here $T = 0.5$ and  $\eta = 1$. The viewpoint is taken from the inside of the SR sphere. In this figure we plot the part of wavefront depicted by the blue square in Fig.~\!\ref{fig:WFcomparison}. If $\xi \neq 1$ both the wavefronts in $\SE$ and in $\SO$ do not intersect at a single point. The first Maxwell sets are depicted by dashed violet lines. The first conjugate points on the SR sphere are depicted by yellow dots. }
\label{fig:WFSingcomparison}
\end{figure*}

One can also observe that the singular points\footnote{By singular points, we mean either conjugate or Maxwell points (recall Remark~\ref{remark:MaxwellConj}).} of $\WF^{\SO}$ are located near singular points of $\WF^{\SE}$. This leads us to a conjecture that the location of conjugate points 
(open problem) can be estimated by the location of conjugate points of $\WF^{\SE}$. Note also that in the general case $\xi>0$ the set of singular points has a complicated shape. In Fig.~\!\ref{fig:WFSingcomparison} we present more detailed plot of the singularities, with the depicted Maxwell set 
and special cases of conjugate points, which are limit points of the Maxwell set on SR sphere (outer part of wavefront).

In Fig.~\!\ref{fig:WFSingcomparison} we observe that the wavefront in $\SO$ has a very special symmetry when $\xi=1$. This is not present for $\xi\neq 1$ and this never happens in the $\SE$ group. The case $\xi = 1$ for SR manifold in $\SO$ was completely examined in~\cite{Boscain_Rossi_B}, where it was shown that locally the conjugate locus is an interval, and globally it is a circle without a point. Changing $\xi$ destroys the symmetry, conjugate and Maxwell points are getting separated, and the conjugate locus has an astroidal shape~\cite{Agrachev_exp_map}.
\subsubsection{Comparison of SR Geodesics in $\SO$ and $\SE$ for $C=1$}\label{subsec:ComparisonGeodesics}
In this subsection, we again consider the case $C=1$ and compare SR geodesics $\gamma^{\SO}(\cdot) = (x(\cdot),y(\cdot), \th(\cdot))$ and $\gamma^{\SE}(\cdot) = (X(\cdot),Y(\cdot), \Theta(\cdot))$ in the image plane. The SR-FM method is used for computation of the geodesics parameterized by SR arclength. Here we prepare background for comparison of the geodesics in retinal images via the schematic eye model, recall Sect.~\ref{subsec:Motivation1}, where as a departure point we use an image (white for $C=1$) on a plane $O_{XY}$, recall Fig.~\!\ref{fig:Twospheres}.

\begin{figure}[ht]
\centering
\includegraphics[width=\linewidth]{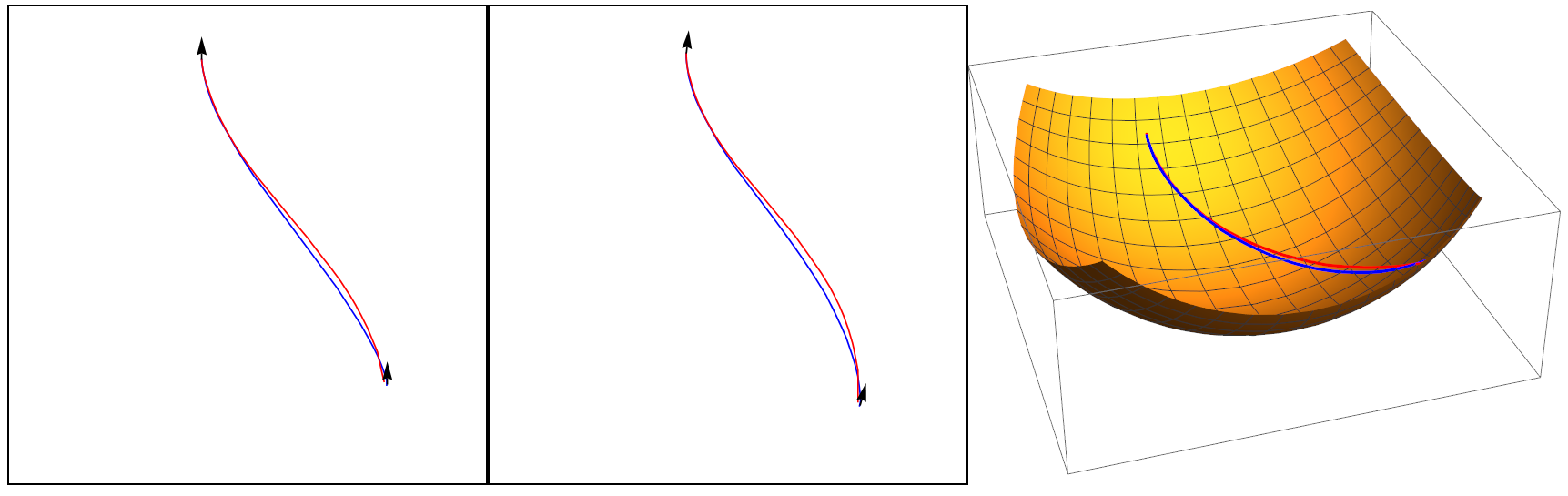}
\includegraphics[width=\linewidth]{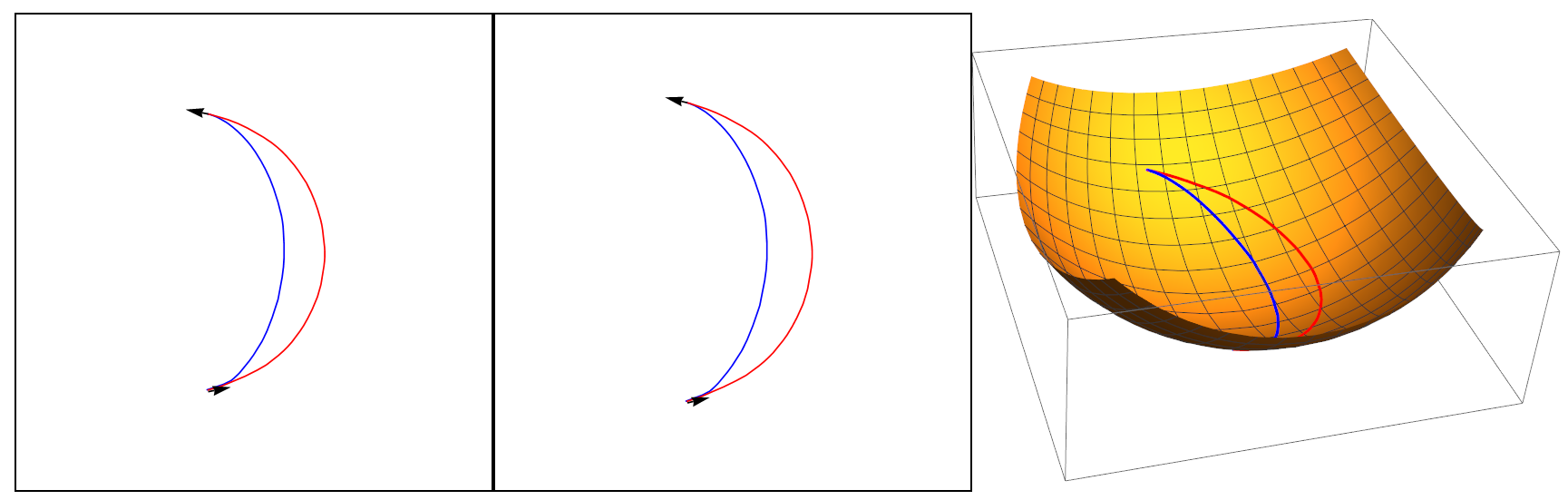}
\caption{Comparison of an $\SE$ SR geodesic (blue) and an $\SO$ SR geodesic (red), from left to right in $(x,y)$ spherical image coordinates, $(X,Y)$ flat image coordinates, and plotted on the sphere $S^2$.}
\label{fig:ExperimentA}
\end{figure}

See Fig.~\!\ref{fig:ExperimentA}, where we compare $\SE$ and $\SO$ SR geodesics projected on the plane and on the sphere (via mappings $\Pi$ and $\Pi^{-1}$). For details, see Appendix~\ref{app:geodcomp}. 
\subsection{Vessel Analysis via $\SO$ SR Geometry and $\SE$ SR Geometry}\label{subsec:Retinal_Example}
As explained in introduction in Sect.~\ref{subsec:Motivation1}, we need to include the spherical geometry of the retina rather than the flat geometry of the flat image. This spherical geometry is encoded in our spherical image model, see Fig.~\!\ref{fig:Twospheres}. Next we will analyze the effect of including this geometry in the SR-FM vessel tracking method along data-driven SR geodesics in $\SO$.

\begin{figure*}[ht]
\centering
\includegraphics[width=0.7\linewidth]{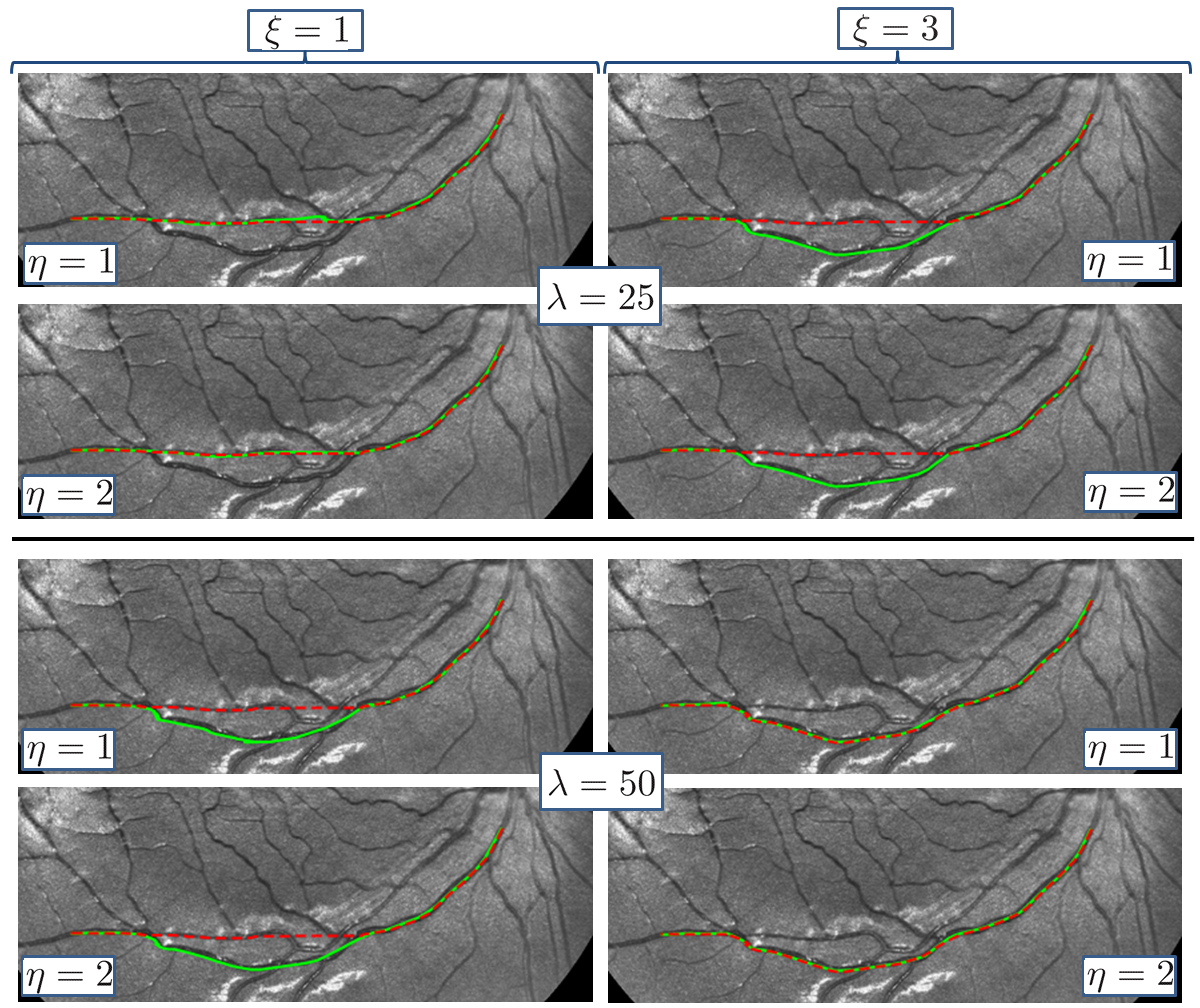}
\caption{Comparison of vessel tracking via $\gamma^{\SO}$ SR geodesics in $\SO$ (green solid lines) in object coordinates and via $\gamma^{\SE}$ SR geodesics in $\SE$ (red dashed lines) in the planar camera coordinates. Here the planar projection $\Gamma^{\SO}$ of $\gamma^{\SO}$,  and spatial projection $\Gamma^{\SE}$ of $\gamma^{\SE}$ are depicted in the same flat image.}
\label{fig:VesselComparison}
\end{figure*}

\begin{figure*}[ht]
\centering
\includegraphics[width=\linewidth]{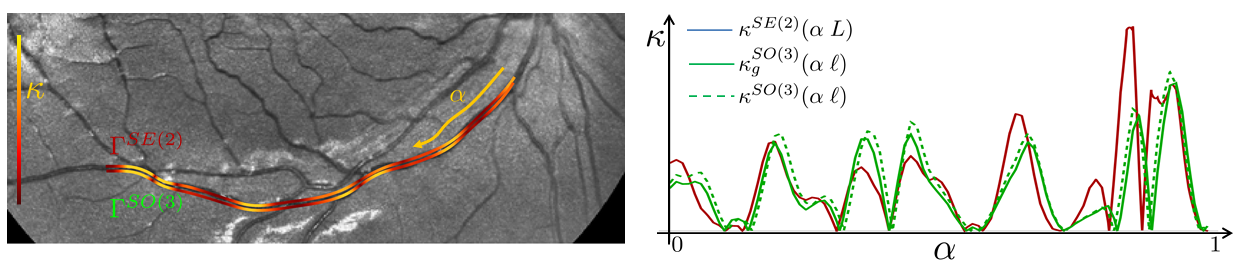}
\caption{\textbf{Left:} Two curves from the experiment in Fig.~\ref{fig:VesselComparison} (right-bottom figure) are depicted with slight shift. The upper curve is a spatial projection $\Gamma^{\SE}$  of the data-driven SR geodesic $\gamma^{\SE}$  with depicted (in color) planar curvature $\kappa^{\SE}$ on top of the curve.
The lower curve $\Gamma^{\SO}$ is the planar projection of the SR geodesic $\gamma^{\SO}$ with depicted geodesic curvature $\kappa^{\SO}_g$ on top of it.
\textbf{Right:} Three graphs are shown in the same plot: planar curvature $\kappa^{\SE}$ of $\Gamma^{\SE}$; geodesic curvature $\kappa_g^{\SO}$ of spherical projection of $\gamma^{\SO}$; planar curvature $\kappa^{\SO}$ of $\Gamma^{\SO}$.
The effect of considering geodesic curvature $\kappa_g^{\SO}$ in object coordinates on $S^2$ rather than planar curvature $\kappa^{\SO}$ in photo coordinates on projection on $\R^2$ is visible (compare the green solid and dashed graphs). A bigger difference comes from using $\SO$ SR-geometry than $\SE$ SR-geometry (compare red and green graphs).}
\label{fig:CurvatureFigure}
\end{figure*}

More precisely, we propose vessel tracking in object coordinates (or spherical image coordinates) via SR geometry in $\SO$ as an extension of vessel tracking in flat images~\cite{SSVM,SIAM} along SR geodesics in $\SE$. Therefore, we want to investigate whether including the correct spherical geometry makes a difference in the vessel tracking in practice. Although a complete detailed comparison on large data sets is left for future work, we present preliminary experiments which indeed indicate considerable differences in both tractography and curvature measurements. These experiments are shown in Figs.~\!\ref{fig:VesselComparison}, and~\!\ref{fig:CurvatureFigure} and next we explain them.

We apply the same scheme as in Sect.~\ref{subsec:ComparisonGeodesics}, but now we compute data-driven geodesics, where the external cost is induced by image data. Next we explain the construction of the external cost that is applied in all experiments. For the sake of simple comparison, we restrict ourselves to a cost depending on the spherical coordinates only, and we set
\begin{equation}
\label{eq:externalcostVessel}
\gothic{G}(x,y) = \left(1+ \frac{VF(\Pi(x,y))}{\lambda\|VF\|^{2}_{\infty}}\right)^{-1},
\end{equation}
where we use standard multiscale vesselness~\cite{Frangi}
\begin{multline*}
(VF)(X,Y) = \\
e^{-\frac{\lambda_1^2(X,Y)}{2 \beta^2 \lambda_2^2(X,Y)}}
 \left(1 - e^{-\frac{\lambda_1^2(X,Y) + \lambda_2^2(X,Y)}{2 c^2}} \right) U(\lambda_2(X,Y)),
\end{multline*}
with $\lambda_i(X,Y)$ eigenvalues of the Gaussian Hessian of image $F:\R^2 \to \R$ (maximized over scales) ordered by $|\lambda_1(X,Y)|<|\lambda_2(X,Y)|$, with $\beta = c = 0.3$, and  with unit step function $U(\lambda) = \begin{cases}1, \text{ for } \lambda \geq 0,\\ 0, \text{ otherwise.} \end{cases}$\\
Here we note that the Gaussian Hessian is given by $\mathcal{H}(G_s \ast F)$ and computed by Gaussian derivatives~\cite{BartBook} at multiscales $s = \frac12 \sigma^2 \in \{2,3,4,5\}$ in term of pixel sizes.

In the experiment in Figure~\ref{fig:VesselComparison}, we show that there is a considerable difference between $\SE$ SR geodesics and $\SO$ SR geodesics. We see that	when internal geometry is dominant over the external cost ($\lambda$ small) the SR geodesics in $\SO$ are more stiff than SR geodesics in $\SE$, and therefore in the boundary value problem they are less eager to take shortcuts and better follow the vessel structure. In case $\lambda$ is large (external cost is dominant than the internal geometry), we see only small differences in the overall locations of the $\SE$ curves and $\SO$ curves. The results are stable w.r.t. choice of $1\leq\eta\leq2$ (which controls the distance from the camera to the eye ball, relative to eye ball radius, recall Fig.~\ref{fig:Twospheres}).

In the next experiment, we measure the curvature of the curves obtained by the vessel tracking method via $\SE$ geometry and via $\SO$ geometry. For this experiment, we used the values $\xi=3$, $\lambda= 50$ and $\eta=2$. Although in this case the result of tractography is very similar for the $\SE$ and $\SO$ curves, we show that there is a considerable difference in their curvature.
\begin{corollary}(from Theorem~\ref{prop:31})
The geodesic curvature of a spherical projection of data-driven geodesic $\gamma^{\SO}(\cdot)$ satisfies $$\kappa_g^{\SO}(\cdot) = \xi^2 \frac{X_2|_{\gamma^{\SO}(\cdot)}(\mathcal{W^{\SO}})}{X_1|_{\gamma^{\SO}(\cdot)}(\mathcal{W^{\SO}})}.$$
\end{corollary}
\textbf{Proof }
The first equality in the chain
$$\kappa_g^{\SO}(\cdot) = \frac{u_2(\cdot)}{u_1(\cdot)} = \xi^2\frac{h_2(\cdot)}{h_1(\cdot)} = \xi^2 \frac{X_2|_{\gamma(\cdot)}\left(\mathcal{W}\right)}{X_1|_{\gamma(\cdot)}(\mathcal{W})}$$
is implied by Eq.~(\ref{Pcurvecontrols}). The second equality follows from application of PMP to problem $\Pmec$, which gives $u_1 = \frac{h_1}{ \xi^2 C^2}$ and $u_2 = \frac{h_2}{C^2}$. The third equality follows from the fact that in the points where the Bellman function  $\mathcal{W}$ is differentiable (almost everywhere in our case) its derivatives are given by components of momentum covector
$h_1 = X_1(\mathcal{W})$,  $h_2 = X_2(\mathcal{W})$, see~\cite{agrachev_sachkov}.
$\hfill \Box$

By the same argument, it can be checked (see~\cite{DuitsJMIV2014} and~\cite{SIAM}) that the planar curvature of spatial projections of SR geodesics in $\SE$ satisfies
$$\kappa^{\SE}(\cdot) = \xi^2 \frac{\mathcal{A}_2|_{\gamma^{\SE}(\cdot)}(\mathcal{W^{\SE}})}{\mathcal{A}_1|_{\gamma^{\SE}(\cdot)}(\mathcal{W^{\SE}})},$$
with $\mathcal{A}_1 = \cos\Theta \partial_X + \sin\Theta \partial_Y$, $\mathcal{A}_2 = \partial_\Theta$ basis left-invariant vector fields $\SE$.

Thus the curvature analysis can be simply done based on vessel tracking, and this shows the benefit of our algorithm. In Fig.~\!\ref{fig:CurvatureFigure}, we show an experiment of vessel curvature measurement based on tracking via $\SO$ geometry  and via $\SE$ geometry. For completeness, we added also a comparison with a planar curvature $\kappa^{\SO}(\cdot)$ of a planar projection $\Gamma^{\SO}(\cdot) :=  \Pi(x(\cdot), y(\cdot))$ of $\gamma^{\SO}(\cdot) = (x(\cdot),y(\cdot),\th(\cdot))$.

We can see a considerable difference in curvature measurement via $\SO$ geometry and $\SE$ geometry. It is also seen, that that the difference between $\kappa^{\SO}$ and $\kappa_g^{\SO}$ is not very significant. Thus we see the importance of including correct spherical geometry in vessel tracking algorithm in retinal images.

\subsection{Vessel Analysis via $\SO$ SR Geometry and $\SO$ Riemannian Geometry}\label{subsec:RiemSRcompar}
\begin{figure}[ht]
\centering
\includegraphics[angle=270,origin=c,width=0.8\linewidth]{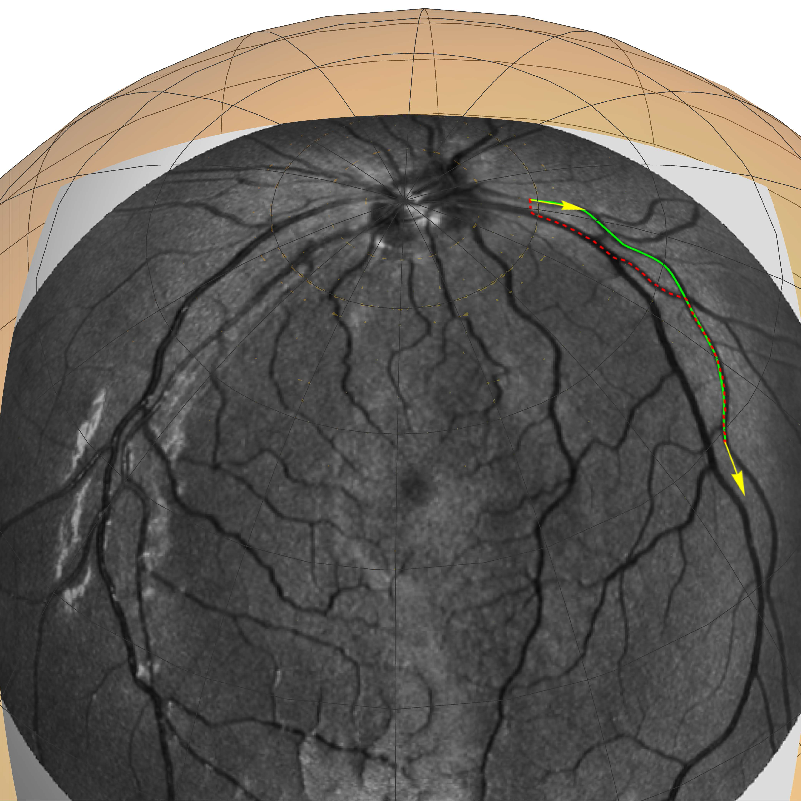}
\caption{Comparison of vessel tracking in a spherical image of the retina via SR geodesic (green solid line) and Riemannian geodesic (red dashed line) in $\SO$.}
\label{fig:SRandRiem}
\end{figure}

Next we address the general benefit of using SR geodesics rather than (isotropic) Riemannian geodesics in the tracking of salient lines (blood vessels) in spherical images. To this end, we show a typical example where the tracking induced by a sub-Riemannian geodesic gives better result than the tracking induced by a Riemannian geodesic. The experiment in Fig.~\ref{fig:SRandRiem} is performed similar to Sect.~\ref{subsec:Retinal_Example},  but now we compare the result of vessel tracking via SR geodesics (obtained by SR-FM) and Riemannian geodesics (isotropic metric with $\varepsilon = 1$). From this experiment, we see that similarly to the $\SE$ case~\cite{SIAM} the tracking via Riemannian geometry suffers from incorrect jumps toward nearly parallel neighboring vessels and yields nonsmooth curves, the tracking via SR geometry gives the desirable result. 

\section{Conclusion}
Data-driven sub-Riemannian geodesics in 3D Lie groups are a suitable tool for tractography of blood vessels in retinal imaging. In previous works on the $\SE$ case~\cite{SIAM,CIARP} practical advantages have been shown in comparison with the (isotropic) Riemannian case, and geodesic methods in the image domain.
However, these models included a SR geometry on $\SE \equiv \R^{2}\rtimes S^{1}$ based on lifts of flat images, which does not match the actual object geometry: The retina is spherical rather than planar, cf.~Fig.~\!\ref{fig:RetinaSphere}, Fig.~\!\ref{fig:Eye2}, and Fig.~\!\ref{fig:Twospheres}. A similar observation holds for models in the psychology of contour perception in human vision~\cite{Boscain2014}.

Therefore, for geometric tracking we propose a frame bundle above $S^{2}$, cf.~Fig.~\!\ref{fig:x}, instead of a frame bundle above $\R^{2}$. Geometric tracking of geodesics is done along globally optimal data-driven SR geodesics in $\SO$ (and their spherical projections) by our new numerical wavefront propagation method. The method was validated for the uniform cost case by comparisons with exact geodesics which we derived in Sect.~\ref{sec:PMEC}, cf.~Figs.~\!\ref{fig:ExactGeodVsNumeric}, and~\ref{fig:CompWithCusp}. Here, in contrast to the $\SE$ case \cite{DuitsJMIV2014}, we do not have a scaling homothety. As a result, the parameter $\xi$ has a considerable effect on the geometry, and for $\xi=1$ we identify the linear case, for $\xi>1$ we identify the hyperbolic case, and for $\xi<1$ we identify the elliptic case, cf.~Figs.~\!\ref{fig:IntegralsVP} and~\!\ref{fig:extremalcurveE}.

For all of these cases, we have computed the first cusp time in Theorem~\ref{thmSmax} (the first time where the spherical projection of a SR geodesic exhibits a cusp). Also, we have presented new formulas for such ``cuspless'' SR geodesics in Theorem~\ref{thmGeogSpar}. These formulas only involve a single elliptic integral thanks to spherical arclength parametrization, and for $\xi \neq 1$ our formulas are simpler than the general formulas for SR geodesics in $\SO$.

Furthermore, we used a specific parametrization of Lie group $\SO$ that allowed us to compare between SR geodesics / wavefronts in $\SO$ to SR geodesics / wavefronts in $\SE$, cf.~Figs.~\!\ref{fig:WFcomparison} and~\ref{fig:WFSingcomparison}. In our comparison we took into account a standard optical model for the mapping between object coordinates on the retina and camera coordinates in the acquired planar retinal image. In our experiments, the differences between the $\SO$ case and the $\SE$ case are considerable, both for the case of uniform cost, cf.~Fig.~\!\ref{fig:ExperimentA}, and for the data-driven case in the retinal image analysis application, cf.~Fig.~\!\ref{fig:VesselComparison}. In general, we see that for realistic parameter settings (in optics) the $\SO$ geodesics have a slower variation in curvature and are less eager to take shortcuts, see, e.g., Figs.~\!\ref{fig:VesselComparison} and~\!\ref{fig:CurvatureFigure}. Furthermore, there are visible differences between geodesic curvature of data-driven SR geodesics on the sphere and the curvature of their planar projections.
As in retinal imaging applications, curvature is considered as a relevant biomarker \cite{Erik_Curv,Kalitzeos,Sasongko} for detection of diabetic retinopathy and other systemic diseases, the data-driven SR geodesic model in $\SO$ is a relevant extension of our data-driven geodesic model in $\SE$. Here we restricted ourselves to feasibility studies. More extensive comparisons between $\SO$ geodesics and $\SE$ geodesics on large retinal imaging benchmark sequences are beyond the scope of this article and are left for future work.

Finally, we note that the computation time for data-driven SR geodesics in $\SO$ is exactly the same as for the $\SE$ case. Our specific choice of coordinates of $\SO$ allowed us to modify the very efficient fast marching approach~\cite{CIARP}, with a simple replacement of the metric tensor matrix.
\section*{Acknowledgments}
The research leading to these results has received funding from the European
Research Council under the European Community's Seventh Framework Programme
(FP7/2007-2013)/ERC grant \emph{Lie Analysis}, a.n.~335555.
The authors gratefully acknowledge EU-Marie Curie project \emph{MANET}
ag.~no.~607643,
and RFBR research Project No. 16-31-60083 mol\_a\_dk
for financial support. 
The publication was supported by the Ministry of Education and Science
of the Russian Federation (the Agreement number 02.a03.21.0008).
The authors gratefully acknowledge Dr. G.R. Sanguinetti and Dr. J-M. Mirebeau on their efforts on anisotropic fast marching implementation of SR geodesics in $\SE$ \cite{CIARP} that allowed for our generalization to the $\SO$ case. Furthermore, the authors thank G.R. Sanguinetti on
helpful comments and suggestions on earlier versions of this article.
We thank Dr. Tos Berendschot for fruitful discussions on retinal imaging and for his encouragement to take into account the spherical nature of the retina in our geometric retinal image processing. The authors thank anonymous reviewers whose valuable comments improved the exposition of this work.
\appendix
\section{Proof of Theorem~\ref{prop:31}} \label{app:C}
The function $s(t) := \int_{0}^{t} u_1 (\tau) {\rm d} \, \tau$ is differentiable and increasing for $t \in [0,T]$; then there exists an inverse function $t(s)$, $s \in [0,l]$, $l = s(T)$, differentiable and increasing. Notice that $\frac{ds}{dt}(t) = u_1(t)$, which gives the second equality in (\ref{Pcurvecontrols}).

Define $\bar{R}(s):=R(t(s))$, $\ul{n}(s) := \bar{R}(s) \ul{e}_1$, $\underline{\ul{n}}(t):=R(t) \ul{e}_1 = \ul{n}(s(t))$.
From the dynamics of $\Pmec$ we obtain identities
$$
\Scale[0.89]{
\begin{array}{l}
\ul{n}'(s) \, =  \dot{\underline{\ul{n}}}(t(s)) \frac{dt}{ds}(s) =
-\bar{R}(s) A_2 \ul{e}_{1} + \frac{u_2(t(s))}{u_1(t(s))}\bar{R}(s) A_1 \ul{e}_1 =  \bar{R}(s) \ul{e}_3,\\
\ul{n}''(s) \!=\!  -\bar{R}(s) A_2 \ul{e}_{3} \!+\! \frac{u_2(t(s))}{u_1(t(s))}\bar{R}(s) A_1 \ul{e}_3 \!=\!\! -\bar{R}(s) \ul{e}_1 \!-\! \frac{u_2(t(s))}{u_1(t(s))}\bar{R}(s) \ul{e}_2.
\end{array}
}
$$
Then the Gauss--Bonnet formula~(\ref{eq:GaussBonnet}) gives
$$
\Scale[0.9]{
k_g(s)\!=\!\left(\!-\bar{R}(s) \ul{e}_1 - \frac{u_2(t(s))}{u_1(t(s))}\bar{R}(s) \ul{e}_2 \right)\cdot \left(\bar{R}(s) \ul{e}_1 \times \bar{R}(s) \ul{e}_3\right)\!=\!\frac{u_2(t(s))}{u_1(t(s))},
}
$$
which implies the third equality in (\ref{Pcurvecontrols}).

The boundary conditions~(\ref{pcurveRbounds}) for the curve $\ul{n}(s)$ follow from its definition and the boundary conditions for $R(t)$.

Recall that by definition $C(R(t)) = \gothic{C}(\ul{n}(s(t)))$. Since the minimizer $R(t)$ is parameterized by SR arclength, we have
$C^2(R(t))\left(\xi^2 u_1^2(t) + u_2^2(t)\right)\equiv 1,$
whence $\frac{d t}{d s}(s) = \frac{1}{u_1(t(s))} = C(R(t(s)))\sqrt{\xi^2 + \frac{u_2^2(t(s))}{u_1^2(t(s))}} = \gothic{C}(\ul{n}(s))\sqrt{\xi^2 + k_g^2(s)},$
which implies~(\ref{eq:ts}) after integration w.r.t. $s$.

Thus, we can see that the optimization functionals of $\Pcurve$ and $\Pmec$ coincide:
$$
\Scale[0.89]{
\int_0^T C(R(t))\sqrt{\xi^2 u_1^2(t)+ u_2^2(t)} {\rm d}\, t = \int_0^l \gothic{C}(\ul{n}(s))\sqrt{\xi^2 + k_g^2(s)} {\rm d}\, s.
}
$$
But also the dynamics coincides, in the following sense. If $\tilde{\ul{n}}(s)$, $s \in [0,\tilde{l}]$, is a smooth curve on the sphere $S^2$ with the initial conditions $\tilde{\ul{n}}(0) = \ul{e}_1$, $\tilde{\ul{n}}'(0) = \ul{e}_3$ and a geodesic curvature $\tilde{k}_g(s)$, then it can be lifted to a curve $\tilde{R}(s)$ in $\SO$ that is a trajectory of control system~(\ref{eq:so3sys}) with the controls $\tilde{u}_1(s) \equiv 1$, $\tilde{u}_2(s) = \tilde{k}_g(s)$ and the initial condition $\tilde{R}(0) = \Id$. The curve $\tilde{R}(s) \ul{e}_1$ on $S^2$ has the same initial conditions and geodesic curvature as $\tilde{\ul{n}}(s)$, thus $\tilde{R}(s) \ul{e}_1 \equiv \tilde{\ul{n}}(s)$.

Summing up, since $R(t)$ is a minimizer of $\Pmec$, then its projection $\ul{n}(s)$ to the sphere $S^2$ is a minimizer of $\Pcurve$.
\section{Integration of the Hamiltonian System $\PmecOne$}\label{app:B}
To integrate the Hamiltonian system~(\ref{eq:hamsys}) we follow the idea of V.~Jurdjevic~\cite{Jurd}, where one can employ left invariance and introduce initial rotation $D_0$ of momentum space, such that the initial momentum transforms to $h(0) = (0,0,M)$, solve the problem in this simple case (i.e., to find the trajectory $\tilde{R}(t)$ that corresponds to $h(0) = (0,0,M)$), and obtain general solution $R(t)$ by a backward transformation $R(t) = D_0^{-1} \tilde{R}(t)$.

\begin{sketchproof}
Substitution of~(\ref{eq:extremalcontrols}) to~(\ref{eq:controlSystem})  gives
\begin{equation*}\label{horizMat}
\Scale[0.97]{
\dot{R} = R \Omega, \quad \Omega =\! \begin{pmatrix}
0& 0 & -h_1 \xi^{-2} \\
0&0 & -h_2\\
h_1 \xi^{-2} &h_2 &0 \\
\end{pmatrix}\! = -\frac{h_1(t)}{\xi^2} A_2+h_2(t)A_1.
}
\end{equation*}
The Hamiltonian system can be written as a Lax type of system~\cite{Ivan}
\begin{eqnarray}
&&\dot{P} = [P,\Omega] \text{ --- vertical part,} \label{vertpartlax}\\
&&\dot{R} = R \Omega \text{--- horizontal part} \label{horpartlax}.
\end{eqnarray}
Here $P(t) = - h_1(t) \omega^2 +h_2(t) \omega^1  + h_3(t) \omega^3$ is a momentum covector expressed in the basis $\omega^i$, dual to $A_i$, i.e., $\langle \omega^i, A_j\rangle = \delta_{ij}$.
The vertical part (\ref{vertpartlax}) has $P(t) = R^{-1}(t) P(0) R(t)$ as a solution.
Thus, the Hamiltonian system (\ref{vertpartlax})---(\ref{horpartlax}) preserves the norm of momentum covector $h(t) \!=\! (h_2(t), -h_1(t),  h_3(t))$, i.e. the value $M^2 = h_1^2(t) + h_2^2(t) + h_3^2(t)$, 
whose isosurface is a coadjoint orbit of $\SO$. The Killing form allows to identify $\so$ and $\so^{*}$ (see~\!\cite{Ivan}), and by isomorphism~(\ref{so3r3iso}), we write
\begin{equation}
\label{pendvecSol}
P(t) = R^{-1}(t) P(0) R(t) \sim h(t) =  h(0)\, R(t).
\end{equation}
Next we apply a left action of $\SO$ on solution curve $R(\cdot)$
\begin{equation}
\label{RtildeviaR}
\tilde{R}(\cdot)= D_{0} R(\cdot),
\end{equation}
where $D_{0}$ is chosen such that the initial momentum matrix $P(0)$ and initial covector $h(0)=(-h_{2}(0),h_{1}(0),h_{3}(0))$ are reduced to a simple form
\begin{equation}\label{eq:P0}
P(0)= M A_{3} \textrm{ and } h(0)=(0,0,M),
\end{equation}
as the nicest possible representant within the coadjoint orbit of $\SO$.

Next we represent the rotation matrix $\tilde{R}$ in form~(\ref{matexp}), i.e., $\tilde{R} = \eexp{\tilde{y} A_3}\eexp{-\tilde{x} A_2} \eexp{\tilde{\th} A_1}$, where we parameterize the rotations via 3 angles $\ds \tilde{x} \in [-\frac{\pi}{2}, \frac{\pi}{2}]$, $\tilde{y}\in  \R/\{2 \pi \Z\}$, $\tilde{\theta} \in \R/\{2 \pi \Z\}$.

Substitution of (\ref{eq:P0}) in (\ref{pendvecSol}) gives the momentum matrix
\begin{eqnarray*}
\Scale[0.94]{
P(t) = M \eexp{-\tilde{\th}(t)A_1} \eexp{\tilde{x}(t) A_2} \overbrace{\eexp{-\tilde{y}(t) A_3} A_3 \eexp{\tilde{y}(t) A_3}}^{A_3}\eexp{-\tilde{x}(t) A_2}\eexp{\tilde{\th}(t)A_1},
}
\end{eqnarray*}
and the equivalent relation for the momentum covector
\begin{eqnarray*}
&&(h_2,-h_1,h_3) = \\
&&\qquad \qquad  (0,0,M)
\left(\begin{array}{ccc}
\cos\tilde{x}&	0 &  -\sin\tilde{x}\\
0&	1&	0\\
\sin\tilde{x}&	0 &	\cos\tilde{x}
\end{array}\right)
\left(\begin{array}{ccc}
1&	0 &	0\\
0&	\cos \tilde{\th}&	-\sin\tilde{\th}\\
0&	\sin\tilde{\th}&	\cos\tilde{\th}
\end{array}\right).
\end{eqnarray*}
By multiplying the matrices in the right-hand side, we obtain
\begin{equation}
\label{barh}
\begin{cases}
h_2(t) =  M \sin\tilde{x}(t),\\
- h_1(t) = M \sin\tilde{\th}(t) \cos\tilde{x}(t),\\
h_3(t) = M \cos\tilde{\th}(t) \cos\tilde{x}(t).
\end{cases}
\end{equation}
From~(\ref{barh}) we immediately have
$
\tilde{\th}(t) = \Arg\left(h_3(t) - \i h_1(t)\right),
$
 and since $\ds \tilde{x}(t) \in [-\frac{\pi}{2}, \frac{\pi}{2}] \Rightarrow \cos \tilde{x}(t) \geq 0$, we can also express
$
\tilde{x}(t)= \arg \left(\sqrt{h^2_{1}(t)+h^2_3(t)}+ \i \, h_{2}(t) \right).
$

To obtain $\tilde{y}(t)$ first notice that due to left invariance, the matrix $\tilde{R}$ satisfies the same equation as $R$, i.e., $\dot{\tilde{R}} = \tilde{R} \Omega$. Thus $\tilde{x}$, $\tilde{y}$, $\tilde{\th}$ satisfy the same equations as $x$, $y$, $\th$, i.e., the horizontal part of~(\ref{eq:hamsys}). Thus we have
\begin{eqnarray}
\label{dbarx}
\dot{\tilde{y}}(t) = - \frac{h_1(t)}{\xi^2} \, \sec \tilde{x}(t) \sin \tilde{\th}(t),
\end{eqnarray}
with the initial condition  $\tilde{y}(0)=0$ following from $\tilde{R}(0)=D_{0}$.

By virtue of~(\ref{barh}) we can express~(\ref{dbarx}) as
$$
 \dot{\tilde{y}}(t) =\frac{M}{\xi^2}\left( \frac{h_1^2(t)}{h_1^2(t) + h_3^2(t)}\right) \Leftrightarrow \tilde{y}(t) = \frac{M}{\xi^{2}}\int_0^t \frac{h_1^2(\tau)}{h_1^2(\tau) + h_3^2(\tau)} d \tau.
$$
Further, since $M^2 = h_1^2(t) + h_2^2(t) + h_3^2(t) = \const$, we have
\begin{equation*}
\label{baryf}
\tilde{y}(t) = 
\frac{M}{\xi^2}\left( t -  \int_0^t \frac{h_3^2(\tau)}{M^2 - h_2^2(\tau)} d \tau\right).
\end{equation*}
Thus we obtained solution for $\tilde{x}(t)$, $\tilde{y}(t)$, $\tilde{\th}(t)$. To finish the proof, it only remains to obtain solution for $x(t), y(t), \th(t)$.
Recall that the expression (\ref{eq:RtotildeR}) for $R(t)$ follows from orthogonality of matrix $D_0$, i.e., $D_0^{-1} = D_0^{T}$, and parameterization of $R(t)$ by $x(t), y(t), \th(t)$ as $R(t) = \eexp{y(t) A_3} \eexp{-x(t) A_2} \eexp{\th(t) A_1}$.

To find $x(t), y(t), \th(t)$ notice that action of $R(t)$ on $\ul{e}_1$ and $R^{-1}(t) = R^T(t)$ on $\ul{e}_3$ gives
\begin{eqnarray}
&&(R_{11}(t),R_{21}(t),R_{31}(t))^T = \nonumber \\
&&\qquad \quad (\cos x(t) \cos y(t), \cos x(t) \sin y(t), \sin x(t))^T, \label{xy}\\
&&(R_{31}(t),R_{32}(t),R_{33}(t))^T = \nonumber \\
&&\qquad \quad (\sin x(t), \cos x(t) \sin \th(t), \cos\th(t) \cos x(t))^T \label{xth}.
\end{eqnarray}
Since $x(t) \in [-\frac{\pi}{2},\frac{\pi}{2}]$, first two equations in (\ref{xythsol}) follow from (\ref{xy}) and the third equation in (\ref{xythsol}) follows from (\ref{xth}).
\end{sketchproof}
\section{Proof of Theorem~\ref{thmSmax}}\label{app:smax}
Recall that we consider the case (see Remark \ref{rm1}) where the curve starts from $\ul{e}_1 \in \S^2$ and
goes in the direction of upper half of the sphere with tangent vector $\ul{e}_3$. Therefore, we have $h_1(0) = u_1(0) \geq 0$, $h'_1(0)>0$.

Via the Hamiltonian $H(s) = \frac{\xi^{-2}h_1^2(s) + h_2^2(s)}{2} = \frac12$, we  express $h_1(s) = \xi\sqrt{1- h_2^2(s)}=\sqrt{\chi^2+1} \sqrt{1- h_2^2(s)}$.
As a result, one has
\begin{equation} \label{eq:h2}
h_1(s)=0 \Leftrightarrow h_2(s) = \pm 1.
\end{equation}
In the linear case $\chi=0$, Eq.\!~(\ref{eq:h2}) reads as $h_2^0 + h_3^0 s = \pm 1$, that has the minimal positive root
\begin{eqnarray*}
\label{eq:smaxL}
\smaxL(h_2^0, h_3^0) = \frac{\sgn(h_3^0)-h_2^0}{h_3^0}.
\end{eqnarray*}
In the elliptic and the hyperbolic cases, we need to find a minimal positive root of
\begin{equation} \label{smaxeq}
h_2^0 \cosh s \chi + \frac{h_3^0}{\chi} \sinh s \chi = \pm 1.
\end{equation}
In the elliptic case, we have $\chi = \i a$, where $a = \sqrt{1-\xi^2} \in (0,1)$, 
yielding solution
\begin{eqnarray*}
\label{eq:smaxE}
&&\smaxE(a, h_2^0, h_3^0) =\\
&& \quad \Scale[1.1]{
-\frac{\i}{a}\, \, \log \left(\frac{\sgn \left(h_3^0\right)\left(\sqrt{(h_3^0)^2 - a^2\, \left(1-(h_2^0)^2\right)}+\i a \right)}{h_3^0+\i a \, h_2^0}\right)\!.
}
\end{eqnarray*}
Denote $\kappa=(h_3^0)^2+\left(1-(h_2^0)^2\right) \chi ^2$. Notice that when $\frac{(h_3^0)^2}{a^2} + (h_2^0)^2 < 1$ $\Leftrightarrow$ $\kappa<0$ we have $\smax(h(0)) = \infty$ and therefore we have cuspless trajectory of infinite SR length, recall (\ref{eq:tfroms}). This fact is also clear from phase portrait of dynamic of vertical part~(\ref{vertpartchi}), see Fig.~\ref{fig:IntegralsVP} (bottom, left), where intersection of the integral curve with the red straight line indicates the moment, when cusp appears.
We see that some geodesics have no cusps in their spherical projection up to infinity. Notice that since $\SO$ is compact, it has bounded diameter (i.e. there exists $D>0$, such that the SR distance between any two elements of $\SO$ does not exceed $D$). Thus in contrast to $\SE$ group, where every cuspless geodesic is optimal (see~\cite{DuitsJMIV2014}), we observe that there exist nonoptimal cuspless geodesics in $\SO$ group.

In the hyperbolic case $\chi = a = \sqrt{\xi^2 -1}$, $\xi>1$ we have $\kappa \geq 0$ and the minimal positive root of (\ref{smaxeq}) reads as
\begin{eqnarray*}
\label{eq:smaxH}
\Scale[1.1]{
\smaxH(a, h_2^0, h_3^0) = \frac{1}{a}\log \left(-\frac{\sqrt{a^2 \left(1-(h_2^0)^2\right)+(h_3^0)^2}+a}{\left|a \, h_2^0+h_3^0\right|}\right).
}
\end{eqnarray*}
Finally, using the parameter $\chi = \sqrt{\xi^2-1}$, we have a single formula~(\ref{eq:smax}) for all cases.
$\hfill \Box$
\section{Outline of the Proof of Theorem~\ref{th:PDE}}\label{app:D}

By the result in \cite[Thm 11.15]{AgrachevBarilariBoscain} points $g \in \SO$ where the 
SR distance $d_{\SO}^{0}(e,\cdot)$ given by (\ref{SRdist}) is nonsmooth are either conjugate points, Maxwell points or abnormal points. So by assumption (and absence of abnormal geodesics due to a 2-bracket generating distribution $\Delta$, cf.~\cite[ch:20.5.1]{agrachev_sachkov}) we
see that $d_{\SO}^{0}(e,\cdot)$ is smooth at $\{\gamma_0^{*}(\tau)\;|\; 0< \tau \leq 1\}$. As a result, the mapping
\[
(0,1] \ni \tau \mapsto d_{\SO}^{0}(e,\gamma^*_{0}(\tau)) \in \R^+
\]
is smooth as it is a smooth composition of maps.
Similarly for the Riemannian case $\varepsilon>0$, recall (\ref{Rdist}), we have by the assumptions that
\[
(0,1] \ni \tau \mapsto d_{\SO}^{\varepsilon}(e,\gamma^*_{\varepsilon}(\tau)) \in \R^+
\]
is smooth for all $\varepsilon>0$.

Let $\mathcal{G}^{-1}$ denote the inverse metric tensor associated with $\mathcal{G}$ given by (\ref{SRMT}), and
let $\mathcal{G}_{\varepsilon}^{-1}$ denote the inverse metric tensor associated with $\mathcal{G}_{\varepsilon}$ given by (\ref{RMT}).

Now we rely on standard results \cite{Lions,Evans} on backtracking of optimal \emph{Riemannian} geodesics in Riemannian manifolds. This means that we find the smooth optimal geodesics via an intrinsic gradient descent on the
Riemannian distance map $\mathcal{W}_{\varepsilon}(g):=d_{\SO}^{\varepsilon}(g,e)$. This yields the identity
\begin{equation}\label{BTR}
\dot{\gamma}^*_{\varepsilon}(\tau)= \mathcal{W}_{\varepsilon}(g)\; \mathcal{G}^{-1}_{\varepsilon} {\rm d}\mathcal{W}_{\varepsilon}(\gamma^*_{\varepsilon}(\tau)),
\end{equation}
where $\mathcal{W}_{\varepsilon}$ is the unique viscosity solution \cite{Evans,Lions} of the corresponding Riemannian eikonal equation given by
\begin{equation} \label{Reik}
\begin{array}{l}
\left.\mathcal{G}_{\varepsilon}\right|_{g}\left((\left.\mathcal{G}_{\varepsilon}\right|_{g})^{-1}{\rm d} \mathcal{W}_{\varepsilon}(g) , (\left.\mathcal{G}_{\varepsilon}\right|_{g})^{-1}{\rm d} \mathcal{W}_{\varepsilon}(g)\right)=1,  \\
\textrm{ for } g\neq e, \textrm{ and }\mathcal{W}_{\varepsilon}(e)=0.
\end{array}
\end{equation}
Next we transfer these standard results toward the sub-Riemannian case, by a limiting procedure.
Firstly, we note that following a fully tangential approach to \cite[App.A,Thm.2]{DuitsMeestersMirebeauPortegies} yields (see Remark~\ref{rem:below} below) the following respectively pointwise and uniform convergence:
\begin{equation} \label{help1}
\lim \limits_{\varepsilon \downarrow 0} \mathcal{W}_{\varepsilon}(g)= \mathcal{W}_{0}(g), \textrm{ and }
\gamma_{\varepsilon}^* \to \gamma_{0}^* \ , \textrm{ as }\varepsilon \downarrow 0.
\end{equation}
Secondly, one has for any smooth function $f:\SO \to \R$ that the intrinsic gradients converge:
\begin{equation}\label{conv}
\begin{array}{l}
\lim \limits_{\varepsilon \downarrow 0}
\mathcal{G}_{\varepsilon}^{-1} {\rm d}f(g)=   \mathcal{G}^{-1} {\rm d}f(g) \desda \\
\lim \limits_{\varepsilon \downarrow 0}\left(
\xi^{-2}\!\left.X_{1}\right|_{g}f \left.X_{1}\right|_{g}\! +\! \left.X_{2}\right|_g f\left.X_{2}\right|_{g}\! + \!
\frac{\varepsilon^2}{\xi^2} \left.X_{3}\right|_{g}f \left.X_{3}\right|_{g}\right)\\
=\xi^{-2}\!\left.X_{1}\right|_{g}f \left.X_{1}\right|_{g}\! +\! \left.X_{2}\right|_g f\left.X_{2}\right|_{g}, \textrm{ for all }g \in \SO.
 \end{array}
\end{equation}
Thirdly, as $\gamma_{\varepsilon}^*$ and $\gamma_{0}^*$ are solutions to the Hamiltonian system of the Pontryagin Maximum Principle, the trajectories are continuously depending on $\varepsilon>0$, and so are their derivatives.
As a result, we have
\begin{equation}\label{finalconv}
\begin{array}{ll}
\dot{\gamma}_0^*(\tau) &= \lim \limits_{\varepsilon \downarrow 0}\dot{\gamma}_\varepsilon^*(\tau) \overset{(\ref{BTR})}{=}
\lim \limits_{\varepsilon \downarrow 0}
\mathcal{W}_{\varepsilon}(g) \, (\mathcal{G}_{\varepsilon}^{-1}{\rm d} \mathcal{W}_{\varepsilon})(\gamma_\varepsilon^*(\tau)) \\
&\overset{(\ref{help1},\ref{conv})}{=} \mathcal{W}(g) \left(\mathcal{G}^{-1}{\rm d} \lim \limits_{\varepsilon \downarrow 0} \mathcal{W}_{\varepsilon}\right)
\left(\lim \limits_{\varepsilon \downarrow 0} \gamma^*_{\varepsilon}(\tau) \right) \\
 &\overset{(\ref{help1})}{=}  \mathcal{W}(g) \,\mathcal{G}^{-1} {\rm d}\mathcal{W}(\gamma^*_{0}(\tau)).
 \end{array}
\end{equation}
Furthermore, from (\ref{Reik}) it follows that
\[
\begin{array}{l}
1=
\lim \limits_{\varepsilon \downarrow 0} \sqrt{\left.\mathcal{G}_{\varepsilon}\right|_{g}\left(\left.\mathcal{G}_{\varepsilon}\right|_{g}^{-1} {\rm d}\mathcal{W}_{\varepsilon}(g),
\left.\mathcal{G}_{\varepsilon}\right|_{g}^{-1} {\rm d}\mathcal{W}_{\varepsilon}(g)\right)} \\
=
\sqrt{\left.\mathcal{G}\right|_{g}\left(\left.\mathcal{G}\right|_{g}^{-1} {\rm d}\mathcal{W}(g),
\left.\mathcal{G}\right|_{g}^{-1} {\rm d}\mathcal{W}(g)\right)}
\end{array}
\]
where we note that the third term under the square root, after substitution of (\ref{conv}), vanishes. 
Finally, the viscosity solution property \cite{Evans,Lions} is also naturally carried over due to continuous $\varepsilon$ dependence
of the Hamiltonian and the convergences (\ref{help1}), (\ref{finalconv}).
$\hfill \Box$
\begin{remark}\label{rem:below}
The idea for the convergence (\ref{help1}) which is tangential to the $\SE=\R^{2} \rtimes S^{1}$ case considered in \\ \cite[App.A,Thm.2]{DuitsMeestersMirebeauPortegies} is based on a similar convergence result by Jean-Marie Mirebeau \& Chen Da of a related elastica model \cite[app.A]{CHENDAPhDthesis}. It relies on closeness of controllable paths \cite[Cor.A.5]{CHENDAPhDthesis} and Arzela-Ascoli's theorem. Another ingredient is the continuity of the mapping of the pair $(\varepsilon,g)$ onto the corresponding indicatrix $\mathcal{B}_{\varepsilon}(g):=\{v \in T_{g}(\SO)\;|\; \left.\mathcal{G}\right|_{g}(v,v)\leq 1\}$, i.e.,
the continuity of
\[
[0,\epsilon_0] \times \SO \ni (\varepsilon,g) \mapsto \mathcal{B}_{\varepsilon}(g) \in \mathcal{P}(T_{g}(\SO)),
\]
where the powerset $\mathcal{P}(T_{g}(\SO))$ of each tangent space\\
$T_{g}(\SO) \equiv \R^{3}$ is equipped with the metric topology induced by the Hausdorff distance. The proof of this continuity is straightforward
and is therefore omitted here.
\end{remark}
\subsection{The SR-eikonal and Backtracking Equations for Cuspless SR Geodesics}\label{app:D1}
Now in this article, we are primarily interested in the SR geodesics whose spherical projection does not have cusps.
This can be taken into account by modifying the standard SR-eikonal equation~(\ref{eq:eikonalsystfull}) as
\begin{equation}\label{eq:eikonalModified}
\begin{cases}
\sqrt{\frac{\max(0,X_{1}|_g(\mathcal{W}))^2}{\xi^2} + \left(X_2|_g(\mathcal{W})\right)^2} = C(g), \text{ for } g\neq e, \\
\mathcal{W}(e) = 0,
\end{cases} \!\!\!\!\!
\end{equation}
and backtracking equation~(\ref{eq:backtrackmain}) as
\begin{equation}\label{eq:backtrackingModified}
\Scale[1]
{
\begin{cases}
\dot{\gamma}_b(t) =  - u_1(t) X_1|_{\gamma_b(t)}  - u_2(t) X_2|_{\gamma_b(t)},  \\
\gamma_b(0) = g_1,
\end{cases}
}
\end{equation}
where $u_1(t) = \frac{\max(0,X_1|_{\gamma_b(t)}(\mathcal{W}))}{\xi^{2} C^2(\gamma_b(t))} $, $u_2(t)= \frac{X_2|_{\gamma_b(t)}(\mathcal{W})}{C^2(\gamma_b(t))}$, and where the system is integrated for $t\in[0,\mathcal{W}(g_1)]$.

The idea behind this is that for cuspless SR geodesics one has $u_1$ positive which holds by if the corresponding momentum component
$h_{1}= X_{1}(\mathcal{W})$ is positive.
Note that in the eikonal equation one substitutes momentum $h={\rm d}\mathcal{W}$ into the Hamiltonian to achieve equidistant wavefront propagation \cite{Evans}.
For further details and analysis on the positive control restriction on a similar problem on $\SE$, see \cite{DuitsMeestersMirebeauPortegies}.
Similar analysis benefits may be expected on the $\SO$ case, but this is beyond the scope of this article.
\section{Comparison of the $\SE$ and $\SO$ Geodesics in the Image Plane}\label{app:geodcomp}
We organize the comparison in the following way:

1) Choose an initial condition $V_0 = (X_0,Y_0,\Theta_0)$ and a terminal  condition $V_1 = (X_1,Y_1,\Theta_1)$ in
$$\Scale[0.93]{D = [X_{min},X_{max}] \times [Y_{min},Y_{max}] \times [-\pi,\pi] \subset \R^2 \times S^1 \cong \SE,}$$
where $X_{max} = Y_{max} = -X_{min} = -Y_{min} =$ $(c + a) \tan\psi_{max}$ $= \frac{149}{105} \tan\frac{\pi}{8}$. These values are obtained from schematic eye model for $\eta=1$, recall Sect.~\ref{subsec:Motivation1}, where $\psi_{max} = \frac{\pi}{8}$ is a maximum scanning angle via a standard fundus cameras.

2) Compute the SR distance $\Scale[0.9]{\mathcal{W}^{\SE}(V) = d^{\SE}(V_0, V)}$ in the volume of interest $D \ni V$ via SR-FM-SE(2) (sub-Riemannian fast marching in $\SE$, see~\cite{CIARP}).

3) Find via backtracking in $\SE$ (see~\cite{CIARP}) a geodesic $\gamma^{\SE}(\cdot)$ satisfying $\gamma^{\SE}(0) = V_0$ and $\gamma^{\SE}(T) = V_1$, with $T = \mathcal{W}^{\SE}(V_1)$.

4) Find the initial condition $\nu_0 = (x_0,y_0,\th_0)$ and the terminal condition $\nu_1 = (x_1,y_1,\th_1)$ in $\SO$, obtained from $V_0$ and $V_1$ via planar projection $\Pi^{-1}$, recall Eq.~(\ref{eq:XYtoxy}), as $(x_i, y_i) = \Pi^{-1}(X_i,Y_i) =: (\Phi_1(X,Y), \Phi_2(X,Y))$ and
$$ \th_i = \arg\left(\dot{\gothic{x}}(X_i,Y_i,\Theta_i)+ \i \cos \Phi_1(X_i,Y_i) \dot{\gothic{y}}(X_i,Y_i,\Theta_i)\right),$$
where \\
$\left. \begin{array}{c}
\displaystyle \dot{\gothic{x}}(X_i,Y_i,\Theta_i) := \frac{\partial \Phi_1}{\partial X}|_{(X_i,Y_i)} \cos\Theta_i + \frac{\partial \Phi_1}{\partial Y}|_{(X_i,Y_i)} \sin\Theta_i, \\[7pt]
\displaystyle \dot{\gothic{y}}(X_i,Y_i,\Theta_i) := \frac{\partial \Phi_2}{\partial X}|_{(X_i,Y_i)} \cos\Theta_i + \frac{\partial \Phi_2}{\partial Y}|_{(X_i,Y_i)} \sin\Theta_i.
\end{array} \right.$

5) Compute the SR distance $\mathcal{W}^{\SO}(\nu) = d^{\SO}(\nu_0, \nu)$ in the domain $\nu \in \mathcal{D} = [x_{min},x_{max}] \times [y_{min},y_{max}] \times [-\pi,\pi]$ via SR-FM in $\SO$. Here we set $x_{max} = y_{max} = -x_{min} = -y_{min} = 0.63$, recall~Eq.~\!(\ref{eq:ymax}).

6) Find via backtracking~(\ref{eq:backtrackmain}) a geodesic $\gamma^{\SO}(\cdot)$ satisfying $\gamma^{\SO}(0) = \nu_0$ and $\gamma^{\SO}(T) = \nu_1$.

7) Plot in the image plane both the spatial projection $\Gamma^{\SE}(\cdot) = (X(\cdot),Y(\cdot))$ of the geodesic $\gamma^{\SE}(\cdot)$ and  the planar projection $\Gamma^{\SO}(\cdot)=\Pi(x(\cdot), y(\cdot))$  of the geodesic $\gamma^{\SO}(\cdot)$.


\begin{thebibliography}{References}
\bibitem{PeyreCohen2010}
G.~Peyr\'{e}, M.~P\'{e}chaud, R.~Keriven, L.D.~Cohen, {\em Geodesic methods in computer vision and graphics.} Found. Trends. Comp. in Comput. Graph. Vis., 5(34), pp.197--397, 2010.

\bibitem{sethian}
J.A.~Sethian, {\em Level Set Methods and Fast Marching Methods.}
Cambridge University Press, Cambridge, 1999.

\bibitem{Kimmel}
V.~Caselles, R.~Kimmel, G.~Sapiro,
\newblock {\em Geodesic active contours},
\newblock Int. J. Comp. Vis. 22(1), p.~61--79, 1997.

\bibitem{osher}
S.~Osher, R.P.~Fedkiw, {\em Level Set Methods and Dynamic Implicit Surfaces. Applied Mathematical Science.} Applied mathematical science, Springer, New York, 2003.

\bibitem{mirebeau}
J-M.~Mirebeau,{\em Anisotropic fast-marching on cartesian grids using Lattice Basis Reduction.}
SIAM J. Num. Anal. pp.1573 ,52(4), 2014.

\bibitem{Ikram}
M.K. Ikram, Y.T. Ong, C.Y. Cheung, and T.Y. Wong. Retinal vascular caliber
measurements: clinical significance, current knowledge and future perspectives.
Ophthalmologica, 229(3):125--136, 2013.

\bibitem{Lions}
P.L.~Lions,
{\em Generalized Solutions of Hamilton-Jacobi Equations.}
Research Notes In Mathematics Series, 69, 1982.

\bibitem{Crandall}
M.G.~Crandall, P.L.~Lions, Viscosity solutions of Hamilton-Jacobi equations. Trans. Am. Math. Soc. 277 (1983), no. 1, 1–-42.

\bibitem{Sasongko}
M.B.~Sasongko, T.Y.~Wong, T.T.~Nguyen, C.Y.~Cheung, J.E.~Shaw, J.J.~Wang, Retinal vascular tortuosity in persons with diabetes and diabetic retinopathy. Diabetologia 2011; 54: pp. 2409--2416

\bibitem{Kalitzeos}
Kalitzeos, A.A., Lip, G.Y., Heitmar, R.: Retinal vessel tortuosity measures and their applications. Exp. Eye Res., 106, 40 -- 46, 2013.

\bibitem{Erik_Curv}
Bekkers, E., Zhang, J., Duits, R., ter Haar Romeny, B.:
Curvature based biomarkers for diabetic retinopathy via exponential curve fits in SE(2).
In Trucco, E., Chen, X., Garvin, M.K., Liu, J.J., Frank, X.Y., eds.: Proceedings of the Ophthalmic Medical Image Analysis Second International Workshop, OMIA 2015, Held in Conjunction with MICCAI 2015, Munchen, Germany, October 9, 2015. Iowa Research Online (2015) 113-120.

\bibitem{Evans}
L.C.~Evans, {\em Partial Differential Equations.} Graduate Studies in Mathematics Vol. 19, AMS, Providence (USA), 1998.

\bibitem{Petitot_B}
 J.Petitot, The neurogeometry of pinwheels as a sub-Riemannian contact structure, J. Physiology - Paris, 97 (2003), 265--309.

\bibitem{Citti_Sarti_B}
G. Citti and A. Sarti, A cortical based model of perceptual completion in the roto-translation
space, J. Math. Imaging Vis., 24 (2006), 307--326.

\bibitem{Boscain_Rossi_B}
 U. Boscain, F. Rossi, Invariant Carnot-Caratheodory metrics on $S^3$, $SO(3)$, $SL(2)$ and Lens Spaces, SIAM, Journal on Control and Optimization, Vol 47 (2008), 1851--1878.

\bibitem{Boscain_Rossi_S2}
U. Boscain, F. Rossi, Projective Reed-Shepp Car on $S^2$ with quadratic cost, ESAIM: COCV 16 (2010), pp. 275--297.



\bibitem{Gauthier}
U. Boscain, J.-P. Gauthier, R. Chertovskih, A. Remizov
Hypoelliptic diffusion and human vision: a semidiscrete new twist
SIAM J. Imaging Sci., 7 (2) (2014), pp. 669-695.

\bibitem{Remco_Eric_B}
 E. Bekkers, R. Duits, T. Berendschot, B.H. Romeny, A multi-orientation analysis approach to retinal vessel tracking. J. Math. Imaging Vis. (2014)

\bibitem{Ivan}
I. Beschastnyi, Y. Sachkov, Sub-Riemannian and almost-Riemannian geodesics on $\SO$ and $S^{2}$, ArXiv, 1409.1559, 2014.

\bibitem{Berestovskiy1}
V. Berestovskii, Geodesics of a left-invariant nonholonomic Riemannian metric on the group of motions of the euclidean plane, Siberian Mathematical Journal, 35 (6), 1994, pp. 1083-1088.

\bibitem{Berestovskiy2}
V.Berestovskii, I. Zubareva, Sub-Riemannian distance in the Lie groups $SU(2)$ and $SO(3)$, Sib. Adv. Math., 26(2), 2016, pp. 77-89.

\bibitem{Schcerbakova}
B. Bonnard, O. Cots, J.-B. Pomet, N. Shcherbakova, Riemannian metrics on 2D-manifolds related to the Euler-Poinsot rigid body motion, ESAIM: Control Optim. Calc. Var., 20 (2014), pp 864-893.

\bibitem{Bonnard}
B. Bonnard, M. Chyba. Two applications of geometric optimal control to the
dynamics of spin particle, 2014, Preprint HAL Id: hal-00956828.

\bibitem{Jurd}
V. Jurdjevic,
Geometric Control Theory,
Cambridge Studies in Advanced Mathematics, (No. 52),
Cambridge: Cambridge University Press,
1996.

\bibitem{Kruglikov}
B. Kruglikov, A. Vollmer, G. Lukes-Gerakopoulos, On integrability of certain rank 2 sub-Riemannian structures,  ArXiv, 1507.03082, 2015.

\bibitem{NMTMA_B}
A.P. Mashtakov, A.A. Ardentov, Y.L. Sachkov, Parallel algorithm and software for image inpainting via sub-Riemannian minimizers on the group of rototranslations, Numer. Math. Theory Methods Appl., Vol. 6, No. 1. (2013), 95-115.

\bibitem{Mashtakov_Sachkov_Integr}
A.~P.~Mashtakov, Yu.~L.~Sachkov,
Superintegrability of left-invariant sub-Riemannian structures on unimodular three-dimensional Lie droups.
Differ. Equ., 2015, Vol. 51, No. 11, pp. 1476--1483.

\bibitem{cut_sre2}
Yu.L. Sachkov,
Cut locus and optimal synthesis in the sub-Riemannian problem  on the group of motions of a
plane, {\em ESAIM: COCV}, 17 (2011), 293--321.

\bibitem{max_sre}
I. Moiseev, Yu.L. Sachkov,
Maxwell strata in sub-Riemannian problem  on the group of motions of a plane,
{\em ESAIM: COCV}, 16 (2010), 380--399.

\bibitem{Boscain2014}
U.~Boscain, R.~Duits, F.~Rossi and Y.~L.~Sachkov,
Curve cuspless reconstruction via sub-Riemannian geometry,
ESAIM: Control Optim. Calc. Var.,
\emph{DOI: 10.1051/cocv/2013082}.
20 (2014), 748--770.


\bibitem{DuitsJMIV2014}
R.~Duits, U.~Boscain, F.~Rossi and Y.~L.~Sachkov,
Association fields via cuspless sub-Riemannian geodesics in SE(2).
JMIV,
\emph{DOI 10.1007/s10851-013-0475-y},
\emph{http://bmia.bmt.tue.nl/people/RDuits/cusp.pdf}
49 (2), (2014), 384--417.

\bibitem{SE3}
R.~Duits, A.~Ghosh, T.C.J.~Dela~Haije, A.~Mashtakov,
On sub-Riemannian geodesics in $SE(3)$ whose spatial projections do not have cusps,
J. Dyn. Control Syst., 22(4), pp.771-805, 2016.


\bibitem{Montgomery}
R.~Montgomery,
A Tour of Subriemannian Geometries, Their Geodesics and Applications. -- Mathematical Surveys and Monographs, 2002.

\bibitem{Agrachev_exp_map}
A.A.~Agrachev,
Exponential mappings for contact sub-Riemannian structures, J. Dyn. Control Syst., 2:3 (1996), 321--358.

\bibitem{AgrachevBarilariBoscain}
A.A.~Agrachev, D.~Barilari and U.~Boscain,
Introduction to Riemannian and Sub-Riemannian Geometry
from the Hamiltonian Viewpoint.
Preprint SISSA 09/2012/M, 2016.


\bibitem{DuitsMeestersMirebeauPortegies}
R.~Duits, S.~Meesters, J.M.~Mirebeau and J.Portegies,
Optimal paths for variants of the 2D and 3D Reeds-Shepp car with applications in image Analysis.
Preprint on Arxiv \emph{www.arxiv.org/pdf/1612.06137.pdf}, 2016.


\bibitem{agrachev_sachkov}
A.A.~Agrachev, Yu.L. Sachkov, Control Theory from the Geometric Viewpoint,
Springer-Verlag, 2004.

\bibitem{bookchapter}
R.~Duits, A.Ghosh, T.C.J.~Dela Haije, Y.L.~Sachkov,
Cuspless Sub-Riemannian Geodesics within the Euclidean Motion Group $SE(d)$, in Neuromathematics of Vision, Springer Series
Lecture Notes in Morphogenesis, vol.1 p.173--240, 2014.

\bibitem{SSVM}
Bekkers, E.J., Duits, R,  Mashtakov, A. and  Sanguinetti, G.R. (joint main authors).
Data-driven Sub-Riemannian Geodesics in $\SE$,
Proc. SSVM 2015 (eds. Aujol, Nikolova, Papadakis), LNCS 9087, p.613-625,  2015.

\bibitem{SIAM}
Bekkers E.J., Duits R., Mashtakov A., Sanguinetti, G.R. (joint main authors).
A PDE Approach to Data-driven Sub-Riemannian Geodesics in $\SE$.
SIAM Journal on Imaging Sciences, 2015, 8:4, PP 2740-2770.

\bibitem{CIARP}
G. Sanguinetti, R. Duits, E. Bekkers, M.H.J. Janssen, A. Mashtakov, J.M. Mirebeau.
Sub-Riemannian Fast Marching in SE(2).
 In book: Progress in Pattern Recognition, Image Analysis, Computer Vision, and Applications. Lecture Notes in Computer Science, Volume 9423, 2015, pp 366-374.






\bibitem{ClinOpt}
American Academy of Ophthalmology, Clinical Optics, 2013-2014. Basic and clinical science course. 2013, San Francisco, CA: American Academy of Ophthalmology.

\bibitem{Rouy}
E.~Rouy and A.~Tourin, {\em A viscosity solutions approach to shape-from-shading}: SIAM J. Num. Anal, 29, 867--884, 1992.

\bibitem{Markina1}
O. Calin, D.-C. Chang, I. Markina,
{\em Sub-Riemannian geometry on the sphere $S^3$}: Canad. J. Math., 2009, Vol. 61(4), pp. 721-739.

\bibitem{Markina2}
D.-C. Chang, I. Markina, A. Vasilev,
{\em Sub-Riemannian geometry on the 3-D sphere}:
 Complex analysis and operator theory, 2009, issue 3, pp. 361-377.

\bibitem{AssocField}
D.J.~Field, A.~Hayes, R.~Hess, Contour integration by the human visual system: Evidence for a local ``association field'', Vision Research, V. 33, 2, 1993, pp. 173-193.

\bibitem{Merabeau_Finsler}
D.~Chen, J.-M.~Mirebeau, L.D.~Cohen. Global Minimum for Curvature Penalized Minimal Path Method. In Xianghua Xie, Mark W. Jones, and Gary K. L. Tam, editors, Proceedings of the British Machine Vision Conference (BMVC), pages 86.1-86.12. BMVA Press, September 2015.

\bibitem{CHENDAPhDthesis}
D.~Chen, New Minimal Path Models for Tubular Structure Extraction and Image Segmentation.
PHD thesis, \'{E}cole Doctorale de Dauphine Paris, 2016.

\bibitem{Rifford}
L.~Rifford, {\em Sub-Riemannian geometry and optimal transport}. Springer Briefs in Mathematics. Springer International Publishing, New York (2014).

\bibitem{Frangi}
A. Frangi, W. Niessen, K. Vincken, and M. Viergever, {\em Multiscale vessel enhancement filtering}, In Proc. of the MICCAI, Lecture Notes in Computer Science, Cambridge, pp. 130-137, 1998.

\bibitem{BartBook}
B.~M. ter Haar Romeny, {\em Front-End Vision and Multi-Scale Image Analysis}, Springer, 484 pp, 2004.

\end{thebibliography}
\end{document}